\documentclass[12pt]{amsart}

\usepackage{longtable}
\usepackage{float}
\restylefloat{table}

\newcommand{\hide}[1]{}

\usepackage{amssymb}
\usepackage{amsbsy}
\usepackage{amscd}
\usepackage{amsmath}
\usepackage{amsthm}
\usepackage{bbold}
\usepackage{mathrsfs}
\usepackage{url}
\usepackage{graphicx}

\numberwithin{equation}{subsection}

\theoremstyle{plain}
\newtheorem{thm}{Theorem}[subsection]
\newtheorem{prop}[thm]{Proposition}

\newtheorem{thm-defi}[thm]{Theorem/Definition}
\newtheorem{cor}[thm]{Corollary}
\newtheorem{lem}[thm]{Lemma}

\theoremstyle{definition}

\newtheorem{defi}[thm]{Definition}
\newtheorem{rem}[thm]{Remark}

\newtheorem{question}[thm]{Question}
\newtheorem{example}[thm]{Example}

\newcommand{\A}{{\mathcal A}}
\newcommand{\B}{{\mathcal B}}

\newcommand{\CC}{{\mathbb C}}

\newcommand{\QQ}{{\mathbb Q}}
\newcommand{\R}{{\mathcal R}}
\newcommand{\RR}{{\mathbb R}}

\newcommand{\T}{{\mathcal T}}

\newcommand{\ZZ}{{\mathbb Z}}
\renewcommand{\P}{{\mathcal P}}
\newcommand{\PP}{{\mathbb P}}
\newcommand{\ann}{{\rm ann}}

\newcommand{\Integers}{{\mathbb Z}}

\newcommand{\Spin}{{\rm Spin}}

\newcommand{\fm}{{\mathfrak m}}

\newcommand{\IsomRightArrow}{\stackrel{\cong}{\rightarrow}}

\newcommand{\RightArrowOf}[1]{\stackrel{#1}{\rightarrow}}

\newcommand{\LongRightArrowOf}[1]{\stackrel{#1}{\longrightarrow}}

\newcommand{\StructureSheaf}[1]{{\mathcal O}_{#1}}
\newcommand{\EndProof}{\hfill  $\Box$}
\newcommand{\restricted}[2]{#1_{\mid_{#2}}}

\newcommand{\disc}{{\rm disc}}

\newcommand{\rank}{{\rm rank}}

\newcommand{\Sym}{{\rm Sym}}
\newcommand{\Ext}{{\rm Ext}}

\newcommand{\Hom}{{\rm Hom}}

\newcommand{\End}{{\rm End}}

\newcommand{\Contract}{\rfloor}

\newcommand{\Choose}[2]
{\left(\!\!\begin{array}{c}#1\\#2\end{array}\!\!\right)}

\renewcommand{\span}{{\rm span}}

\input xy
\xyoption{all}

\begin{document}
\title[Secant sheaves and Weil classes on abelian $2n$-folds of CM-type]
{Secant sheaves on abelian $n$-folds with real multiplication and Weil classes on abelian $2n$-folds with complex multiplication}
\author{Eyal Markman}
\address{Department of Mathematics and Statistics, 
University of Massachusetts, Amherst, MA 01003, USA}
\email{markman@umass.edu}

\date{\today}
%\rightline{\today}
%\begin{center}
%\begin{Large}
%{\bf 
%\noindent
%The monodromy of generalized Kummer varieties}
%\end{Large}
%\\
%Eyal Markman
%\end{center}

\begin{abstract}
Let $K$ be a CM-field, i.e., a totally complex quadratic extension of a totally real field $F$. Let $X$ be an abelian variety admitting an algebra embedding $\hat{\eta}:F\rightarrow \End_\QQ(X)$, and let 
$\hat{X}$ be the dual abelian variety. We construct an embedding $\eta:K\rightarrow \End_\QQ(X\times\hat{X})$ associated to a choice of a polarization $\Theta$ in $\wedge^2_FH^1(X,\QQ)$ and an element $q\in F$, such that $K=F(\sqrt{-q})$. 
We get the $[K:\QQ]$-dimensional subspace $HW(X\times\hat{X},\eta)$ 
of Hodge Weil classes in $H^d(X\times\hat{X},\QQ)$, where $d:=4\dim(X)/[K:\QQ]$. 

The even cohomology $S^+_\CC:=H^{ev}(X,\CC)$ is the half-spin representation of the group $\Spin(H^1(X\times\hat{X},\CC))$ and so $\PP(S^+_\CC)$ contains the even spinorial variety. The even spinorial variety is a component of the Grassmannian of maximal isotropic subspaces of $H^1(X\times\hat{X},\CC)$. We associate to $(\Theta,q)$ a rational $2^{[F:\QQ]}$-dimensional subspace $B$ of $S^+_\QQ$ such that $\PP(B)$ is secant to the even spinorial variety. Associated to two coherent sheaves $F_1$ and $F_2$ on $X$ with Chern characters in $B$ we obtain the object $E:=\Phi(F_1\boxtimes F_2^\vee)$ in the derived category $D^b(X\times\hat{X})$, where $\Phi:D^b(X\times X)\rightarrow D^b(X\times\hat{X})$ is Orlov's equivalence.
The flat deformations of the normalized Chern class $\kappa(E):=ch(E)\exp\left(-\frac{c_1(E)}{\rank(E)}\right)$ of $E$ remain of Hodge type under every deformation of $(X\times\hat{X},\eta)$ as an abelian variety $(A',\eta')$ of Weil type. 

We provide a criterion for $ch(F_1)\otimes ch(F_2)$ to belong to the open subset in $B\otimes B$ for which 
the algebraicity of the flat deformation of
$\kappa(E)$ implies the algebraicity of all classes in $HW(A',\eta')$. The algebraicity would thus follow if E is semiregular in the appropriate sense. 
Examples of such $F_1$ and $F_2$ are provided for $X$ the Jacobian with real multiplication by $F$ of a genus $4$ curve when $[F:\QQ]=2$ and $Gal(K/\QQ)\cong \ZZ/2\ZZ\times\ZZ/2\ZZ$, but the search for semiregular such $E$ is left for future work.
\end{abstract}

\maketitle

\setcounter{tocdepth}{3}
\tableofcontents

%****************************************************************
% 
%****************************************************************
\section{Introduction}
\label{sec-introduction}

%****************************************************************
% 
%****************************************************************
\subsection{A summary of the construction}
\label{sec-summary-of-the-construction}
Let $K$ be a CM-field and $F$ its totally real subfield.
The Galois involution $\iota$  in $Gal(K/F)$ is a central element of $Gal(K/\QQ)$ and for every embedding $\sigma:K\rightarrow \CC$ we have $\sigma(\iota(t))=\overline{\sigma(t)}$, for all $t\in K$, where the right hand side is complex conjugation \cite[Sec. 4]{deligne-milne}. 
Set 
\[
e:=\dim_\QQ(K).
\] 
Let $\Sigma:=\Hom(K,\CC)$ be the set of complex embeddings. Note that the cardinality of $\Sigma$ is $e$.
The group $Gal(K/\QQ)$ acts on $\Sigma$ by $g(\sigma)=\sigma\circ g^{-1}$, for all $g\in Gal(K/\QQ)$ and $\sigma\in\Sigma$. The latter action is transitive if $K$ is a Galois extension of $\QQ$.

Let $A$ be an abelian variety and set $\End_\QQ(A):=\End(A)\otimes_\ZZ\QQ$.
Given an embedding
$\eta:K\rightarrow \End_\QQ(A)$, we get a decomposition
\[
H^1(A,\CC)=\oplus_{\sigma\in\Sigma}H^1_\sigma(A),
\]
where $\eta(K)$ acts on $H^1_\sigma(A)$ via the embedding $\sigma$. The complex multiplication $\eta$ corresponds to an embedding $\eta:K\rightarrow \End(H^1(A,\QQ))$ preserving the Hodge structure.
Set 
\[
d:=\dim_K(H^1(A,\QQ)).
\] 
We have  $\dim_\CC(H^1_\sigma(A))=d$, for all $\sigma\in \Sigma$. 
Set $H^{1,0}_\sigma(A):=H^1_\sigma(A)\cap H^{1,0}(A)$ and
$H^{0,1}_\sigma(A):=H^1_\sigma(A)\cap H^{0,1}(A)$. Assume that for every $\sigma\in\Sigma$, we have 
\begin{equation}
\label{eq-condition-for-HW-to-consist-of-Hodge-classes}
\dim(H^{1,0}_\sigma(A))=\dim(H^{0,1}_\sigma(A))=\frac{d}{2}.
\end{equation}
Then $\wedge_K^dH^1(A,\QQ)$ is an $e$-dimensional subspace $HW(A,\eta)$ of $H^{\frac{d}{2},\frac{d}{2}}(A,\QQ)$, which is a one-dimensional $K$ vector space
\cite[Prop. 4.4]{deligne-milne}.
The challenge is then to prove the algebraicity of the classes in $HW(A,\eta)$. 
%$\wedge_K^dH^1(A,\QQ)$. 
It suffices to prove the algebraicity of one non-zero class in $HW(A,\eta)$,
%$\wedge_K^dH^1(A,\QQ)$, 
as $K$ acts via algebraic correspondences. If $d=2$, then $HW(A,\eta)$ consists of rational $(1,1)$ classes, which are algebraic.
Hence, we may assume that $d>2$. 

A pair $(A,\eta)$ satisfying (\ref{eq-condition-for-HW-to-consist-of-Hodge-classes}) is called an abelian variety of {\em Weil type} and the
classes in $HW(A,\eta)$ are called {\em Weil classes} in the literature (see \cite{Moonen-Zarhin-Weil-classes}).
When $F=\QQ$, so that $K$ is an imaginary quadratic number field, the Weil classes were introduced in \cite{weil}. They were proved algebraic for all abelian fourfolds of Weil type and for abelian sixfolds of split Weil type when $K=\QQ(\sqrt{-3})$ in \cite{schoen1,schoen}, when $K=\QQ(\sqrt{-1})$ in \cite{van-Geemen,koike},
and for all imaginary quadratic number fields in \cite{markman-sixfolds}. Proofs of the algebraicity of the Weil classes on abelian fourfolds of split Weil type for all imaginary quadratic number fields were obtained using hyper-K\"{a}hler techniques in \cite{FF,van-Geemen-rapagneta,markman-generalized-kummers}.
The current paper presents the natural generalization of the strategy of \cite{markman-sixfolds} for more general CM-fields.

%A polarized abelian variety of Weil type consists of a triple $(A,\eta,h)$, where $(A,\eta:K\rightarrow \End_\QQ(A))$ is an abelian variety of Weil type 
%and $h\in H^{1,1}(A,\QQ)$ is the class of a polarization, such that for $k\in K$, $[\eta(k)]^*(h)=Nn_{K/}$ (???)

Let $d$ be a positive even integer, let $X$ be an abelian variety of dimension 
\[
n:=de/4,
\] 
and set $A:=X\times\hat{X}$. 
We will construct the embedding $\eta:K\rightarrow \End_\QQ(A)$ 
%depending on a choice of a CM-type $T$,  
assuming $X$ admits an embedding $\hat{\eta}:F\rightarrow \End_\QQ(X)$ of the totally real subfield $F$ of $K$.

The free abelian group $V:=H^1(X,\ZZ)\oplus H^1(\hat{X},\ZZ)$ admits a natural unimodular symmetric bilinear pairing. Using the natural isomorphism $H^1(\hat{X},\ZZ)\cong H^1(X,\ZZ)^*$, the pairing is given by 
\begin{equation}
\label{eq-pairing-on-V}
((w_1,\theta_1),(w_2,\theta_2))_V:=\theta_1(w_2)+\theta_2(w_1).
\end{equation}
The $K$ vector space 
$
V_{\hat{\eta},K}:=[H^1(X,\QQ)\oplus H^1(\hat{X},\QQ)]\otimes_FK
$
admits a natural $K$-valued non-degenerate symmetric bilinear form (\ref{eq-pairing-on-V-hat-eta-K}). Let $W$  be a maximal isotropic subspace of $V_{\hat{\eta},K}$. 
%image of $H^1(\hat{X},\QQ)\otimes_FK$ under an element 
%\begin{equation}
%\label{eq-g}
%g_0
%\end{equation} 
%of a $K$-linear spin group, defined below in (\ref{eq-Spin-V-K}), so that $W$ is maximal isotropic.
% and consider its orbit under $Gal(K/\QQ)$. As the subspace $H^1(X,\QQ)$ is $Gal(F/\QQ)$ invariant, the orbit will consist of two subspace $W_T$ and 
Let $\bar{W}$ be the image of $W$ via the non-trivial element $\iota$ in $Gal(K/F)$. 
Assume that $W\cap\bar{W}=(0)$. Then $W$ is naturally isomorphic to 
$H^1(A,\QQ)$, as a vector space over $\QQ$, and so determines a $K$-vector space structure on $H^1(A,\QQ)$ (Section \ref{sec-complex-multiplication--by-CM-field-K}). We construct such a $W$ associated to a choice of a polarization $\Theta\in [H^{1,1}(X,\QQ)\cap\wedge^2_FH^1(X,\QQ)]$ and a non-zero element $q\in K$, such that $\iota(q)=-q$,  and show that scalar multiplication by elements of $K$ preserves the Hodge structure of $H^1(A,\QQ)$, hence yields an embedding
$\eta:K\rightarrow \End_\QQ(A)$, such that $(A,\eta)$ is of Weil type
(Section \ref{sec-example-of-K-secant-for-CM-field})

The grassmannian $IGr^+(2n,V_\CC)$ of even maximal isotropic subspaces of $V_\CC:=V\otimes_\ZZ\CC$ naturally embeds in the projectivization $\PP(S^+_\CC)$ of the even half-spin representation $S^+_\CC:=H^{ev}(X,\CC)$ of $\Spin(V_\CC)$. We refer to the image of $IGr^+(2n,V_\CC)$ in $\PP(S^+_\CC)$ as the {\em even spinorial variety}.
The point $\ell_W\in \PP(S^+_\CC)$ associated to a an even maximal isotropic subspace $W$ of $V_\CC$ is called an {\em even pure spinor}, as are the classes in the corresponding line in $S^+_\CC$ \cite{chevalley}. 

Let $\hat{\Sigma}:=\Sigma/\iota$ be the set of embeddings of $F$ in $\RR$. 
A {\em complex multiplication type} (CM-type) $T:\hat{\Sigma}\rightarrow\Sigma$ consists of a choice of an embedding $T(\hat{\sigma})\in\Sigma$ restricting to $F$ as $\hat{\sigma}$, for each  $\hat{\sigma}\in\hat{\Sigma}$ (see \cite[Sec. 5]{deligne-milne}). 
A maximal isotropic subspace $W$ of $V_{\hat{\eta},K}$ determines a maximal isotropic subspace $W_T$ of $V_\CC:=H^1(A,\CC)$,
given in (\ref{eq-W-T}),
hence 
%we get an even maximal isotropic subspace $W_T$ of $H^1(A,\CC)$ and so 
a pure spinor $\ell_{W_T}\in \PP(S^+_\CC)$,
for each CM-type $T$.
We denote $\ell_{W_T}$ by $\ell_T$ for simplicity.
The subspace of $S^+_\CC$ spanned by the set $\{\ell_T\}$, as $T$ varies over all CM-types, is defined over $\QQ$ and corresponds to
a subspace $B$ of $S^+_\QQ$ of dimension $2^{e/2}$ (Lemma \ref{lemma-B-is-rational}). By definition, the linear subspace $\PP(B)$ is secant to the even spinorial variety.

A {\em $B$-secant sheaf} is a coherent sheaf $G$ on $X$ with $ch(G)$ in $B$. Given two $B$-secant sheaves $F_1$ and $F_2$ we get the object $E:=\Phi(F_1\boxtimes F_2^\vee)$ in $D^b(X\times\hat{X})$ via Orlov's derived equivalence $\Phi:D^b(X\times X)\rightarrow D^b(X\times\hat{X})$ 
given in (\ref{eq-Orlov-derived-equivalence}). 
Assume that $r:=\rank(E)\neq 0$. We prove that the characteristic class $\kappa(E):=ch(E)\exp(-c_1(E)/r)$ remains of Hodge type under every deformation of $(A,\eta)$, $\eta:K\rightarrow \End_\QQ(A)$, as an abelian variety of Weil type (Proposition \ref{prop-kappa-class-of-image-of-secant-class-yields-a-HW-class}(\ref{lemma-item-kappa-phi-c-is-invariant})). Let $\A^{2k}$ be the subspace of $H^{k,k}(A,\QQ)$ consisting of classes which remain of Hodge type under 
 every deformation of $(A,\eta)$ as an abelian variety of Weil type. The subspace $\A^d$ is the direct sum
 \begin{equation}
 \label{eq-decomposition-of-generic-Hodge-classes-in-H-d}
\A^d= Im(\Sym^{d/2}(\A^2))\oplus HW(A,\eta),
 \end{equation}
where  $Im(\Sym^{d/2}(\A^2))$ is the image of $\Sym^{d/2}(\A^2)$
in $H^d(A,\QQ)$, by Proposition \ref{prop-generator-for-the-invariant-subalgebra}.

We provide a criterion for the graded summand $\kappa_{d/2}(E)\in \A^d$ to project to a non-zero class in the summand $HW(A,\eta)$ of (\ref{eq-decomposition-of-generic-Hodge-classes-in-H-d}) in terms of the position of $ch(F_1)\otimes ch(F_2)$ with respect to 
a natural grading on $B\otimes B$ introduced  in the next paragraph.
%(Proposition \ref{prop-kappa-class-of-image-of-secant-class-yields-a-HW-class}(\ref{lemma-item-difference-projects-to-HW})).
A generic class in $B\otimes B$ satisfies the criterion, which is an open condition. 
A refinement of the latter criterion, in terms of a grading of $B$, is obtained in 
Lemma \ref{lemma-a-criterion-for-alpha-0-tensor-beta-to-map-via-orlov-equivalence-to-a-class-inducing-non-zero-HW-class}, 
when the Galois group $Gal(\tilde{K}/\QQ)$ of the Galois closure $\tilde{K}$ 
of $K$ splits as $\ZZ/2\ZZ\times Gal(\tilde{F}/\QQ)$, where $\tilde{F}$ is the Galois closure of $F$.
In Example \ref{example-coherent-sheaf-on-Jacobian-of-genus-4-with-ch-a-multiple-of-beta-prime} we exhibit two $B$-secant sheaves $F_1$ and $F_2$, satisfying the criterion for $\kappa_{d/2}(E)$ to project to a non-zero class in $HW(A,\eta)$,
when $X$ is the Jacobian of a genus $4$ curve with real multiplication $\hat{\eta}:F\rightarrow \End_\QQ(X)$ by a real quadratic number field $F$, and $K=F(\sqrt{-q})$,
where $q$ is a positive integer (Lemma \ref{lemma-box-tensor-product-of-sheaves-maps-to-one-with-equivariant-kappa-class-case-of-degree-4-extension}). 

We state next the criterion for the graded summand $\kappa_{d/2}(E)\in H^d(A,\QQ)$ to project to a non-zero class in $HW(A,\eta)$.
Given two CM-types $T$, $T'$, let $|T\cap T'|$ be the cardinality of the subset $T(\hat{\Sigma})\cap T'(\hat{\Sigma})$ of $\Sigma$.
Then $B\otimes B$ decomposes as a direct sum $\oplus_{k=0}^{e/2}BB_k$, where
$
BB_k:=\oplus_{\{(T,T') \ : \ |T\cap T'|=k\}} \ell_T\otimes\ell_{T'},
$
and $BB_k$ is defined over $\QQ$.  We get the decomposition  $ch(F_1)\otimes ch(F_2)=\sum_{i=0}^{e/2}c_k$, where $c_k\in BB_k$.
Let $\phi:H^*(X\times X)\rightarrow H^*(X\times\hat{X})$ be the isomorphism induced by Orlov's derive equivalence $\Phi$ and let $\tau$ act on $H^i(X)$ by multiplication by $(-1)^{i(i-1)/2}$. The image of $BB_1$ via $\phi\circ(id\otimes\tau)$ is contained in $\oplus_{i\geq d}H^i(X\times\hat{X},\QQ)$ and
the image of the composition 
\begin{equation}
\label{eq-composition-from-BB-1-to-HW}
BB_1\RightArrowOf{\phi\circ(id\otimes\tau)}\oplus_{i\geq d}H^i(X\times\hat{X},\QQ)\rightarrow H^d(X\times\hat{X},\QQ)
\end{equation}
is $HW(X\times\hat{X},\eta)$, by Lemma \ref{lemma-on-the-phi-image-of-lines-tensor-ell-T-ell-T-prime}. The above composition maps the line $\ell_T\otimes\ell_{T'}$ 
isomorphically onto the line $\wedge^d V_\sigma$, where $\{\sigma\}=T(\hat{\Sigma})\cap T'(\hat{\Sigma})$.
The dimension of $BB_1$ is\footnote{In contrast, $\dim(B\otimes B)=2^e.$ When $[K:\QQ]=4$, for example, $\dim HW(A,\eta)=4$, $\dim(BB_1)=8$, and $\dim(B\otimes B)=16$.
} 
$e2^{\frac{e}{2}-1}$ and $\dim HW(A,\eta)=e$, 
 so the above composition is an isomorphism if and only if $e:=[K:\QQ]=2$.
%Choose a non-zero element $t_\sigma\in\wedge^d V_\sigma$, for each $\sigma\in\Sigma$, and denote by $\lambda_{T,T'}\in\ell_T\otimes \ell_{T'}$ 
%the element mapping to $t_\sigma$. 
Write $c_1=\sum_{\{(T,T') \ : \ |T\cap T'|=1\}} c_{T,T'}$, where $c_{T,T'}$ is the direct summand in $\ell_T\otimes\ell_{T'}$. Let $\Pi$ be the composition (\ref{eq-composition-from-BB-1-to-HW}).

\begin{prop} 
\label{prop-introduction}
(Proposition \ref{prop-kappa-class-of-image-of-secant-class-yields-a-HW-class}(\ref{lemma-item-difference-projects-to-HW})). 
Assume that $d>2$, the rank of $\Phi(F_1\boxtimes F_2^\vee)$ is non-zero, and 
${\displaystyle \sum_{\{(T,T') \ : \ |T\cap T'|=\{\sigma\}\}} \Pi(c_{T,T'})
%\phi((id\otimes\tau)(c_{T,T'}))
\neq 0}$, 
for some $\sigma\in \Sigma$. Then $\kappa_{d/2}(E)$ projects to a non-zero class in the summand $HW(A,\eta)$ of $\A^d$ in (\ref{eq-decomposition-of-generic-Hodge-classes-in-H-d}). 
\end{prop}

In Example \ref{example-coherent-sheaf-on-Jacobian-of-genus-4-with-ch-a-multiple-of-beta-prime} and Lemma \ref{lemma-a-criterion-for-alpha-0-tensor-beta-to-map-via-orlov-equivalence-to-a-class-inducing-non-zero-HW-class}, where the above criterion is verified, we choose the sheaves $F_1$ and $F_2$ so that all but one summand $c_{T,T'}$ vanish, for each $\sigma\in\Sigma$, and so the non-vanishing of the sum in the above proposition, for some $\sigma\in \Sigma$, is equivalent to the non-vanishing of the summand $c_1$ in $BB_1$.

We summarize the above results as follows. Let $F_1$ and $F_2$ be $B$-secant sheaves satisfying the genericity criterion of Proposition 
\ref{prop-introduction}.

\begin{thm} 
\label{thm-main-CM-field}
(Corollary \ref{corollary-VHC})
A flat deformation of the class $\kappa(\Phi(F_1\boxtimes F_2^\vee))$ remains of Hodge type, on every abelian variety of Weil type $(A',\eta')$ deformation equivalent to $(X\times\hat{X},\eta)$.
If the flat deformation of the class $\kappa(\Phi(F_1\boxtimes F_2^\vee))$ remains furthermore algebraic, then every Weil class in $HW(A',\eta')$ is algebraic.
\end{thm}

In  \cite[Theorem 1.5.1]{markman-sixfolds} the flat deformations of the class $\kappa(\Phi(F_1\boxtimes F_2^\vee))$ were shown to remain algebraic, when $F=\QQ$, $X$ is the Jacobian of a genus $3$ curve,
and the $B$-secant sheaves $F_1$ and $F_2$ on $X$ were chosen to be equivariantly semi-regular in the sense of Buchweitz-Flenner \cite{buchweitz-flenner}. We postpone  for future work the search for semiregular $B$-secant sheaves for CM-fields with $[K:\QQ]>2$.

%****************************************************************
% 
%****************************************************************
\subsection{Organization of the paper}
\label{sec-organization}
A table of notation is provided in Section \ref{sec-CM-notation}.
Throughout the paper $F$ is a totally real field and
$X$ is an abelian variety endowed with an embedding $\hat{\eta}:F\rightarrow\End_\QQ(X)$.
In Section \ref{sec-exterior-algebra-over-F} we review the exterior subalgebra $\wedge^*_F H^1(X,\QQ)$ of $\wedge^*_\QQ H^1(X,\QQ)$. In addition, we review the relation between the $\QQ$-valued and $F$-valued quadratic forms on $V_\QQ:=H^1(X,\QQ)\oplus H^1(\hat{X},\QQ)$. We describe the tensor product factorization of even pure spinors in $S^+_\CC$ associated to the direct sum decomposition $V_\CC=\oplus_{\hat{\sigma}\in\hat{\Sigma}}V_{\hat{\sigma},\CC}$.

We denote by $V_{\hat{\eta}}$ the space $V_\QQ$ considered as a vector space over $F$ with its natural $F$-valued quadratic form.
In Section \ref{sec-spin-groups-over-F} we construct a natural homomorphism from $\Spin(V_{\hat{\eta}})$ to $\Spin(V_\QQ)$, which is injective if $[F:\QQ]$ is odd, and has kernel $\{\pm 1\}$ if $[F:\QQ]$ is even. We denote by $\Spin(V_\QQ)_{\hat{\eta}}$ the image of $\Spin(V_{\hat{\eta}})$ in $\Spin(V_\QQ)$.

In Section \ref{sec-complex-multiplication--by-CM-field-K} we construct an embedding $\eta:K\rightarrow \End(V_\QQ)$ of the CM-field $K$ associated to a maximal isotropic subspace $W$ of $V_{\hat{\eta}}\otimes_FK$. We construct a homomorphism $\Xi:K_-\rightarrow \wedge^2_FV_{\hat{\eta}}$ from the $-1$-eigenspace $K_-$ of the Galois involution $\iota\in Gal(K/F)$. We construct a maximal isotropic subspace $W_T$ of $V_\CC:=V\otimes_\ZZ\CC$ associated to $W$ and a choice of a CM-type $T$.

Let $\tilde{K}$ be the Galois closure of the extension $K$ of $\QQ$.
In Section \ref{sec-Galois-action-on-set-of-CM-types} we describe the $Gal(\tilde{K}/\QQ)$-action on the set $\{W_T \ : \ T \ \mbox{a CM-type}\}$ and the associated set of pure spinors $\{\ell_T\}$. We conclude that the subspace of $S^+_\CC$ spanned by the lines $\ell_T$ is defined over $\QQ$ and corresponds to a subspace $B$ of $S^+_\QQ$.

In Section \ref{sec-secant-space-B} we consider the subgroup $\Spin(V_\QQ)_\eta$ of $\Spin(V_\QQ)_{\hat{\eta}}$ respecting the direct sum decomposition
$V_{\tilde{K}}:=V_\QQ\otimes_\QQ\tilde{K}=\oplus_{\sigma\in\Sigma}V_\sigma.$ We denote by $\Spin(V_\QQ)_{\eta,B}$
the subgroup of $\Spin(V_\QQ)_\eta$ fixing all points in $B$. 
The character $\ell_{T_1}\otimes\ell_{T_2}$ of $\Spin(V_\QQ)_\eta$, associated to two CM-types $T_1$ and $T_2$, is proved to be the tensor product of
the characters $\wedge^d{V_\sigma}$, as $\sigma$ varies in the subset $T_1(\hat{\Sigma})\cap T_2(\hat{\Sigma})$ of $\Sigma$.
We conclude that the set of lines $\{\ell_T  \ : \ T \ \mbox{a CM-type}\}$ is linearly independent and  $\dim(B)=2^{e/2}.$
We construct a $\Spin(V_\QQ)_\eta$-invariant $K$-valued hermitian form $H_t:V_\QQ\times V_\QQ\rightarrow K$, associated to non-zero $t\in K_-$.
Finally we prove that if $B$ is spanned by Hodge classes, then so does the space $HW(X\times\hat{X},\eta)$ of Weil classes. 
The group 
%$\Spin(V_\QQ)_\eta$ is the Mumford-Tate group  and 
$\Spin(V_\QQ)_{\eta,B}$ is the special Mumford-Tate group of the generic abelian variety of Weil type deformation equivalent to $(X\times\hat{X},\eta)$.

In Section \ref{sec-Mumford-Tate-groups} we show that the subalgebra $\A:=(\wedge^*V_\QQ)^{\Spin(V_\QQ)_{\eta,B}}$ of the exterior algebra $\wedge^*V_\QQ$
of $\Spin(V_\QQ)_{\eta,B}$-invariant classes is generated by the $(e/2)$-dimensional $(\wedge^2V_\QQ)^{\Spin(V_\QQ)_{\eta,B}}$ and the $e$-dimensional $HW(X\times\hat{X},\eta).$

In Section \ref{sec-adjoint-orbit-n-Spin-V-RR-eta-B} we show that the subgroup $\Spin(V_\RR)_B$ of $\Spin(V_\RR)$, fixing all points of $B$, is the product
of the special unitary groups $SU(V_{\hat{\sigma},\RR})$, $\hat{\sigma}\in \hat{\Sigma}$ (see Remark \ref{rem-SU-V-hat-sigma-RR} for the definition of $SU(V_{\hat{\sigma},\RR})$). The group $\Spin(V_\RR)_B$ maps injectively into the image $SO_+(V_\RR)$ of $\Spin(V_\RR)$ in $SO(V_\RR)$. 
We show that a finite union $\Omega_B$ of adjoint orbits in $\Spin(V_\RR)_B$ of complex structures on $V_\RR$ is 
a period domain of polarized abelian varieties of Weil type with complex multiplication by $K$. 
The dimension of each adjoint orbit in $\Omega_B$ is $d^2e/2$, so each is a connected component of $\Omega_B$.
The set of connected components is in bijection with the set of CM-types.

In Section \ref{sec-HW-classes-pure-spinors-and-Orlov-equivalence} we introduce the grading $B\otimes B=\oplus_{i=0}^{e/2}BB_i$, where $BB_i\otimes_\QQ\tilde{K}$ is the direct sum of $\ell_T\otimes\ell_{T'}$, where the ordered pairs $(T,T')$ consist of CM-types $T,T':\hat{\Sigma}\rightarrow\Sigma$ 
having precisely $i$ common values. We show that $HW(X\times\hat{X},\eta)$ is the image of a natural homomorphism $BB_1\rightarrow \wedge^dV_\QQ$, given in (\ref{eq-composition-from-BB-1-to-wedge-d}), and describe the kernel. Finally we prove Theorem \ref{thm-main-CM-field}.

Section \ref{sec-example-isometry-g-0} consists of Examples. We construct a maximal isotropic subspace $W$ of $[H^1(X,\QQ)\oplus H^1(\hat{X},\QQ)]\otimes_FK$ as the image of $H^1(\hat{X},\QQ)\otimes_FK$ under the element of $\Spin(V_{\hat{\eta}}\otimes_FK)$ of cup product with $\exp(\sqrt{-q}\Theta)$,
where $\Theta$ is a polarization of $X$ in $\wedge^2_FH^1(X,\QQ)$, and $q\in F$ satisfies $K=F(\sqrt{-q})$. We get the embedding $\eta:K\rightarrow \End_\QQ(X\times\hat{X})$ 
and the rational subspace $B$ of $S^+_\QQ:=H^{ev}(X,\QQ)$ described above. 
We show that the adjoint orbit in $\Spin(V_\RR)_B$ of the complex structure $I_{X\times\hat{X}}$ of $X\times\hat{X}$ is a period domain of abelian varieties of split Weil type (Definition \ref{def-split-type}) deformation equivalent to $(X\times\hat{X},\eta)$ (Lemma \ref{lemma-I-X-hat-X-belongs-to-Omega-B}).

In Section \ref{sec-example-CM-field-with-order-4-Galois-group} we consider the case, where $[F:\QQ]=2$ and $K=F(\sqrt{-q})$, where $q$ is a positive integer. 
We choose $X$ to be a Jacobian of a genus $4$ curve with an embedding $\hat{\eta}:F\rightarrow \End_\QQ(X)$. 
We give an example of two coherent sheaves $F_1$ and $F_2$ on $X$ satisfying the hypothesis of Theorem \ref{thm-main-CM-field}, so that flat deformations of $\kappa(\Phi(F_1\boxtimes F_2^\vee))$ remain of Hodge type under all deformations of $(X\times\hat{X},\eta)$ in the connected component of the period domain $\Omega_B$, but we do not know if the class $\kappa(\Phi(F_1\boxtimes F_2^\vee))$ remains algebraic (Lemma \ref{lemma-box-tensor-product-of-sheaves-maps-to-one-with-equivariant-kappa-class-case-of-degree-4-extension}). 
If it does remain algebraic, then so do the Weil classes, by Lemma \ref{lemma-box-tensor-product-of-sheaves-maps-to-one-with-equivariant-kappa-class-case-of-degree-4-extension}. 

In Section \ref{sec-Gal-K-is-cartesion-product-of-Gal-F-and-Z-2} we consider more general totally real fields $F$, but specialize to the case where $Gal(\tilde{K}/\QQ)$ splits as the direct product of $Gal(\tilde{F}/\QQ)$ and $\{id,\iota\}$, where $\iota$ is the extension to $\tilde{K}$ of the involution in $Gal(K/F)$. We prove a genericity criterion for the Chern characters of coherent sheaves $F_1$ and $F_2$, so that $\kappa_{d/2}(\Phi(F_1\boxtimes F_2^\vee)$ is the sum of a non-zero element of $HW(X\times\hat{X},\eta)$ and a polynomial in Hodge (necessarily algebraic)  classes in $H^2(X\times\hat{X},\QQ)^{\Spin(V_\QQ)_{\eta,B}}$
(Lemma \ref{lemma-a-criterion-for-alpha-0-tensor-beta-to-map-via-orlov-equivalence-to-a-class-inducing-non-zero-HW-class}). 
We do not provide examples of such $F_1$ and $F_2$, when $[F:\QQ]>2$. Example of semiregular such $F_1$ and $F_2$ are likely to yield the algebraicity of the Weil classes of all abelian varieties of Weil type deformation equivalent to $(X\times\hat{X},\eta).$

\begin{rem}
We have $2^{[F:\QQ]}$ maximal isotropic subspaces $W_{T,\CC}$ in the quadratic vector space $V_\CC:=H^1(X\times\hat{X},\CC)$, one for each CM-type $T$. Nevertheless, each direct summand $V_{\hat{\sigma},\CC}$, $\hat{\sigma}\in \hat{\Sigma}$, of $V_\CC$  decompose as the direct sum of two maximal isotropic subspaces $V_{\hat{\sigma},\CC}=V_{\sigma,\CC}\oplus V_{\bar{\sigma},\CC}$, where $\sigma$ and $\bar{\sigma}$ are the two embeddings of $K$ in $\CC$ restricting
to the embedding $\hat{\sigma}$ of $F$. As a result, many of the proofs in Sections \ref{sec-secant-space-B} to \ref{sec-HW-classes-pure-spinors-and-Orlov-equivalence} involve a reduction of a statement to a direct summand $V_{\hat{\sigma},\CC}$ in $V_\CC$ or to a pure spinor in $S^+_{\hat{\sigma}}$, which is a tensor  factor of a pure spinor $\ell_T$, and then invoke an argument from \cite{markman-sixfolds} 
in the quadratic field extension case. 
\end{rem}
%\begin{rem}
%\label{rem-notation-T}
%A word about the notation. The subspace $W_T$ of $H^1(A,\QQ)\otimes_FK$ is independent of the CM-type $T$, but $T$ determines an embedding 
%(\ref{eq-embedding-of-V-otimes-FK-in-V-CC})
%of 
%$H^1(A,\QQ)\otimes_FK$ in $H^1(A,\CC)$, and so the subspace $W_{T,\CC}$ of $H^1(A,\CC)$. 
%The subspace $W_T$ will facilitate the construction of the complex multiplication $\eta$, which is independent of $T$.
%Once constructed, we will get a maximal isotropic subspace $W_{T'}$ of $H^1(A,\CC)$, for every CM-type $T'\in \T_K$, 
%so that $W_{T,\CC}$ is the one associated to $T$. 
%\end{rem}

\section{Notation} 
\label{sec-CM-notation}
\hspace{0ex}

\begin{longtable}{l l l}
%\begin{table}[H]
%\caption{Glossary of Notation}
%\label{table}
%\begin{tabular}{l l l}
%\hline
%\hline
\\
%\hline
$X$ & a complex abelian variety & $\dim_\CC(X)=n=de/4$
\\
$V$ & the lattice $H^1(X,\Integers)\oplus H^1(\hat{X},\Integers)$ & $\rank(V)=4n=de$
\\
$K$ & a CM field & $[K:\QQ]=e$
\\
$e$ & $[K:\QQ]$ & Sec. \ref{sec-summary-of-the-construction}
\\
%\hline
$F$ & the totally real subfield of $K$ & $[F:\QQ]=e/2$
\\
$\hat{\eta}$ & an algebra homomorphism $F\rightarrow \End_{Hdg}(V_\QQ)$ & (\ref{eq-hat-eta})
\\
$V_{\hat{\eta}}$ & $V_\QQ$ considered as a vector space over $F$ & $\dim_F(V_{\hat{\eta}})=2d$
\\
$V_{\hat{\eta},K}$ & $V_{\hat{\eta}}\otimes_F K$  & Sec. \ref{sec-complex-multiplication--by-CM-field-K}
\\
$(\bullet,\bullet)_V$ & the bilinear pairing on the lattice $V$ & (\ref{eq-pairing-on-V})
\\
$(\bullet,\bullet)_{V_{\hat{\eta}}}$ & $F$-valued $F$-bilinear pairing & (\ref{eq-pairing-on-V-hat-eta})
\\
$(\bullet,\bullet)_{\hat{\eta},K}$ & the $K$-valued pairing on $V_{\hat{\eta},K}$ & (\ref{eq-pairing-on-V-hat-eta-K})
\\
$C(V)$ & the Clifford algebra & Sec. \ref{sec-factorization-of-pure-spinors}
\\
$\rho$ & the vector representation $\Spin(V)\rightarrow SO_+(V)$ & (\ref{eq-rho})
\\
$S$ & $H^*(X,\ZZ)$ considered as the spin representation & Sec. \ref{sec-factorization-of-pure-spinors}
\\
$S^+$ &  the half spin representation $H^{ev}(X,\ZZ)$ & Sec. \ref{sec-factorization-of-pure-spinors}
\\
$S^-$ &  the half spin representation $H^{odd}(X,\ZZ)$& Sec. \ref{sec-factorization-of-pure-spinors}
\\
$m$ & the spin representation $m:C(V)\rightarrow \End(S)$ & (\ref{eq-m}) 
\\
$m$ & the spin representation $m:\Spin(V)\rightarrow GL(S)$ & (\ref{eq-spin-representation-m})
\\
$\eta$ & an algebra homomorphism $K\rightarrow \End_{Hdg}(V_\QQ)$ & (\ref{eq-eta})
\\
$V_\eta$ & $V_\QQ$ considered as a vector space over $K$ & $\dim_K(V_{\eta})=d$
\\
$\tilde{F}$ & the Galois closure of $F$ over $\QQ$ &
\\
$\Sigma$ & set of embeddings of the CM-field $K$ in $\CC$ &
\\
$\hat{\Sigma}$ & set of embeddings of the field $F$ in $\RR$
\\
$\hat{\sigma}$ & an embedding $F\rightarrow \RR$ &
\\
$V_{\hat{\sigma}}$ & the direct summand of $V_{\hat{\eta}}\otimes_F\tilde{F}\cong\oplus_{\hat{\sigma}\in\hat{\Sigma}}V_{\hat{\sigma}}$ & $\dim_{\tilde{F}}(V_{\hat{\sigma}})=2d$
\\
$V_{\hat{\sigma},\RR}$ & & (\ref{eq-V-hat-sigma-RR})
\\
$\tilde{K}$ & the Galois closure of $K$ over $\QQ$ &
\\
$\sigma$ & an embedding $K\rightarrow \CC$ &
\\
$V_{\sigma}$ & the direct summand of $V_{\eta}\otimes_K\tilde{K}\cong\oplus_{\sigma\in\Sigma}V_\sigma$ & $\dim_{\tilde{K}}(V_{\sigma})=d$
\\
$V^{1,0}_\sigma$, $V^{0,1}_\sigma$ & the Hodge direct summands of $V_{\sigma,\CC}$ & $\dim_\CC(V^{1,0}_\sigma)=d/2$
\\
$W$ & a maximal isotropic subspace of $V_{\hat{\eta}}\otimes_FK$ & $\dim_K(W)=d$
\\
$T$ & a CM-type for the CM-field $K$ & Sec. \ref{sec-summary-of-the-construction}
\\
$\T_K$ & The set of CM-types for the CM-field $K$ & Sec. \ref{sec-Galois-action-on-set-of-CM-types}
\\
$e_T$ & embedding $V_{\hat{\eta}}\otimes_FK\rightarrow V_{\tilde{K}}$ determined by $T$ & (\ref{eq-e-T-over-tilde-K})
\\
$W_T\subset V_{\tilde{K}}$ & $\span_{\tilde{K}}(e_T(W))$ & $\dim_{\tilde{K}}(W_T)=2n=\frac{de}{2}$
\\
$W_{T,\hat{\sigma}}$ & the direct summand of $W_T$ & $\dim_{\tilde{K}}(W_{T,\hat{\sigma}})=d$
\\
$\Spin(V_\bullet)$ & Spin group of $V\otimes_\ZZ\bullet$ for a field $\bullet$ &
\\
$SO_+(V_\bullet)$ & the image of $\Spin(V_\bullet)$ in $SO(V_\bullet)$ &
\\
$\Spin(V)$ & the arithmetic subgroup of $\Spin(V_\QQ)$ & %preserving the lattice $V$
\\
$\Spin(V_{\hat{\eta}})$ & spin group of $V_\QQ$ as an $F$ vector space &
\\
$\Spin(V_\eta)$ & spin group of $V_\QQ$ as a $K$ vector space &
\\
$\Spin(V_\QQ)_{\hat{\eta}}$ & subgroup of $\Spin(V_\QQ)$ defined in (\ref{eq-Spin-V-QQ-hat-eta}) &
\\
$\Spin(V)_\eta$ & subgroup of $\Spin(V)_{\hat{\eta}}$ defined in (\ref{eq-Spin-V-eta}) &
\\
$\Spin(V_{\tilde{K}})_\eta$ &&(\ref{eq-Spin-V-tilde-K-eta})
\\
$\zeta$ & homomorphism $\zeta:\Spin(V_{\hat{\eta}})\rightarrow \Spin(V_\QQ)_{\hat{\eta}}$ & (\ref{eq-homomorphism-from-Spin-V-hat-eta-to-commutator-in-Spin-V-QQ})
%\\
%$S^+$ & $H^{ev}(X,\ZZ)$ as the half-spin representation & $\rank(S^+)=2^{(de)/2-1}$
\\
$\ell_T$ & the pure spinor of $W_T$  in $S^+_{\tilde{K}}$& Sec. \ref{sec-Galois-action-on-set-of-CM-types}
\\
$B$ & subspace of $H^{ev}(X,\QQ)$ spanned by the $\ell_T$'s & $\dim_\QQ(B)=2^{e/2}$
\\
$B_{\sqrt{-q}\Theta}$ & $B$ associated to $\Theta\!\!\in\!\!\wedge^2_FH^1\!(\!X,\!\QQ)$ and $\sqrt{-q}\!\!\in\!\! K$ & Sec. \ref{sec-example-CM-field-with-order-4-Galois-group}
\\
$K_-$ & subspace of $K$ of purely imaginary elements & Sec. \ref{sec-F-invariant-2-forms}, \ref{sec-Hermitian-form-H-t}
\\
$\tilde{\Xi}_t$ & an $F$-valued $2$-form in $\wedge^2_FV^*_{\hat{\eta}}$ ,  $t\in K_-$ & (\ref{eq-Xi-t})
\\
$\Xi_t$ & the $2$-form $tr_{F/\QQ}\circ\Xi_t$ in $\wedge^2_\QQ V^*_\QQ$,  $t\in K_-$ & (\ref{eq-Xi-t})
\\
$H_t$ & a $K$-valued hermitian form on $V_\eta$,  $t\in K_-$ & (\ref{eq-H-t})
\\
$\Spin(V)_{\eta,B}$ & subgroup of $\Spin(V)_\eta$ fixing all points of $B$ & (\ref{eq-Spin-V-eta-B})
\\
$\det_\sigma$ & character of $\Spin(V_{\tilde{K}})_\eta$, $\sigma\in\Sigma$ & (\ref{eq-det-sigma}) 
\\
$\Spin(V_\RR)_B$ & subgroup of $\Spin(V_\RR)$ fixing all points of $B$ & Sec. \ref{sec-Spin-V-R-B}
\\
$\Omega_{B,t}$ & period domain of Weil type abelian varieties  & (\ref{eq-omega-B-t})
%Sec. \ref{sec-period-domain-for-abelian-varieties-of-Weil-type-CM-case}
\\
$\Omega_{B,T}$ & $\Omega_{B,t}$, $T$ is the CM-type associated to $t\in K_-$ & Sec. \ref{sec-period-domain-for-abelian-varieties-of-Weil-type-CM-case}
\\
$\Omega_B$ & $\cup_{T\in \T_K}\Omega_{B,T}$ & (\ref{eq-omega-B}) %Sec. \ref{sec-organization}
\\
$S^+_{\hat{\sigma}}$ & the half-spin representation of $\Spin(V_{\hat{\sigma}})$ &
\\
$\ell_\sigma$ & the pure spinor of $V_\sigma$ in $S^+_{\hat{\sigma}}$ & 
\\
$P_{\hat{\sigma}}$ & the plane $\span\{\ell_\sigma,\ell_{\hat{\sigma}}\}$ in $S^+_{\hat{\sigma}}$ &
\\
$HW(X\!\!\times\!\!\hat{X}\!,\!\eta)$ & $\wedge^d_{\eta(K)}H^1(X\times\hat{X},\QQ)$ & (\ref{eq-HW-CM-field-case})
\\
$\A$ & $\Spin(V)_{\eta,B}$-invariant subalgebra of $\wedge^*V_\QQ$ & Sec. \ref{sec-Mumford-Tate-groups}
\\
$\R$ & $\Spin(V)_\eta$-invariant subalgebra of $\wedge^*V_\QQ$ & Sec. \ref{sec-Mumford-Tate-groups}
\\
$BB_i$ & a graded summand of $B\otimes_\QQ B$ & Sec. \ref{sec-grading-of-BB}
\\
$KB_1$ & The kernel of $BB_1\rightarrow HW(X\times\hat{X},\eta)$ & (\ref{eq-composition-from-BB-1-to-wedge-d}) 
\\
$F^k(\wedge^*V)$ & the increasing filtration $\oplus_{i\leq k}\wedge^iV$ of $\wedge^*V$ & Sec. \ref{sec-grading-of-BB}
\\
$F_k(\wedge^*V)$ & the decreasing filtration $\oplus_{i\geq k}\wedge^iV$ of $\wedge^*V$ & Sec. \ref{sec-grading-of-BB}
\\
$\tau$ & multiplication of $H^i(X,\ZZ)$ by $(-1)^{i(i-1)/2}$ & (\ref{eq-tau})
%\end{tabular}
\end{longtable}
%****************************************************************
% 
%****************************************************************
\section{Field embeddings and their associated decompositions}
\label{sec-exterior-algebra-over-F}
%If $\sigma:K\rightarrow \CC$ is an embedding and we denote by $\hat{\sigma}:F\rightarrow \RR$ its restriction, then 
We have the decomposition 
$
F\otimes_\QQ\RR \cong \prod_{\hat{\sigma}\in\hat{\Sigma}}\RR.
$
Given an F vector space $U$ we get the decomposition $U\otimes_\QQ\RR=U\otimes_F(F\otimes_\QQ\RR)=\oplus_{\hat{\sigma}\in\hat{\Sigma}}U_{\hat{\sigma}}$,
where $U_{\hat{\sigma}}$ is the $\hat{\sigma}$-isotypic subspace, i.e., the subspace of $U\otimes_\QQ\RR$, where $F$ acts via the character $\hat{\sigma}:F\rightarrow\RR$. 
The embedding $\hat{\eta}:F\rightarrow \End_{Hdg}(H^1(X,\QQ))$ endows the singular cohomology $H^1(X,\QQ)$ with the structure of a vector space over $F$. We get the decomposition
\[
H^1(X,\RR)=H^1(X,\QQ)\otimes_\QQ\RR\cong \oplus_{\hat{\sigma}\in\hat{\Sigma}} H^1_{\hat{\sigma}}(X),
\]
where $H^1_{\hat{\sigma}}(X)$ denotes the $\hat{\sigma}$-isotypic  subspace.
%Let 
%\[
%\hat{\eta}:F\rightarrow \End(H^1(X,\QQ))
%\] 
%be the embedding. 
The dimensions $\dim_\RR H^1_{\hat{\sigma}}(X)$ are the same for all $\hat{\sigma}\in \hat{\Sigma}$ and are equal to $d=\dim_FH^1(X,\QQ)$.
We have the equality
\[
\dim(H^{1,0}_{\hat{\sigma}}(X))=\dim(H^{0,1}_{\hat{\sigma}}(X))=d,
\]
for all $\hat{\sigma}\in\hat{\Sigma}$, since the $\hat{\eta}(F)$-action commutes with the complex structure of $X$, by assumption, and with complex conjugation, as $F$ is totally real.

%****************************************************************
% 
%****************************************************************
\subsection{Trace forms}
Denote by $H^1(X,\QQ)_{\hat{\eta}}$ the vector space $H^1(X,\QQ)$ considered as an $F$ vector space.
We have the isomorphism of vector spaces over $\QQ$
\[
\Hom_F(H^1(X,\QQ)_{\hat{\eta}},F) \rightarrow \Hom_\QQ(H^1(X,\QQ),\QQ)
\]
sending $h\in \Hom_F(H^1(X,\QQ)_{\hat{\eta}},F)$ to $tr_{F/\QQ}\circ h$.
We denote by 
\begin{equation}
\label{eq-hat-eta}
\hat{\eta}:F\rightarrow \End(V_\QQ)
\end{equation}
the scalar multiplication action of $F$ on $H^1(X,\QQ)_{\hat{\eta}}\oplus \Hom_F(H^1(X,\QQ)_{\hat{\eta}},F)$. We denote the image of $a\in F$ by $\hat{\eta}_a$.
Explicitly, $\hat{\eta}_a$ acts on $\Hom_F(H^1(X,\QQ)_{\hat{\eta}},F)$ by 
$\hat{\eta}_a(\theta)=\theta\circ \hat{\eta}_a$, where on the left $\hat{\eta}_a$ is $\End(V_\QQ)$ and on the right it is in $\End(H^1(X,\QQ))$.
We defined $V_{\hat{\eta}}$ as $V_\QQ$ with its $F$-vector space structure. 
Define the $\hat{\eta}(F)$-bilinear pairing 
%\footnote{Here we denote by $(\bullet,\bullet)_{V_{\hat{\eta}}}$ both the the map from $V_\QQ\times V_\QQ$ to $F$ and the unique 
%linear homomorphism $V_\QQ\otimes_{\hat{\eta}(F)} V_\QQ\rightarrow F$ determined by the former.}  
\begin{equation}
\label{eq-pairing-on-V-hat-eta}
(\bullet,\bullet)_{V_{\hat{\eta}}} : V_\QQ\times V_\QQ\rightarrow F,
\end{equation}
by $((w_1,\theta_1),(w_2,\theta_2))_{V_{\hat{\eta}}}=\theta_1(w_2)+\theta_2(w_1)$, where $\theta_1$ and $\theta_2$ are regarded as a elements of 
$\Hom_F(H^1(X,\QQ)_{\hat{\eta}},F)$.
The $\QQ$-valued pairing 
$(\bullet,\bullet)_{V_\QQ}$ 
is, by definition, 
the composition $tr_{F/\QQ}\circ (\bullet,\bullet)_{V_{\hat{\eta}}}.$ The $\hat{\eta}(F)$-bilinearity of $(\bullet,\bullet)_{V_{\hat{\eta}}}$ implies that $\hat{\eta}(a)$ is 
self-dual with respect to the pairing
$(\bullet,\bullet)_{V_\QQ}$.

Consider the $\RR$-bilinear $\RR$-valued pairing on $F\otimes_\QQ\RR$ given by
\[
(f_1\otimes r_1,f_2\otimes r_2)=r_1r_2tr_{F/\QQ}(f_1f_2),
\]
for $f_i\in F$ and $r_i\in \RR$, $i=1,2$.
The decomposition $F\otimes_\QQ\RR=\oplus_{\hat{\sigma}\in\hat{\Sigma}}\RR$ consists of pairwise orthogonal lines with respect to the above pairing. Set 
\begin{eqnarray}
\label{eq-V-hat-sigma-RR}
V_{\hat{\sigma},\RR}&:=&H^1_{\hat{\sigma}}(X)\oplus H^1_{\hat{\sigma}}(\hat{X}).
%\\
%V_{\hat{\sigma}}&:=&[H^1_{\hat{\sigma}}(X)\oplus H^1_{\hat{\sigma}}(\hat{X})]\cap H^1(X\times\hat{X},\hat{\sigma}(F)).
\end{eqnarray}
%The composition $V_\QQ\rightarrow V_\RR\rightarrow V_{\hat{\sigma},\RR}$ maps $V_\QQ$ isomorphically onto $V_{\hat{\sigma}}$. 

\begin{lem}
\label{lemma-orthogonal-direct-sum-decomposition-V-hat-sigma}
The decomposition $V\otimes_\ZZ\RR=\oplus_{\hat{\sigma}\in\hat{\Sigma}} V_{\hat{\sigma},\RR}$ consists of pairwise orthogonal subspaces with respect to $(\bullet,\bullet)_{V_\RR}.$
\end{lem}

\begin{proof}
This follows from the self-duality of
 $\hat{\eta}(a)$, for $a\in F$. Explicitly, the pairing $(\bullet,\bullet)_{V_\RR}:V_\RR\otimes V_\RR\rightarrow \RR$  factors as an $\hat{\eta}(F)$-bilinear pairing 
$(\bullet,\bullet)_{V_{\hat{\eta}}}$ with values in $F\otimes_\QQ\RR$ composed with $tr_{F/\QQ}\otimes id_\RR:F\otimes_\QQ\RR\rightarrow \RR$. Now, given $x\in V_{\hat{\sigma}_1,\RR}$, $y\in V_{\hat{\sigma}_2,\RR}$, and $q\in F$ a primitive element,  we have
\begin{eqnarray*}
q(x,y)_{V_{\hat{\eta}}}&=&(\hat{\eta}_q x,y)_{V_{\hat{\eta}}}=(\hat{\sigma}_1(q)x,y)_{V_{\hat{\eta}}}=
\hat{\sigma}_1(q)(x,y)_{V_{\hat{\eta}}}
\\
q(x,y)_{V_{\hat{\eta}}}&=&(x,\hat{\eta}_q y)_{V_{\hat{\eta}}}= (x,\hat{\sigma}_2(q)y)_{V_{\hat{\eta}}}=\hat{\sigma}_2(q)(x,y)_{V_{\hat{\eta}}},
\end{eqnarray*}
where multiplication by $q$ on the left above corresponds to multiplication in the left tensor factor of $F\otimes_\QQ\RR$ and multiplication by 
$\hat{\sigma}_i(q)$ on the right above corresponds to multiplication in the right tensor factor of $F\otimes_\QQ\RR$. Hence, if $\hat{\sigma}_1\neq\hat{\sigma}_2$, then
$(x,y)_{V_{\hat{\eta}}}=0.$ Consequently, its trace $(x,y)_{V_\RR}$ vanishes as well.
\end{proof}

The $j$-th cartesian power $(F^\times)^j$ of the multiplicative group $F^\times$ acts on the $j$-th tensor power over $\QQ$ of $H^1(X,\QQ)$ and hence
on $\wedge^j_\QQ H^1(X,\QQ)$ and on $\wedge^j_\QQ H^1(X,\QQ)\otimes_\QQ\RR\cong\wedge^j_\RR H^1(X,\RR)$.
The latter decomposes as a direct sum of characters subspaces 
\[
\wedge^j_\RR H^1(X,\RR)=
\oplus_{(\hat{\sigma}_1, \dots, \hat{\sigma}_j)\in (\hat{\Sigma})^j}\left(\wedge^j_\RR H^1(X,\RR)\right)_{(\hat{\sigma}_1, \dots, \hat{\sigma}_j)}
\]
The vector space $\wedge^j_\QQ H^1(X,\QQ)$ decomposes into subrepresentations of $(F^\times)^j$ defined over $\QQ$. The subrepresentation $\wedge^j_FH^1(X,\QQ)$ of $\wedge^j_\QQ H^1(X,\QQ)$  corresponds to the subrepresentation
\[
\oplus_{\hat{\sigma}\in\hat{\Sigma}}\wedge^j_\RR H^1_{\hat{\sigma}}(X)= 
\oplus_{\hat{\sigma}\in\hat{\Sigma}}\left(\wedge^j_\RR H^1(X,\RR)\right)_{(\hat{\sigma}, \dots, \hat{\sigma})}.
\]
of $\wedge^j_\RR H^1(X,\RR)$. Given $\alpha_i\in \wedge^{j_i}_F H^1(X,\QQ)$, $i=1,2$, then $\alpha_1\wedge_F\alpha_2$ 
is the projection of $\alpha_1\wedge_\QQ\alpha_2$ to the direct summand $\wedge_F^{j_1+j_2}H^1(X,\QQ)$.

The above action of $(F^\times)^j$ on $\wedge^j_FH^1(X,\QQ)$ factors through the natural action of $F$ on it as an $F$-vector space via the product homomorphisn $(F^\times)^j\rightarrow F^\times\subset F$, $(f_1,\dots, f_j)\mapsto \prod_{i=1}^j f_j$.

%****************************************************************
% 
%****************************************************************
\subsection{Tensor product factorizations of pure spinors}
\label{sec-factorization-of-pure-spinors}
Set $V:=H^1(X,\ZZ)\oplus H^1(\hat{X},\ZZ)$ endowed with the unimodular even symmetric bilinear pairing (\ref{eq-pairing-on-V}). Set $Q(v):=(v,v)_V/2$. 
Let $C(V)$ be the Clifford algebra $\oplus_{k=0}^\infty V^{\otimes k}/\langle 
%v_1v_2+v_2v_1-(v_1,v_2) \ : \ v_1, v_2
v^2-Q(v) \ : \ v
\in V\rangle$, where the integer $Q(v)$ is in $V^{\otimes 0}=\ZZ$. $C(V)$ is $\ZZ/2\ZZ$ graded, $C(V)=C(V)^{ev}\oplus C(V)^{odd}$. 
Set $S:=H^*(X,\ZZ)$, $S^+:=H^{ev}(X,\ZZ)$, and $S^-:=H^{odd}(X,\ZZ)$. We identify $V$ with the image of $V^{\otimes 1}$ in $C(V)$.
The Clifford group is the subgroup $G(V)$ of the group of invertible elements $g$ in $C(V)$ satisfying $gVg^{-1}=V$.
The latter conjugation is an isometry of $V$ and 
we get the {\em vector representation} $\rho:G(V)\rightarrow O(V)$.
By definition, $\Spin(V)$ is an index $2$ subgroup of $G(V)\cap C(V)^{ev}$ (kernel of the norm character) \cite{golyshev-luntz-orlov}.  The vector representation restricts to 
\begin{equation}
\label{eq-rho}
\rho:\Spin(V)\rightarrow O(V)
\end{equation}
and we denote its image by $SO_+(V)$. The kernel of $\rho$ is of order $2$ generated by $-1$ (the element $-1\in C(V)$).
We have the isomorphism 
\begin{equation}
\label{eq-m}
m:C(V)\rightarrow \ \End(S),
\end{equation}
by \cite[ Prop. 3.2.1(e)]{golyshev-luntz-orlov}, which restricts to $m:V\rightarrow \End(S)$ 
and to the spin representation 
\begin{equation}
\label{eq-spin-representation-m}
m:\Spin(V)\rightarrow GL(S).
\end{equation}

Set $S_{\hat{\sigma}}:=\wedge^*H^1_{\hat{\sigma}}(X)$. 
Note the isomorphisms $S_\RR:=\wedge^*H^1(X,\RR)\cong\otimes_{\hat{\sigma}\in \hat{\Sigma}}S_{\hat{\sigma}}$ and 
\[
S_\RR\otimes S_\RR\cong \otimes_{\hat{\sigma}\in \hat{\Sigma}}(S_{\hat{\sigma}}\otimes_\RR S_{\hat{\sigma}}),
\]
both in the categories of $\ZZ/2\ZZ$ graded algebras. 
The tensor product $V_{\hat{\eta}}\otimes_\QQ\RR$, with the induced bilinear pairing, decomposes as an orthogonal direct sum $\oplus_{\hat{\sigma}\in \hat{\Sigma}}V_{\hat{\sigma},\RR}$, by Lemma \ref{lemma-orthogonal-direct-sum-decomposition-V-hat-sigma}. Hence, $C(V_\QQ)\otimes_\QQ\RR$ is isomorphic to the tensor product 
$\otimes_{\hat{\sigma}\in\hat{\Sigma}}C(V_{\hat{\sigma},\RR})$, where $V_{\hat{\sigma},\RR}$ is endowed with the restriction of the bilinear pairing of $V_\RR$ and
 the tensor product is over $\RR$ in the category of $\ZZ_2$-graded algebras \cite[Lemma V.1.7]{lam}.
 In the current section \ref{sec-factorization-of-pure-spinors} we will denote $V_{\hat{\sigma},\RR}$ by $V_{\hat{\sigma}}$.
%and $C(V_\RR)\cong \otimes_{\hat{\sigma}\in\hat{\Sigma}}C(V_{\hat{\sigma}})$, 
The isomorphism $\varphi:S_\RR\otimes S_\RR\rightarrow C(V_\RR)$, given in \cite[III.3.1]{chevalley},  factors as the tensor product of 
\[
\varphi_{\hat{\sigma}}:S_{\hat{\sigma}}\otimes S_{\hat{\sigma}}\rightarrow %\otimes_{\hat{\sigma}\in\hat{\Sigma}}
C(V_{\hat{\sigma}}).
\] 

Given maximal isotropic subspaces $W_{\hat{\sigma}}$ of $V_{\hat{\sigma}}$, for all $\hat{\sigma}\in \hat{\Sigma}$, set $W:=\oplus_{\hat{\sigma}\in \hat{\Sigma}}W_{\hat{\sigma}}$. Let $\ell_{\hat{\sigma}}\subset S^+_{\hat{\sigma}}$ be the line spanned by a pure spinor of $W_{\hat{\sigma}}$ and by $\ell\subset S^+_\RR$ 
the line spanned by a pure spinor of $W$. Then 
\begin{equation}
\label{eq-tensor-product-decomposition-of-pure-spinor}
\ell=\bigotimes_{\hat{\sigma}\in \hat{\Sigma}}\ell_{\hat{\sigma}}.
\end{equation}
Given another set of maximal isotropic subspaces $W'_{\hat{\sigma}}$ of $V_{\hat{\sigma}}$, for all $\hat{\sigma}\in \hat{\Sigma}$, with pure spinors $\ell'_{\hat{\sigma}}$, we have
$\varphi(\ell\otimes\ell')=\otimes_{\hat{\sigma}\in \hat{\Sigma}}\varphi_{\hat{\sigma}}(\ell_{\hat{\sigma}}\otimes\ell'_{\hat{\sigma}}).$

Similarly, the isomorphism $\psi:C(V_\RR)\rightarrow \wedge^*V_\RR$, given in \cite[Eq. (2.3.1)]{markman-sixfolds}, decomposes as a tensor product of
$\psi_{\hat{\sigma}}:C(V_{\hat{\sigma}})\rightarrow \wedge^*V_{\hat{\sigma}}$ and so 
\begin{equation}
\label{eq-tilde-varphi}
\tilde{\varphi}:=\psi\circ\varphi
\end{equation} 
decomposes as the tensor product of 
\begin{equation}
\label{eq-tilde-varphi-hat-sigma}
\tilde{\varphi}_{\hat{\sigma}}:=\psi_{\hat{\sigma}}\circ \varphi_{\hat{\sigma}}:S_{\hat{\sigma}}\otimes S_{\hat{\sigma}}\rightarrow \wedge^*V_{\hat{\sigma}}.
\end{equation}

%**********
% Hide
%**********
\hide{
maps $S_{\hat{\sigma}_0}\otimes S_{\hat{\sigma}_0}$ to 
$C(V_{\hat{\sigma}_0})\otimes \otimes_{\{\hat{\sigma} \ : \ \hat{\sigma}\neq\hat{\sigma}_0\}}\wedge^{top}H^1_{\hat{\sigma}}(\hat{X})$,
where we regard $\wedge^{top}H^1_{\hat{\sigma}}(\hat{X})$ as a line in $C(V_{\hat{\sigma}})$.
Given a pure spinor $\ell$ in $S^+_{\hat{\sigma}_0}$ of a maximal isotropic subspace $W$ in $V_{\hat{\sigma}_0}$, the isomorphism $\varphi$ maps 
$\ell\otimes\ell$ to 
the tensor product of $\wedge^{top}W$ (??? $\varphi_{\hat{\sigma}_0}(\ell\otimes\ell)$ ???) with $\otimes_{\{\hat{\sigma} \ : \ \hat{\sigma}\neq\hat{\sigma}_0\}}\wedge^{top}H^1_{\hat{\sigma}}(\hat{X})$, which projects, in the appropriate graded summand, onto
the top wedge product of the maximal isotropic subspace $W\oplus \oplus_{\{\hat{\sigma} \ : \ \hat{\sigma}\neq\hat{\sigma}_0\}} H^1_{\hat{\sigma}}(\hat{X})$.
The pure spinor $\ell$, as a line on the left hand side, is mapped to its tensor product with the pure spinors $1_{\hat{\sigma}}$, $\hat{\sigma}\neq\hat{\sigma}_0$, on the right hand side.

\begin{rem} 
%\label{rem-image-via-varphi-of-top-wedge-over-F-of-a-maximal-isotropic-F-subspace}
Let $W$ be a maximal isotropic $F$-subspace of $V_{\hat{\eta}}$. Then $W_\RR=\oplus_{\hat{\sigma}\in\hat{\Sigma}}W_{\hat{\sigma}}$, and 
$W_{\hat{\sigma}}$ is maximal isotropic in $V_{\hat{\sigma}}$. The $e/2$-dimensional $\QQ$-subspace $\wedge^d_FW$ of $\wedge^d_\QQ V_\QQ$ spans over $\RR$ the subspace
spanned by the lines $\wedge^dW_{\hat{\sigma}}$. Let $\ell_{\hat{\sigma}}$, $\hat{\sigma}\in\hat{\Sigma}$,  be the pure spinor of $W_{\hat{\sigma}}$ in $S_{\hat{\sigma}}$. The line in $C(V_{\hat{\sigma}})$ spanned by $\varphi_{\hat{\sigma}}(\ell_{\hat{\sigma}}\otimes \ell_{\hat{\sigma}})$
is the product in $C(V_{\hat{\sigma}})$ of the elements of a basis of $W_{\hat{\sigma}}$, by \cite[III.3.2]{chevalley}.
Hence, the isomorphism $\tilde{\varphi}$ maps the subspace spanned by $\ell_{\hat{\sigma}}\otimes \ell_{\hat{\sigma}}$, , $\hat{\sigma}\in\hat{\Sigma}$,  
to a subspace of $\wedge^*V_\RR$, whose projection to $\wedge^dV_\RR$ is equal to the subspace spanned over $\RR$ by $\wedge^d_FW$.
In particular, considering $W=H^1(\hat{X},\QQ)$, the projection of 
$\tilde{\varphi}_{\hat{\sigma}}(1\otimes 1)$ to $\wedge^{top}H^1_{\hat{\sigma}}(\hat{X})$ is non-zero. Similarly, considering $W=H^1(X,\QQ)$, the projection of 
$\tilde{\varphi}_{\hat{\sigma}}([pt]\otimes [pt])$ to $\wedge^{top}H^1_{\hat{\sigma}}(X)$ is non-zero. 
\end{rem}

Set $S_{\hat{\eta}}:=S_\QQ\cap[\sum_{\hat{\sigma}\in\hat{\Sigma}}S_{\hat{\sigma}}]$, $S^+_{\hat{\eta}}:=S^+_\QQ\cap[\sum_{\hat{\sigma}\in\hat{\Sigma}}S^+_{\hat{\sigma}}]$, and define $S^-_{\hat{\eta}}$ analogously. These are subspaces of $S_\QQ:=\wedge^*_\QQ H^1(X,\QQ)$. The latter does not have a natural structure as an $F$-vector space, but it is naturally a representation of the  group $F^\times$. 
If we let $f\in F^\times$ act on $\theta\in H^1(X,\QQ)^*$ by  $(f,\theta(\bullet))\mapsto \theta(f^{-1}(\bullet))$, then $\wedge^{top}H^1(X,\RR)^*$ has the character $\prod_{\hat{\sigma}\in\hat{\Sigma}}\hat{\theta}^{-d}$. Identify $C(V_\RR)$ with $\End(S_\RR)$ via the isomorphism $m$.
The isomorphism $\varphi:S_\RR\otimes S_\RR\rightarrow C(V_\RR)$ has weight $\prod_{\hat{\sigma}\in\hat{\Sigma}}\hat{\theta}^{-d}$, which is the lowest weight in $C(V_\RR)$. The weights in $S_{\hat{\sigma}}\otimes_\RR S_{\hat{\sigma}}$ are 
$\{\hat{\theta}^j \ : \ 0\leq j \leq d\}$ and it is characterized as the subspace of $S_\RR\otimes S_\RR$ containing all characters with these weights.
%**********
% End Hide
%**********
}

%*************
% Hide
%*************
\hide{
Set $(S\otimes S)_{\hat{\eta}}:=(S_\QQ\otimes_\QQ S_\QQ)\cap [\sum_{\hat{\sigma}\in\hat{\Sigma}}S_{\hat{\sigma}}\otimes_\RR S_{\hat{\sigma}}]$.

\begin{lem}
\label{lemma-Chevalley-isomorphism-over-F} (??? needs to be corrected using Remark \ref{rem-image-via-varphi-of-top-wedge-over-F-of-a-maximal-isotropic-F-subspace} ???)
The isomorphism 
\[
\phi\circ (id\otimes\tau):S_\QQ\otimes_\QQ S_\QQ\rightarrow \wedge^*V_\QQ,
\] 
where $\phi$ is
given in (\ref{eq-phi-introduction}), 
%is also an isomorphism of $F$-vector spaces and it 
restricts to an isomorphism of $(S\otimes S)_{\hat{\eta}}$ onto (???).
\end{lem}

\begin{proof}
The isomorphism $\phi\circ (id\otimes\tau)$ is $\Spin(V)$-equivariant, by Proposition \ref{prop-the-orlov-image-of-HW-P-projects-into-the-3-dimensional-space-of-HW-classes}, and is hence also $\Spin(V_{\hat{\eta}})$-equivariant. Hence, it maps the irreducible $\Spin(V_{\hat{\eta}})$-sub-representation
$S_{\hat{\eta}}\otimes_F S_{\hat{\eta}}$ to a $\Spin(V_{\hat{\eta}})$-sub-representation of $\wedge^*_FV_{\hat{\eta}}$. (??? show that the restriction is an isomorphism of $F$-vector spaces and that the multiplicity of $\wedge^*_FV_{\hat{\eta}}$ in $\wedge^*_\QQ V_\QQ$ is $1$ ???)
Note that $\wedge^*_FV_{\hat{\eta}}\cong \wedge^*_FH^1(X,\QQ)\otimes_F \wedge^*_FH^1(\hat{X},\QQ)$
\end{proof}

%*************
% End Hide
%*************
}

%*************************************************************
%
%*************************************************************
\subsection{The embedding of the even isotropic grassmannian in the even half-spin representation}
\label{sec-embedding-of-isotropic-grassmannian}
Given an element $\lambda\in S^+_\CC$, define $\mu_\lambda:V_\CC\rightarrow S^-_\CC$ by $\mu_\lambda(v)=m_v(\lambda)$, where $m_v$ is given in (\ref{eq-m}) and we consider $V_\CC$ as a subspace of $C(V_\CC)$.
An element $\lambda\in S^+_\CC$ is called an {\em even pure spinor}, if the kernel $W_\lambda$ of $\mu_\lambda$ is a maximal isotropic subspace \cite[III.1.4]{chevalley}. In that case, $\span\{\lambda\}=\{x\in S^+_\CC \ : \ m_w(x)=0, \forall w\in W_\lambda\}$. Furthermore, the set of maximal isotropic subspaces that arise this way consists of one of the two connected components of the grassmannian $IGr(2n,V_\CC)$ of maximal isotropic subspaces of $V_\CC$. 
The component is denoted by $IGr^+(2n,V_\CC)$ and is called the even component.
We get an embedding
\begin{equation}
\label{eq-ell}
\ell:IGr^+(2n,V_\CC)\rightarrow \PP(S^+_\CC)
\end{equation}
sending a maximal isotropic subspace $W$ of $V_\CC$ in $IGr^+(2n,V_\CC)$ to the unique point $\ell_W$ associated to the unique line in $S^+_\CC$, such that for every non-zero class $\lambda$ in $\ell_W$, the equality $W=\ker(\mu_\lambda)$ holds. 
The image  $\ell(IGr^+(2n,V_\CC))$ is called the {\em even spinorial variety}. We will refer to points of $\ell(IGr^+(2n,V_\CC))$ as even pure spinors as well.
The embedding $\ell$ is $\Spin(V_\CC)$-equivariant and it is defined over $\QQ$.

\begin{example}
\label{example-pure-spinors}
The pure spinor of the subspace $H^1(X,\QQ)$ of $V_\QQ$ is $\PP(H^{4n}(X,\QQ))$.
The pure spinor of the subspace $H^1(\hat{X},\QQ)$ of $V_\QQ$ is $\PP(H^0(X,\QQ))$. In particular, $1\in H^0(X,\QQ)\subset S^+_\QQ$ is a pure spinor and for every $g\in \Spin(V_\CC)$, the class $m_g(1)$ is a pure spinor. For each class $\theta\in H^2(X,\CC)$, there exists a unique element $g\in \Spin(V_\CC)$, such that $m_g(x)=\exp(\theta)\cup x$, for all $x\in S_\CC=H^*(X,\CC)$, and so $exp(\theta)=m_g(1)$ is a pure spinor. 
\end{example}

Given an embedding $\hat{\eta}:F\rightarrow \End(V_\QQ)$ of a totally real field $F$ of a CM-field $K$, so that $V_{\hat{\eta}}$ is $2d$-dimensional over $F$, we get the analogous embedding of the even component
$IGr^+(2d,V_{\hat{\eta}}\otimes_FK)$ of maximal isotropic $K$-subspaces of $V_{\hat{\eta}}\otimes_FK$ in the projectivization of the $K$-vector space 
$\wedge^{ev}_K[H^1(X,\QQ)\otimes_FK]$ (the even exterior algebra).
Similarly, given a character $\hat{\sigma}$ of $F$, we get the analogous embedding of the even component
$IGr^+(d,V_{\hat{\sigma},\RR})$ in $\PP(\wedge_\RR^{ev} V_{\hat{\sigma},\RR})$ as well the analogous embedding for $V_{\hat{\sigma},\RR}\otimes_\RR\CC$.

%*************************************************************
%
%*************************************************************
\section{Spin groups over $F$ and real embeddings of $F$}
\label{sec-spin-groups-over-F}
We define 
\[
\Spin(V_{\hat{\eta}})
\] 
as the Spin group associated to the $F$-valued $F$-bilinear pairing $(\bullet,\bullet)_{V_{\hat{\eta}}}$.
%We define 
%\[
%\Spin(V)_{\hat{\eta}}
%\] 
%as the subgroup of $\Spin(V_{\hat{\eta}})$ preserving the lattice $V$. 

%*************************************************************
%
%*************************************************************
\subsection{A homomorphism from $\Spin(V_{\hat{\eta}})$ to $\Spin(V_\RR)$}
\label{sec-homomorphism-from-Spein-V-hat-eta-to-Spin-V-RR}

%\begin{rem}
%\label{rem-factorization-of-Clifford-algebra}
%\begin{enumerate}
%\item 
 In the current section \ref{sec-homomorphism-from-Spein-V-hat-eta-to-Spin-V-RR} we continue to denote $V_{\hat{\sigma},\RR}$ by $V_{\hat{\sigma}}$.
The composition
\[
C(V_\QQ)\rightarrow C(V_\RR)\IsomRightArrow \otimes_{\hat{\sigma}\in\hat{\Sigma}}C(V_{\hat{\sigma}})
\]
is a $\QQ$-algebra homomorphism. 
%The field $F$ acts on the $\hat{\sigma}$-factor via the embedding $\hat{\sigma}:F\rightarrow \RR$.

%\item
%\label{rem-item-C-V-hat-eta}
Consider the Clifford algebra 
\[
C(V_{\hat{\eta}}):=\oplus_{j\geq 0} V_{\hat{\eta}}^{\otimes_F j}/\langle u\otimes_F v+ v\otimes_F u - (u,v)_{V_{\hat{\eta}}}
\rangle,
\]
where the scalar $(u,v)_{V_{\hat{\eta}}}$ is in $F$. 
The embedding
\[
(\hat{\sigma})_{\hat{\sigma}\in\hat{\Sigma}}:
V_{\hat{\eta}}\rightarrow V_\RR=\oplus_{\hat{\sigma}\in\hat{\Sigma}}V_{\hat{\sigma}}
\]
is isometric, if the pairing on the right hand side is given by $(u,v)\mapsto((u_{\hat{\sigma}},v_{\hat{\sigma}})_{V_{\hat{\sigma}}})_{\hat{\sigma}\in\Hat{\Sigma}}$
as an element in $F\otimes_\QQ\RR\cong\oplus_{\hat{\sigma}\in\Hat{\Sigma}}\RR$. 
The diagonal embedding of $V_\QQ$ in $V_\RR$ sends $v\in V_\QQ\subset C(V_\QQ)$ to $\sum_{\hat{\sigma}\in\hat{\Sigma}}\hat{\sigma}(v)\in \otimes_{\hat{\sigma}\in\hat{\Sigma}}C(V_{\hat{\sigma}})$ (see the proof of \cite[Lemma V.1.7]{lam}). It does not extends to an algebra homomorphism
$
C(V_{\hat{\eta}})\rightarrow C(V_\RR).
$
%The vector space $V_\RR$ is also a module over $F\otimes_\QQ\RR\cong\prod_{\hat{\sigma}\in\hat{\Sigma}}\RR$ and we can form the Clifford algebra 
%of $V_\RR$ over the latter ring, in which $C(V_{\hat{\eta}})$ embeds. 

We construct next a diagonal multiplicative group homomorphism from the group 
$C(V_{\hat{\eta}})^{even,\times}$, of invertible even\footnote{Note that even tensor factors  commute in the tensor product $\otimes_{\hat{\sigma}\in\hat{\Sigma}}C(V_{\hat{\sigma}})$ in the category of $\ZZ_2$-graded algebras.} elements in $C(V_{\hat{\eta}})$, into $C(V_\RR)^{even,\times}.$
We have the $F$-algebra homomorphism
\begin{equation}
\label{eq-diagonal-embedding-of-Clifford-algebra-of-F-vector-space-V}
C(V_{\hat{\eta}})\rightarrow C(V_{\hat{\eta}})\otimes_\QQ\RR\cong C(V_{\hat{\eta}})\otimes_F(F\otimes_\QQ\RR)\cong 
\prod_{\hat{\sigma}\in\hat{\Sigma}}C(V_{\hat{\sigma}}),
\end{equation}
where the left homomorphism sends $x$ to $x\otimes 1$.
For each $\hat{\sigma}\in\hat{\Sigma}$ we have an embedding $\hat{\sigma}:C(V_{\hat{\eta}})\rightarrow C(V_{\hat{\sigma}})$ obtained by composing the above homomorphism with projection to the $\hat{\sigma}$-factor.
%$\hat{\sigma}:C(V_{\hat{\eta}})\rightarrow C(V_\RR)\IsomRightArrow \otimes_{\hat{\sigma}\in\hat{\Sigma}}C(V_{\hat{\sigma}})$, 
%sending $g\in C(V_{\hat{\eta}})$ to $g$ (???) tensored with $1_{\hat{\sigma}_1}\in C(V_{\hat{\sigma}_1})$, for all $\hat{\sigma}_1\neq\hat{\sigma}$. 
We get the diagram
\begin{equation}
\label{eq-diagram-of-spin-groups}
\xymatrix{
C(V_{\hat{\eta}})^{even,\times} \ar[r]^-{(\ref{eq-diagonal-embedding-of-Clifford-algebra-of-F-vector-space-V})} &
\prod_{\hat{\sigma}\in\hat{\Sigma}}C(V_{\hat{\sigma}})^{even,\times}
\ar[r] &
\left[\otimes_{\hat{\sigma}\in\hat{\Sigma}}C(V_{\hat{\sigma}})\right]^{even,\times} \ar[r]^-\cong
&
C(V_\RR)^{even,\times}
\\
\Spin(V_{\hat{\eta}}) \ar[r] \ar[u]^{\cup} &
\prod_{\hat{\sigma}\in\hat{\Sigma}}\Spin(V_{\hat{\sigma}})
\ar[rr] \ar[u]^{\cup} & &
\Spin(V_\RR) \ar[u]^{\cup}
}
\end{equation}
where the top middle arrow sends $(g_{\hat{\sigma}})_{\hat{\sigma}\in\hat{\Sigma}}$ to $\otimes_{\hat{\sigma}\in\hat{\Sigma}}g_{\hat{\sigma}}$.
The composition 
\begin{equation}
\label{eq-homomorphism-of-even-miltiplicative-Clifford-groups}
C(V_{\hat{\eta}})^{even,\times}\rightarrow C(V_\RR)^{even,\times},
\end{equation}
of the top horizontal homomorphisms above, restricts to 
$Nm_{F/\QQ}:F^\times\rightarrow \QQ^\times\subset\RR^\times$. In particular, it is not injective.
%\item

%The subspace $\wedge^*H^1_{\hat{\sigma}}(X)$ of $S_\RR:=\wedge^*H^1(X,\RR)$ is $C(V_{\hat{\eta}})^{even,\times}$-invariant, 
%as the tensor factorization of $C(V_\RR)$ is compatible with the factorization 
%$\wedge^*H^1(X,\RR)\cong \otimes_{\hat{\sigma}\in\hat{\Sigma}}\wedge^*H^1_{\hat{\sigma}}(X)$, 
%where again the tensor product is taken in the category of $\ZZ_2$-graded algebras.
%\end{enumerate}
%\end{rem}

%***********
% Hide
%***********
\hide{
The direct sum
$\oplus_{\hat{\sigma}\in\hat{\Sigma}}V_{\hat{\sigma}}$, considered as a subspace of $V_\QQ\otimes_\QQ(F\otimes_\QQ\RR)$, is equal to the image of
$V_\QQ\otimes_\QQ F$ via
\begin{equation}
\label{eq-embedding-of-V-tensor-F}
V_\QQ\otimes_\QQ F\LongRightArrowOf{(id_{V_\QQ}\otimes \hat{\sigma})_{\hat{\sigma}\in\hat{\Sigma}}} V_\QQ\otimes_\QQ \left(\prod_{\hat{\sigma}\in\hat{\Sigma}}\RR\right)
\cong V_\QQ\otimes_\QQ (F\otimes_\QQ \RR).
%\rightarrow V_\QQ\otimes_F (F\otimes_\QQ \RR).
\end{equation}
%***********
% End  Hide
%***********
}

%*************************************************************
%
%*************************************************************
\subsection{A homomorphism from $\Spin(V_{\hat{\eta}})$ to $\Spin(V_\QQ)$}
Let $SO(V_{\hat{\eta}})$ be the subgroup of $SO(V_\QQ)$ preserving the $F$-valued pairing $(\bullet,\bullet)_{V_{\hat{\eta}}}.$ 
%Elements of $SO(V_{\hat{\eta}})$ are precisely the elements of $SO(V_\QQ)$ leaving invariant $V_{\hat{\sigma}}$, for all $\hat{\sigma}\in\hat{\Sigma}$. 
Let  $SO_+(V_{\hat{\eta}})$  be the image of $\Spin(V_{\hat{\eta}})$ in $SO(V_{\hat{\eta}})$.
Let 
\begin{equation}
\label{eq-Spin-V-QQ-hat-eta}
\Spin(V_\QQ)_{\hat{\eta}}
\end{equation}
be the image of $\Spin(V_{\hat{\eta}})$ in $\Spin(V_\QQ)$ via the homomorphism (\ref{eq-homomorphism-of-even-miltiplicative-Clifford-groups}).
Lemma \ref{lemma-Spin-V-hat-eta} below states that $\Spin(V_\QQ)_{\hat{\eta}}$ is isomorphic to either the domain or the codomain
of the degree $2$ homomorphism $\rho:\Spin(V_{\hat{\eta}})\rightarrow SO_+(V_{\hat{\eta}})$. The lemma also compute $\Spin(V_\QQ)_{\hat{\eta}}$ explicitly.

Let $\{e_{\hat{\sigma},1}, \dots, e_{\hat{\sigma},d}\}$ be an orthogonal basis of $V_{\hat{\sigma},\RR}$ satisfying $(e_{\hat{\sigma},i},e_{\hat{\sigma},i})\in\{2,-2\}$. Then $e_{\hat{\sigma},i}\in C(V_\RR)$ belongs to the Clifford group and $\rho_{e_{\hat{\sigma},i}}$ acts on $V_\RR$ as minus the reflection in the hyperplane orthogonal to $e_{\hat{\sigma},i}$. Note that $e_{\hat{\sigma},i}e_{\hat{\sigma},j}=-e_{\hat{\sigma},j}e_{\hat{\sigma},i}$, if $i\neq j$, but the images 
of $e_{\hat{\sigma},i}$ and $e_{\hat{\sigma},j}$
in $O(V_\RR)$ commute.
Set 
\begin{equation}
\label{eq-g-hat-sigma}
g_{\hat{\sigma}}:=\prod_{i=1}^de_{\hat{\sigma},i}.
\end{equation} 
Choose an embedding $\hat{\sigma}_0:F\rightarrow \RR$ and let $\tilde{F}\subset \RR$ be the Galois closure of $\hat{\sigma}_0(F)$ over $\QQ$. 
The subfield $\tilde{F}$ of $\RR$ is independent of $\hat{\sigma}_0$ and any embedding $\hat{\sigma}\in\hat{\Sigma}$ factors through an embedding of $F$ in $\tilde{F}$, which we denote by $\hat{\sigma}$ as well.
Diagram (\ref{eq-diagram-of-spin-groups}) holds, with $V_\RR$ replaced by $V_{\tilde{F}}$ and $V_{\hat{\sigma}}$ now denotes the corresponding direct summand of $V_{\tilde{F}}$.
The element $g_{\hat{\sigma}}$ belongs to the Clifford subgroup of $C(V_{\tilde{F}})^{even,\times}$, which belongs to $\Spin(V_{\tilde{F}})$ if and only if $d/2$ is even, since the signature of $V_{\hat{\sigma}}$ is $(d/2,d/2)$. Note that $g_{\hat{\sigma}}^2=(-1)^{d/2}$.
The element 
$g_{\hat{\sigma}}$ maps to $SO(V_{\tilde{F}})$ acting on $V_{\hat{\sigma}}$ via multiplication by  $-1$ and it acts as the identity on $V_{\hat{\sigma}'}$, 
for $\hat{\sigma}'\neq \hat{\sigma}$. Under the natural embedding of $C(V_{\hat{\sigma}_1})$ in $C(V_{\tilde{F}})$, the element $g_{\hat{\sigma}_2}$
commutes with $C(V_{\hat{\sigma}_1})^{ev}$ and anti-commutes with $C(V_{\hat{\sigma}_1})^{odd}$, for all $\hat{\sigma}_1$, $\hat{\sigma}_2$ in $\hat{\Sigma}$, by \cite[II.2.4]{chevalley}.

Let $SO_+(V_{\hat{\sigma}})$ be the image of $\Spin(V_{\hat{\sigma}})$ in $SO(V_{\hat{\sigma}})$. We have the short exact sequences
\begin{eqnarray*}
0\rightarrow SO_+(V_\QQ)\rightarrow &SO(V_\QQ)&\rightarrow \QQ^\times/\QQ^{\times,2}\rightarrow 0,
\\
0\rightarrow SO_+(V_{\hat{\eta}})\rightarrow &SO(V_{\hat{\eta}})&\rightarrow F^\times/F^{\times,2}\rightarrow 0,
\\
0\rightarrow SO_+(V_{\hat{\sigma}})\rightarrow &SO(V_{\hat{\sigma}})&\rightarrow \tilde{F}^\times/\tilde{F}^{\times,2}\rightarrow 0,
\end{eqnarray*}
where $\QQ^{\times,2}$, $F^{\times,2}$, and $\tilde{F}^{\times,2}$ are the groups of squares of non-zero elements 
(see \cite[II.3.7]{chevalley}). Furthermore, $SO_+(\bullet)$ is the commutator subgroup of $SO(\bullet)$, for $\bullet=V_\QQ$, $V_{\hat{\eta}}$, and $V_{\hat{\sigma}}$ \cite[II.3.9]{chevalley}. 
Hence, for every $g\in \Spin(V_\QQ)_{\hat{\eta}}$, $\rho_g$ restricts to the subspace $V_{\hat{\sigma}}$ of $V_{\tilde{F}}$ 
as an element of $SO_+(V_{\hat{\sigma}})$.

Let $\Spin(V_\QQ)_{\{g_{\hat{\sigma}}\}}$ be the subgroup 
of $\Spin(V_\QQ)$ consisting of elements, which commute in $C(V_{\tilde{F}})^{even,\times}$ with $g_{\hat{\sigma}}$, for all $\hat{\sigma}\in\hat{\Sigma}$.
Let $\Spin(V_\QQ)''_{\{g_{\hat{\sigma}}\}}$
%\footnote{A better notation would be $\Spin(V_\QQ)_{\hat{\eta}}$, but only after we change the notation in section \ref{sec-isotypic-decomposition}.} 
be the commutator subgroup of the commutator subgroup of $\Spin(V_\QQ)_{\{g_{\hat{\sigma}}\}}$ (the second derived subgroup).
%the subgroup of $\Spin(V_\QQ)$ consisting of elements $g$, which commute in $C(V_{\tilde{F}})^{even,\times}$ with $g_{\hat{\sigma}}$, 
%for all $\hat{\sigma}\in\hat{\Sigma}$, and such that $\rho_g$ restricts to the subspace $V_{\hat{\sigma}}$ of $V_{\tilde{F}}$ 
%as an element of $SO_+(V_{\hat{\sigma}})$.
The group $\Spin(V_\QQ)_{\{g_{\hat{\sigma}}\}}$ is independent of the choice of the bases $\{e_{\hat{\sigma},i}\}$,
since $g_{\hat{\sigma}}$ is determined, up to sign,  by its image $\rho_{g_{\hat{\sigma}}}$ in $SO(V_{\tilde{F}})$ and $\rho_{g_{\hat{\sigma}}}$ is independent of the choice of the bases. 

%***********
% Hide
%***********
\hide{
Following is another characterization of $\Spin(V_\QQ)_{\hat{\eta}}$. The Galois group $Gal(\tilde{F}/\QQ)$ acts on the Clifford algebra $C(V_{\tilde{F}})\cong C(V_\QQ)\otimes_\QQ\tilde{F}$ via its action on the second tensor factor. The group $\Spin(V_{\tilde{F}})$ is $Gal(\tilde{F}/\QQ)$-invariant.
The group $\Spin(V_\QQ)_{\hat{\eta}}$ is the subgroup of $\Spin(V_{\tilde{F}})$ consisting of $Gal(\tilde{F}/\QQ)$-invariant elements $g$, such that $\rho_g$
leaves invariant each summand in the decomposition $V_{\tilde{F}}=\oplus_{\hat{\sigma}\in\hat{\Sigma}}V_{\hat{\sigma}}$ and its restriction to each $V_{\hat{\sigma}}$ belongs to $SO_+(V_{\hat{\sigma}})$. 
Indeed, 
let $g$ be an element of $\Spin(V_\QQ)_{\hat{\eta}}$. For each $\hat{\sigma}\in\hat{\Sigma}$ choose an element $h_{\hat{\sigma}}\in\Spin(V_{\hat{\sigma}})$ mapping to the restriction of $\rho_{g}$ to $V_{\hat{\sigma}}$. Set $g':=\otimes_{\hat{\sigma}\in\hat{\Sigma}}h_{\hat{\sigma}}\in C(V_{\tilde{F}})$. Then $g'$ belongs to $\Spin(V_{\tilde{F}})$ and $\rho_g=\rho_{g'}$ in $SO_+(V_{\tilde{F}})$. Thus, $g'=\pm g$ and 
so $g'$ belongs to $\Spin(V_\QQ)_{\hat{\eta}}$.
%$\rho_{g'}$ belongs to $SO_+(V_\QQ)$, since $\rho_g$ belongs to $SO_+(V_\QQ)$. 
It follows that $g$ is a $Gal(\tilde{F}/\QQ)$-invariant element of $\Spin(V_{\tilde{F}})$ such that $\rho_g$ preserves the decomposition $V_{\tilde{F}}=\oplus_{\hat{\sigma}\in\hat{\Sigma}}V_{\hat{\sigma}}$. Conversely, let $g$ be a $Gal(\tilde{F}/\QQ)$-invariant element of $\Spin(V_{\tilde{F}})$, which restriction to each $V_{\hat{\sigma}}$ belongs to $SO_+(V_{\hat{\sigma}})$. Choose an element $h_{\hat{\sigma}}\in\Spin(V_{\hat{\sigma}})$ mapping to the restriction of $\rho_{g}$ to $V_{\hat{\sigma}}$. Set $g':=\otimes_{\hat{\sigma}\in\hat{\Sigma}}h_{\hat{\sigma}}\in C(V_{\tilde{F}})$. Then $g'$ belongs to $\Spin(V_{\tilde{F}})$ and $\rho_g=\rho_{g'}$ in $SO_+(V_{\tilde{F}})$. Thus, $g'=\pm g$. Now $h_{\hat{\sigma}}$ commutes with $g_{\hat{\sigma}}$, for all $\hat{\sigma}\in\hat{\Sigma}$.
Hence, so does $g$. We conclude that $g$ belongs to $\Spin(V_\QQ)_{\hat{\eta}}$. 
%If $\dim_\QQ(F)$ is odd, then $-1$ belongs to the image of $\Spin(V_{\hat{\eta}})$ and so 
%************
% End Hide
%************
}

\begin{lem}
\label{lemma-Spin-V-hat-eta}
%\begin{enumerate}
%\item
%\label{lemma-item-Spin-V-hat-eta-maps-into-Spin-V-QQ}
The  homomorphism (\ref{eq-homomorphism-of-even-miltiplicative-Clifford-groups}) maps $\Spin(V_{\hat{\eta}})$ onto 
$\Spin(V_\QQ)''_{\{g_{\hat{\sigma}}\}}$, so the latter is equal to $\Spin(V_\QQ)_{\hat{\eta}}$.
%the intersection in $C(V_\RR)^{even,\times}$ of $\Spin(V_\QQ)$. 
%with the image of the subgroup $C(V_{\hat{\eta}})^{even,\times}$ of the group of even invertible elements.  
The resulting homomorphism 
\begin{equation}
\label{eq-homomorphism-from-Spin-V-hat-eta-to-commutator-in-Spin-V-QQ}
\zeta:\Spin(V_{\hat{\eta}})\rightarrow \Spin(V_\QQ)_{\hat{\eta}}
\end{equation}
%$\Spin(V_{\hat{\eta}})\rightarrow \Spin(V_\QQ)$ 
is an isomorphism, if $\dim_\QQ(F)$ is odd, and its kernel is $-1$, if $\dim_\QQ(F)$ is even.
The outer square of the diagram
\[
\xymatrix{
\Spin(V_{\hat{\eta}})\ar[r]^{\zeta} \ar[d]_\rho&
\Spin(V_\QQ)_{\hat{\eta}}\ar[r]^{\subset} \ar[dl]^{\bar{\zeta}}& \Spin(V_\QQ) \ar[d]_\rho
\\
SO_+(V_{\hat{\eta}}) \ar[rr]_-\subset & & SO_+(V_\QQ)
}
\]
is commutative. Consequently,  the right vertical homomorphism $\rho$ factors through a homomorphism $\bar{\zeta}:\Spin(V_\QQ)_{\hat{\eta}}\rightarrow SO_+(V_{\hat{\eta}})$, which is an isomorphism, if $\dim_\QQ(F)$ is even, and is equal to $\rho\circ\zeta^{-1}$, if $\dim_\QQ(F)$ is odd.
%\item
%\label{lemma-item-Spin-V-hat-eta-commutes-with-g-hat-sigma}
%The group $\Spin(V_{\hat{\eta}})$ is mapped via (\ref{eq-diagram-of-spin-groups}) into 
%$\Spin(V_\QQ)_{\hat{\eta}}$.
%The resulting homomorphism 
%\begin{equation}
%\label{eq-homomorphism-from-Spin-V-hat-eta-to-commutator-in-Spin-V-QQ}
%\Spin(V_{\hat{\eta}})\rightarrow \Spin(V_\QQ)_{\hat{\eta}}
%\end{equation}
%is an isomorphism, if $\dim_\QQ(F)$ is odd, and its image is an index $2$ subgroup if $\dim_\QQ(F)$ is even.
%\end{enumerate}
\end{lem}

\begin{proof}
%(\ref{lemma-item-Spin-V-hat-eta-maps-into-Spin-V-QQ}) 
The homomorphism  (\ref{eq-homomorphism-of-even-miltiplicative-Clifford-groups})  factors through $C(V_{\tilde{F}})^{even,\times}$.
The group $Gal(\tilde{F}/\QQ)$ acts on $C(V_{\tilde{F}})=C(V_\QQ)\otimes_\QQ\tilde{F}$ via its action on the second tensor factor and the subset $C(V_{\tilde{F}})^{even,\times}$ is $Gal(\tilde{F}/\QQ)$-invariant. The homomorphism (\ref{eq-homomorphism-of-even-miltiplicative-Clifford-groups}) is $Gal(\tilde{F}/\QQ)$-eqvariant, as the top left homomorphism in diagram (\ref{eq-diagram-of-spin-groups}) is the diagonal homomorphism. Hence, the image of $C(V_{\hat{\eta}})^{even,\times}$   in $C(V_{\tilde{F}})^{even,\times}$, via the homomorphism (\ref{eq-homomorphism-of-even-miltiplicative-Clifford-groups}),
is contained in $C(V_\QQ)^{even,\times}$. The image of $\Spin(V_{\hat{\eta}})$ is thus contained in $\Spin(V_\QQ)$, by the commutativity of the diagram (\ref{eq-diagram-of-spin-groups}).
The last statement about the kernel follows from the equality $Nm_{F/\QQ}(-1)=(-1)^{\dim_\QQ(F)}$. 

%(\ref{lemma-item-Spin-V-hat-eta-commutes-with-g-hat-sigma}) 
The group $\Spin(V_{\hat{\eta}})$ is equal to its commutator subgroup, being a simple Lie group \cite[Def. 12.45 and Cor. 21.50]{milne-algebraic-groups}. 
The element $g_{\hat{\sigma}}$ belongs to the center of 
\[
\left[\otimes_{\hat{\sigma}\in\hat{\Sigma}}C(V_{\hat{\sigma}})\right]^{even,\times}\cap\left[\otimes_{\hat{\sigma}\in\hat{\Sigma}}C(V_{\hat{\sigma}})^{ev}\right], 
\]
by \cite[II.2.4]{chevalley}.
Hence, the group $\Spin(V_{\hat{\eta}})$ is mapped via (\ref{eq-homomorphism-of-even-miltiplicative-Clifford-groups}) into 
$\Spin(V_\QQ)''_{\{g_{\hat{\sigma}}\}}$.
%************
% Hide
%************
\hide{
Let $SO_+(V_\QQ,(\bullet,\bullet)_{V_{\hat{\eta}}})$ be the subgroup of $SO_+(V_\QQ)$ preserving the pairing $(\bullet,\bullet)_{V_{\hat{\eta}}}.$ Elements 
of $SO_+(V_\QQ,(\bullet,\bullet)_{V_{\hat{\eta}}})$ are precisely the elements of $SO_+(V_\QQ)$ leaving invariant $V_{\hat{\sigma}}$, for all $\hat{\sigma}\in\hat{\Sigma}$. The image of the homomorphism $\Spin(V_{\hat{\eta}})\rightarrow SO_+(V_\QQ)$ is\footnote{The norm homomorphism from
the Clifford group of $(\bullet,\bullet)_{V_{\hat{\eta}}}$
has values in $F^\times/F^{\times,2}$, where  $F^{\times,2}$ is the group of squares of invertible elements, 
$\Spin(V_{\hat{\eta}})$ is the intersection of the kernel of the norm homomorphism with the even Clifford group,
and  $SO_+(V_\QQ,(\bullet,\bullet)_{V_{\hat{\eta}}})$  is the image of  $\Spin(V_{\hat{\eta}})\rightarrow SO(V_\QQ,(\bullet,\bullet)_{V_{\hat{\eta}}})$.
} 
$SO_+(V_\QQ,(\bullet,\bullet)_{V_{\hat{\eta}}})$, by definition (??? the group $SO_+(V_\QQ,(\bullet,\bullet)_{V_{\hat{\eta}}})$ is defined above differently ???)). Note that $SO(V_\QQ,(\bullet,\bullet)_{V_{\hat{\eta}}})/SO_+(V_\QQ,(\bullet,\bullet)_{V_{\hat{\eta}}})$
is isomorphic to $F^\times/F^{\times,2}$, 
by \cite[II.3.7]{chevalley}, and its kernel is generated by $-1$. The homomorphism $\Spin(V_\QQ)_{\hat{\eta}}\rightarrow SO_+(V_\QQ)$ has the same image $SO_+(V_\QQ,(\bullet,\bullet)_{V_{\hat{\eta}}})$ and its kernel is generated by $-1$. The homomorphism $\Spin(V_{\hat{\eta}})\rightarrow SO_+(V_\QQ,(\bullet,\bullet)_{V_{\hat{\eta}}})$
factors through $\Spin(V_\QQ)_{\hat{\eta}}$. Hence, the homomorphism (\ref{eq-homomorphism-from-Spin-V-hat-eta-to-commutator-in-Spin-V-QQ}) is an isomorphism, if $\dim_\QQ(F)$ is odd and its image is an index $2$ subgroup if $\dim_\QQ(F)$ is even.
%************
% End Hide
%************
}

It remains to prove that the subgroup $\Spin(V_\QQ)_{\hat{\eta}}$ is the whole of $\Spin(V_\QQ)''_{\{g_{\hat{\sigma}}\}}$.
%surjectivity of (\ref{eq-homomorphism-from-Spin-V-hat-eta-to-commutator-in-Spin-V-QQ}).
The pairing $(\bullet,\bullet)_{\hat{\eta}}$ is the unique $F$-valued symmetric $F$-bilinear pairing on $V_{\hat{\eta}}$, such that $tr_{F/\QQ}\circ (\bullet,\bullet)_{\hat{\eta}}=(\bullet,\bullet)_V$. Indeed, $f\mapsto tr_{F/\QQ}\circ f :\Hom_F(V_{\hat{\eta}},F)\rightarrow \Hom_\QQ(V_\QQ,\QQ)$ is an isomorphism of $\QQ$-vector spaces and the pairing $(\bullet,\bullet)_{\hat{\eta}}$ corresponds to an element of the subspace 
$\Hom_F(V_{\hat{\eta}},\Hom_F(V_{\hat{\eta}},F))$ of $\Hom_\QQ(V_{\hat{\eta}},\Hom_F(V_{\hat{\eta}},F))$. Hence, if $g$ is an isometry of $(V_\QQ,(\bullet,\bullet)_V)$, which commutes with $\hat{\eta}(F)$, then $g$ is an $F$-linear isometry of $(V_{\hat{\eta}},(\bullet,\bullet)_{V_{\hat{\eta}}})$. Indeed, in that case $tr_{F/\QQ}(g(x),g(y))_{V_{\hat{\eta}}}=(g(x),g(y))_V=(x,y)_V$, and 
$(g(\bullet),g(\bullet))_{V_{\hat{\eta}}}$ is an $F$-valued symmetric $F$-bilinear pairing on $V_{\hat{\eta}}$.
Hence, for $g\in \Spin(V_\QQ)_{\{g_{\hat{\sigma}}\}}$, the isometry $\rho_g$ of $(V_\QQ,(\bullet,\bullet)_V)$ is also an isometry of $(\bullet,\bullet)_{V_{\hat{\eta}}}$. Thus, if $g$ belongs to $\Spin(V_\QQ)''_{\{g_{\hat{\sigma}}\}}$, then $\rho_g$ belongs to the commutator subgroup $SO_+(V_{\hat{\eta}})$ of $O(V_{\hat{\eta}})$,
where $O(V_{\hat{\eta}})$ is the subgroup of $O(V_\QQ)$ of $F$-linear automorphisms preserving the $F$-valued pairing $(\bullet,\bullet)_{V_{\hat{\eta}}}.$ 
Consequently, $\rho_g=\rho_{\tilde{g}}$, for some  $\tilde{g}\in\Spin(V_{\hat{\eta}})$, by \cite[II.3.8]{chevalley}. So $g=\pm\zeta(\tilde{g})$. 
If $\dim_\QQ(F)$ is odd, then $\zeta(-1)=-1$, and so $g$ belongs to the image of $\zeta$. 
If  $\dim_\QQ(F)$ is even, then the above argument shows that the image of $\zeta$ has index at most $2$ in the first derived subgroup of $\Spin(V_\QQ)_{\{g_{\hat{\sigma}}\}}$. If the index is $2$, then the image of $\zeta$ is contained in the second derived subgroup and is a normal subgroup of the first derived subgroup with an abelian quotient, and so the image of $\zeta$ is equal to the second derived subgroup $\Spin(V_\QQ)''_{\{g_{\hat{\sigma}}\}}$.
\end{proof}

\hide{
In other words, $\Spin(V)_{\hat{\eta}}$ is isomorphic to the subgroup of 
$\Spin(H^1_{\hat{\sigma}}(X)\oplus H^1_{\hat{\sigma}}(\hat{X}))$ preserving the lattice $V$, where $V_{\hat{\sigma}}:=H^1_{\hat{\sigma}}(X)\oplus H^1_{\hat{\sigma}}(\hat{X})$ is endowed with the obvious $\RR$-bilinear symmetric pairing. The composition 
\[
V_\QQ\rightarrow V_\QQ\otimes_\QQ F \rightarrow V_{\hat{\sigma}},
\]
of the inclusion with the projection onto the eigenspace,
should be (??? define the righthand side as a $\QQ$-vector space ???) a $\QQ$-linear isomorphism of  vector spaces 
such that the bilinear pairing on $V_\QQ$ is the composition of the $F$-bilinear (???) pairing on the right hand side with $tr_{F/\QQ}:F\rightarrow \QQ$. 
%********
% End Hide
%********
}

%********
% Hide
%********
\hide{
\begin{cor}
\begin{enumerate}
\item
\label{cor-item-odd-extension}
If $\dim_\QQ(F)$ is odd, then $\zeta:\Spin(V_{\hat{\eta}})\rightarrow \Spin(V_\QQ)_{\hat{\eta}}$ is an isomorphism and 
$\Spin(V_\QQ)_{\hat{\eta}}$ is equal to the commutator subgroup of $\Spin(V_\QQ)_{\{g_{\hat{\sigma}}\}}$.
\item
\label{cor-item-even-extension}
If $\dim_\QQ(F)$ is even, then $\rho$ restricts to an isomorphism from $\Spin(V_\QQ)_{\hat{\eta}}$ onto $SO_+(V_{\hat{\eta}})$ and the latter is the commutator subgroup of $SO(V_{\hat{\eta}})$.
\end{enumerate}
\end{cor}

\begin{proof}
The surjective homomorphism $\rho:\Spin(V_\QQ)\rightarrow SO_+(V_\QQ)$ maps $\Spin(V_\QQ)_{\{g_{\hat{\sigma}}\}}$
onto $SO_+(V_\QQ)\cap SO(V_{\hat{\eta}})$. $SO_+(V_{\hat{\eta}})$ is the commutator subgroup of $SO(V_{\hat{\eta}})$, by \cite[II.3.8]{chevalley}. 
The group $SO_+(V_\QQ)$ is equal to its commutator subgroup.
Hence, $\rho$ maps the commutator subgroup $\left(\Spin(V_\QQ)_{\{g_{\hat{\sigma}}\}}\right)'$ onto $SO_+(V_{\hat{\eta}})$.
Set $\B:=\{\pm 1\}\cap\left(\Spin(V_\QQ)_{\{g_{\hat{\sigma}}\}}\right)'$.
We get the commutative diagram 
\[
\xymatrix{
0\ar[r]& \{\pm 1\} \ar[r] \ar[dd]_{Nm_{F/\QQ}} &
\Spin(V_{\hat{\eta}}) \ar[r]^\rho \ar[d]_{\zeta} &
SO_+(V_{\hat{\eta}}) \ar[d]^{=} \ar[r] & 0
\\
& & \Spin(V)_{\hat{\eta}} \ar[d]_{\cap} \ar[r]^\rho & SO_+(V_{\hat{\eta}}) \ar[d]^{=}
\\
0\ar[r] &\B \ar[r] & \left(\Spin(V_\QQ)_{\{g_{\hat{\sigma}}\}}\right)'\ar[r]_{\rho} & SO_+(V_{\hat{\eta}}) \ar[r] & 0
}
\]
with short exact top and bottom rows. 

(\ref{cor-item-odd-extension}) If $\dim_\QQ(F)$ is odd, then $Nm_{F/\QQ}:\{\pm 1\}\rightarrow \{\pm 1\}$ is an isomorphism, and so $\zeta$ 
maps $\Spin(V_{\hat{\eta}})$ isomorphically onto $\left(\Spin(V_\QQ)_{\{g_{\hat{\sigma}}\}}\right)'$.

(\ref{cor-item-even-extension}) If $\dim_\QQ(F)$ is even, then $Nm_{F/\QQ}(-1)=\zeta(-1)=1$, and $\zeta$ is surjective, by Lemma \ref{lemma-Spin-V-hat-eta},
hence the middle $\rho$ is an isomorphism.
\end{proof}
%********
% End Hide
%********
}

%****************************************************************
% 
%****************************************************************
\section{Complex multiplication by $K$ from a maximal isotropic subspace of $H^1(X\times\hat{X},\QQ)\otimes_F K$}
\label{sec-complex-multiplication--by-CM-field-K}

%****************************************************************
% 
%****************************************************************
\subsection{Complex multiplication by the CM-field $K$}
Set $V_{\hat{\eta},K}:=H^1(X\times\hat{X},\QQ)\otimes_F K$.
We have the isomorphism
\[
\Hom_K(H^1(X,\QQ)\otimes_FK,K)\IsomRightArrow \Hom_\QQ(H^1(X,\QQ)\otimes_FK,\QQ), 
\]
sending $h\in\Hom_K(H^1(X,\QQ)\otimes_FK,K)$ to $tr_{K/\QQ}\circ h$. We get the $K$ bilinear pairing 
\begin{equation}
\label{eq-pairing-on-V-hat-eta-K}
(\bullet,\bullet)_{V_{\hat{\eta},K}}:V_{\hat{\eta},K}\otimes_K V_{\hat{\eta},K}\rightarrow K
\end{equation} 
analogous to $(\bullet,\bullet)_{V_{\hat{\eta}}}$. 
We define 
\begin{equation}
\label{eq-Spin-V-K}
\Spin(V_{\hat{\eta},K})
\end{equation}
as the associated spin group. Let
\begin{equation}
\label{eq-Spin-V-Q-hat-eta-K}
\Spin(V_{\hat{\eta},K})_\QQ
\end{equation}
be the subgroup of $\Spin(V_{\hat{\eta},K})$ leaving the subset $V_\QQ$ of $V_{\hat{\eta},K}$  invariant.
%We keep the subscript $\hat{\eta}$ to indicate that $V_K$ is $V_\QQ\otimes_{\hat{\eta}{(F)}}K$ and not $V_\QQ\otimes_\QQ K$.

\begin{rem}
The group $\Spin(V_{\hat{\eta}})$ is a subgroup of $\Spin(V_{\hat{\eta},K})_\QQ$. 
%The group $\Spin(V_{\hat{\eta},K})_\QQ$ is the $Gal(K/F)$-invariant subgroup of $\Spin(V_{\hat{\eta}}\otimes_F K)$, since .
The difference is explained by the following diagram. The Galois involution $\iota$ acts on $C(V_{\hat{\eta}}\otimes_FK)\cong C(V_{\hat{\eta}})\otimes_FK$
via its action on the second tensor factor $K$. The subset $\Spin(V_{\hat{\eta},K})$ of $C(V_{\hat{\eta}}\otimes_FK)$ is $\iota$-invariant and is thus endowed with a $Gal(K/F)$-action.
\[
\xymatrix{
0 \ar[r] & \{\pm 1\} \ar[r] \ar[d]^\cong & \Spin(V_{\hat{\eta}}) \ar[r] \ar[d]^\cong & SO(V_{\hat{\eta}}) \ar[r] \ar[d]^\cong &F^\times/F^{\times, 2} \ar[r] \ar[d] & 0
\\
0 \ar[r] & \{\pm 1\} \ar[r] \ar[d]^\cong & \Spin(V_{\hat{\eta},K})^\iota \ar[r] \ar[d]^{\cap} & SO(V_{\hat{\eta}}\otimes_F K)^\iota \ar[r] \ar[d]^= &
[K^\times/K^{\times,2}]^\iota \ar[d]^=
\\
0 \ar[r] & \{\pm 1\} \ar[r] \ar[d]^\cong &\Spin(V_{\hat{\eta},K})_\QQ \ar[r] \ar[d]^{\cap} & SO(V_{\hat{\eta}}\otimes_F K)^\iota \ar[r] \ar[d]^\cap &
[K^\times/K^{\times,2}]^\iota \ar[d]^\cap
\\
0 \ar[r] & \{\pm 1\} \ar[r] &\Spin(V_{\hat{\eta},K}) \ar[r] & SO(V_{\hat{\eta}}\otimes_F K) \ar[r] & K^\times/K^{\times,2} \ar[r] & 0
}
\]
The top and bottom rows are exact, by \cite[II.3.7]{chevalley}. The third row is exact, by definition of $\Spin(V_{\hat{\eta},K})_\QQ$. 
The second row is not exact at $SO(V_{\hat{\eta}}\otimes_F K)^\iota$, since the homomorphism $F^\times/F^{\times, 2}\rightarrow K^\times/K^{\times,2}$ is not injective. Hence, $\Spin(V_{\hat{\eta}})$ is a proper subgroup of $\Spin(V_{\hat{\eta},K})_\QQ$ of index equal to the cardinality of the kernel of $F^\times/F^{\times, 2}\rightarrow K^\times/K^{\times,2}$.
\end{rem}

Choose an element $g_0$ 
%in (\ref{eq-g}) to be an element 
of $\Spin(V_{\hat{\eta},K})$. 
%So we would get a set of isometries 
%\[
%\{g_\sigma \ : \ \sigma\in\Sigma\}.
%\]
Set 
\begin{equation}
\label{eq-W}
W:=\rho_{g_0}(H^1(\hat{X},\QQ)\otimes_FK).
\end{equation}
Its pure spinor in $\wedge^*_K [H^1(X,\QQ)\otimes_FK]\cong [\wedge^*_F [H^1(X,\QQ)]\otimes_FK]$ is $m_{g_0}(1)$ (see Example \ref{example-pure-spinors}). 
Define $\eta:K\rightarrow \End(V_{\hat{\eta},K})$, so that $\eta(t)$
acts on $W$ via its $K$-subspace structure and $\eta(t)$ acts on $\bar{W}$ via scalar multiplication by $\iota(t)$, where $\iota$ is the generator of $Gal(K/F)$. The vector space 
$V_{\hat{\eta}}=H^1(X,\QQ)\oplus H^1(X,\QQ)^*$ is an $F$-subspace of $V_{\hat{\eta},K}$ of half the dimension. 
It is equal to\footnote{Here we abuse notation and write $\bar{w}$ for $(id_{V_\QQ}\otimes \iota)(w)$, where $\iota$ is the generator of $Gal(K/F)$. 
} 
\begin{equation}
\label{eq-V-QQ-is-isomorphic-to-W-T}
V_{\hat{\eta}}=\{w+\bar{w} \ : \ w\in W\},
\end{equation}
provided 
\begin{equation}
\label{eq-W-T-and-its-conjugate-are-complementary}
W\cap \iota(W)=(0),
\end{equation} 
which we assume. The above transversality holds in our main example in Lemma \ref{lemma-I-X-hat-X-belongs-to-Omega-B} ($W$ is described in equation (\ref{eq-m-g-0-H-1-hat-X}) and $\bar{W}$ in equation (\ref{eq-bar-m-g-0-H-1-hat-X})). 
Hence, the subspace $V_{\hat{\eta}}$ is invariant under $\eta(t)$ and we get an embedding
\begin{equation}
\label{eq-eta}
\eta:K\rightarrow \End_\QQ\left[H^1(X,\QQ)\oplus H^1(X,\QQ)^*\right].
\end{equation}
Extending coefficients linearly from $\QQ$ to $\CC$ we get the embedding
\begin{equation}
\label{eq-eta-CC}
\eta:K\rightarrow \End_\CC\left[H^1(X,\CC)\oplus H^1(\hat{X},\CC)\right].
\end{equation}

\begin{lem}
\label{lemma-units-are-isometries}
The endomorphism $\eta(\iota(t))$ of $V_\QQ$ is the adjoint of $\eta(t)$ with respect to both $(\bullet,\bullet)_{V_{\hat{\eta},K}}$ and $(\bullet,\bullet)_{V_\QQ}$, for all $t\in K$. Consequently, $\eta(t)$ belongs to the image of $\Spin(V_{\hat{\eta},K})_\QQ$ in $SO(V_\QQ)$, if and only if $t\iota(t)=1$.
\end{lem}

\begin{proof}
It suffices to prove the statement for $\eta:K\rightarrow V_{\hat{\eta},K}$ and the $K$-valued bilinear pairing $(\bullet,\bullet)_{V_{\hat{\eta},K}}$, since
$(x,y)_{V_\QQ}=tr_{K/\QQ}((x,y)_{V_{\hat{\eta},K}})$, for all $x,y\in V_\QQ$. Let $x\in W$ and $y\in\bar{W}$. We have
\begin{eqnarray*}
(\eta(t)x,y)_{V_{\hat{\eta},K}}&=&(tx,y)_{V_{\hat{\eta},K}}=t(x,y)_{V_{\hat{\eta},K}}=(x,ty)_{V_{\hat{\eta},K}}=(x,\iota^2(t)y)_{V_{\hat{\eta},K}}
\\
&=& (x,\eta(\iota(t))y)_{V_{\hat{\eta},K}},
\end{eqnarray*}
where the first and last equality are by definition of $\eta$. The statement follows for all $x,y\in V_{\hat{\eta},K}$, since $W$ and $\bar{W}$ are complementary and isotropic.

If $\eta(t)$ belongs to the image of $\Spin(V_{\hat{\eta},K})$, then it is an isometry and so $\eta(t)\eta(\iota(t))=id$. Hence, $t\iota(t)=1$. Conversely, if 
$t\iota(t)=1$, then $\eta(t)$ act on the complementary isotropic subspaces $W$ and $\bar{W}$ via scalar multiplication by $t$ and $t^{-1}$ and so belongs to the image of 
$\Spin(V_{\hat{\eta},K})$. As $\eta(t)$ leaves $V_\QQ$ invariant, it belongs to the image of $\Spin(V_{\hat{\eta},K})_\QQ$.
\end{proof}

\begin{rem}
\label{rem-isometry-as-primitive-element}
Note that there exists an element $t\in K$ satisfying $K=\QQ(t)$ and $t\iota(t)=1$, by \cite[Theorem 3.7]{huybrechts-k3-book}. The statement of \cite[Theorem 3.7]{huybrechts-k3-book} assumes that the CM-field is associated to the transcendental sublattice of a lattice of $K3$-type, but the proof of the existence of such a primitive element is valid for any CM-field.
\end{rem}

%****************************************************************
% 
%****************************************************************
\subsection{$F$-invariant $2$-forms on $V_\QQ$ from totally imaginary elements of $K$}
\label{sec-F-invariant-2-forms}
Let $K_-$ be the $-1$-eigenspace of $\iota:K\rightarrow K$. Given $t\in K_-$ and $x,y\in V_\QQ$, set
\begin{eqnarray}
\label{eq-Xi-t}
\tilde{\Xi}_t(x,y)&:=& (\eta(t)x,y)_{V_{\hat{\eta}}},
\\
\nonumber
\Xi_t(x,y)&:=& (\eta(t)x,y)_V.
\end{eqnarray}
We have $\Xi_t=tr_{F/\QQ}\circ \tilde{\Xi}_t$.

\begin{cor}
\label{cor-Neron-Severi-group}
$\tilde{\Xi}_t$ is anti-symmetric and for $s\in K$ we have
\[
\tilde{\Xi}_t(\eta(s)x,\eta(s)y)=\tilde{\Xi}_t(\hat{\eta}(s\iota(s))x,y).
\]
We get the injective
homomorphism $\Xi:K_-\rightarrow \wedge^2_FV^*_\QQ\subset \wedge^2_\QQ V^*_\QQ$, given by $t\mapsto \Xi_t.$
\end{cor}

\begin{proof} We have
\[
\tilde{\Xi}_t(x,y)=(\eta(t)x,y)_{V_{\hat{\eta}}}\stackrel{\mbox{Lemma \ref{lemma-units-are-isometries}}}{=}
(x,\eta(\iota(t))y)_{V_{\hat{\eta}}}=-(\eta(t)y,x)_{V_{\hat{\eta}}}=-\tilde{\Xi}_t(y,x).
\]
Hence $\tilde{\Xi}_t$ is anti-symmetric. 
We have $\Xi_t(\eta(f)x,y)=\Xi_t(x,\eta(f)y)$, for all $f\in F$, by Lemma \ref{lemma-units-are-isometries} again.
Hence $\Xi_t$ belongs to $\wedge^2_FV^*$. 
\end{proof}

Composing the homomorphism $\Xi$ in Corollary \ref{cor-Neron-Severi-group} 
with the isomorphism $V\cong V^*$ we get an injective homomorphism denoted again by $\Xi$
\begin{equation}
\label{eq-Xi}
\Xi:K_-\hookrightarrow \wedge^2_FH^1(X\times \hat{X},\QQ).
\end{equation}
%from $K_-$ to $\wedge^2_FH^1(X\times \hat{X},\QQ)$.
Note that the Neron-Severi group of a simple abelian variety with complex multiplication by $K$ has rank\footnote{
We will see in Proposition \ref{prop-generator-for-the-invariant-subalgebra}(\ref{lemma-item-Spin-V-eta-B-invariant-2-forms})
that the image of $\Xi$ consists precisely of the $(1,1)$-classes that are invariant with respect to the special Mumford-Tate group of the generic polarized abelian variety of Weil type $(A',\eta')$ in a family deforming $(X\times\hat{X},\eta)$ as an abelian variety of Weil type. The equality $\dim_\QQ(F)=\dim_\QQ(K_-)$ 
means that the flat deformations of $\Xi(K_-)$ are equal to 
$H^{1,1}(A',\QQ)$.
} 
$\dim_\QQ(F)=\dim_\QQ(K_-)$, by \cite[Prop. 5.5.7]{BL}.

%****************************************************************
% 
%****************************************************************
\subsection{A maximal isotropic subspace $W_T$ of $V_\CC$ associated to a CM-type $T$}
\label{sec-W-T}
We have two isomorphisms $V_{\hat{\sigma},\RR}\otimes_FK\rightarrow V_{\hat{\sigma},\RR}\otimes_\RR\CC$, one via $id_{V_{\hat{\sigma},\RR}}\otimes\sigma$, the other via $id_{V_{\hat{\sigma}},\RR}\otimes\bar{\sigma}$, where $\sigma$ is an embedding restricting to $F$ as $\hat{\sigma}$. 
A CM-type $T$ thus provides an embedding of $V_\RR\otimes_FK=\oplus_{\{\hat{\sigma}\in\hat{\Sigma}\}}V_{\hat{\sigma},\RR}\otimes_FK$ into $V_\RR\otimes_\RR\CC$. Precomposing with the obvious embedding of $V_\QQ\otimes_FK$ into $V_\RR\otimes_FK$ and post composing with the isomorphism $V_\RR\otimes_\RR\CC\cong V_\QQ\otimes_\QQ\CC$ we get the embedding
\begin{equation}
\label{eq-embedding-of-V-otimes-FK-in-V-CC}
e_T:V_\QQ\otimes_F K\rightarrow V_\QQ\otimes_\QQ \CC.
\end{equation}
Set
\begin{equation}
\label{eq-W-T}
W_{T,\CC}:=\span_\RR(e_T(W)).
\end{equation}
Note that $W_{T,\CC}$ is a complex subspace, as $\span_\RR(e_T(W))=\span_\CC(e_T(W))$.

The subspace 
$V_{\hat{\sigma},\RR}\otimes_FK$ of $(V_{\hat{\eta}}\otimes_\QQ\RR)\otimes_F K$ decomposes as the direct sum of its intersections 
$W_{\hat{\sigma}}:=[V_{\hat{\sigma},\RR}\otimes_FK]\cap[W\otimes_\QQ\RR]$ and 
$\bar{W}_{\hat{\sigma}}:=[V_{\hat{\sigma},\RR}\otimes_FK]\cap[\bar{W}\otimes_\QQ\RR]$, since $W$ and $\bar{W}$ are complementary $\hat{\eta}(F)$-subspaces of $V_{\hat{\eta}}\otimes_F K$.
%Consider the isomorphism
%\begin{equation}
%\label{eq-embedding-of-V-tensor-K-into-V-tensor-C}
%\sigma_{W_T}:(V_{\hat{\eta}}\otimes_\QQ\RR)\otimes_F K\rightarrow V_{\hat{\eta}}\otimes_\QQ\CC
%\end{equation}
%given by mapping 
Set
\begin{eqnarray}
\label{W-T-hat-sigma}
W_{T,\hat{\sigma}}&:=&[id_{V_{\hat{\sigma}},\RR}\otimes T(\hat{\sigma})](W_{\hat{\sigma}}),
%=[V_{\hat{\sigma},\RR}\otimes_FK]\cap [W\otimes_\QQ\RR]
\\
\nonumber
\bar{W}_{T,\hat{\sigma}}&:=&[id_{V_{\hat{\sigma}},\RR}\otimes T(\hat{\sigma})](\bar{W}_{\hat{\sigma}}).
%[V_{\hat{\sigma},\RR}\otimes_FK]\cap [\bar{W}\otimes_\QQ\RR]
\end{eqnarray}
% injectively into $V_{\hat{\sigma},\RR}\otimes_\RR\CC$ via $id_{V_{\hat{\sigma}},\RR}\otimes\sigma$, and mapping 
% \[
% \bar{W}_{T,\hat{\sigma}}&:=&[V_{\hat{\sigma},\RR}\otimes_FK]\cap [\bar{W}\otimes_\QQ\RR]
% \] 
% injectively into $V_{\hat{\sigma},\RR}\otimes_\RR\CC$ via $id_{V_{\hat{\sigma}},\RR}\otimes\bar{\sigma}$.
% We have two $K$ actions on the domain 
 %$(V_{\hat{\eta}}\otimes_\QQ\RR)\otimes_F K$ 
 %of $\sigma_{W_T}$, the left action on the left tensor factor $(V_{\hat{\eta}}\otimes_\QQ\RR)$ via $\eta$ and the right action on the right tensor factor $K$.
 %The isomorphism $\sigma_{W_T}$ conjugates the right action of $K$ on its domain to the action 
 %$\eta:K\rightarrow \End_\CC(V_{\hat{\eta}}\otimes_\QQ\CC)$, given above.
 If $\sigma:=T(\hat{\sigma})$, $v\in V_{\hat{\sigma},\RR}$, and $t\in K$, then 
 \[
 v\otimes_F\sigma(\iota(t))=v\otimes_F\overline{\sigma(t)}=\overline{v\otimes_F\sigma(t)}.
 \] 
 Hence,
 the subspace $\bar{W}_{T,\hat{\sigma}}$ of $V_\CC$ is the complex conjugate of $W_{T,\hat{\sigma}}$.
 We have $W_{T,\CC}=\oplus_{\hat{\sigma}\in\hat{\Sigma}}W_{T,\hat{\sigma}}$ and
 $\bar{W}_{T,\CC}=\oplus_{\hat{\sigma}\in\hat{\Sigma}}\bar{W}_{T,\hat{\sigma}}=\oplus_{\hat{\sigma}\in\hat{\Sigma}}\overline{W_{T,\hat{\sigma}}}$. By the definition of the complex multiplication $\eta:K\rightarrow \End(V_\CC)$ in (\ref{eq-eta-CC}), $\eta(t)$ acts on the complex subspace $W_{T,\hat{\sigma}}$ of $V_\CC$ via multiplication by $\sigma(t):=T(\hat{\sigma})(t)$ and on $\bar{W}_{T,\hat{\sigma}}$ via multiplication by $\bar{\sigma}(t)=\overline{\sigma(t)}$.
If $T(\hat{\sigma})=\sigma$, then under the decomposition $V_\CC=\oplus_{\sigma\in \Sigma}V_\sigma$ associated to $\eta$, we have
 $V_\sigma=W_{T,\hat{\sigma}}$, and 
 $V_{\bar{\sigma}}=\bar{W}_{T,\hat{\sigma}}$. So $W_T=\oplus_{\hat{\sigma}\in\hat{\Sigma}}V_{T(\hat{\sigma})}.$
 
%*************
% Hide
%*************
\hide{
We get the commutative diagram
\[
\xymatrix{
K\ar[r]^-\eta \ar[dr]_{\eta}& \End_\QQ\left(V_{\hat{\eta}}\right) \ar[d]
\\
& \End_\CC(V_{\hat{\eta}}\otimes_\QQ\CC),
}
\]
where the vertical homomorphism sends $a$ to $a\otimes id_\CC$ 
%We need to compare $V_{\hat{\eta}}\otimes_FK$ and $V_{\hat{\eta}}\otimes_\QQ K\cong V_{\hat{\eta}}\otimes_\QQ F\otimes_F K$. 
%We have the natural surjective and injective homomorphisms
%\begin{eqnarray*}
%(id_{V_{\hat{\eta}}}\otimes tr_{F/\QQ})\otimes id_K:(V_{\hat{\eta}}\otimes_\QQ F)\otimes_F K&\rightarrow& V_{\hat{\eta}}\otimes_F K,
%\\
%(id_{V_{\hat{\eta}}}\otimes 1)\otimes id_K:V_{\hat{\eta}}\otimes_F K &\rightarrow & (V_{\hat{\eta}}\otimes_\QQ F)\otimes_F K
%\end{eqnarray*}
Hence, 
%if we define $\eta_T:K\rightarrow \End(V_{\hat{\eta},K})$, so that $\eta_T(t)$
%acts on $W_T$ via  via the CM-type $T$, i.e., $\eta_T(t)$ acts on $g(H^1_{\hat{\sigma}}(X)\otimes_FK)$ 
%via the the representative $\sigma$ in $T$ restricting to $F$ as $\hat{\sigma}$, and let it act on $\bar{W}_T$ via the complex conjugate action, then 
$V_{\hat{\eta}}$ is invariant under $\eta(K)$. 
%and we get an embedding
%\[
%\eta_T:K\rightarrow \End_\QQ\left[H^1(X,\QQ)\oplus H^1(X,\QQ)^*\right].
%\]
In other words, $H^1(X,\QQ)\oplus H^1(\hat{X},\QQ)$ is invariant under the $\eta(K)$ action via
\[
\eta:K\rightarrow \End_\CC\left[H^1(X,\CC)\oplus H^1(\hat{X},\CC)\right].
\]
%*************
% End Hide
%*************
}

%****************************************************************
% 
%****************************************************************
\section{The Galois group action on the set $\T_K$ of CM-types}
\label{sec-Galois-action-on-set-of-CM-types}
Let $\sigma_0:K\rightarrow\CC$ be an embedding. Let $\tilde{K}$ be the Galois closure of $\sigma_0(K)$ in $\CC$. The subfield $\tilde{K}$ of $\CC$ is independent of $\sigma_0$. The set of embedding of $K$ in $\tilde{K}$ are in bijection to the set $\Sigma$ of its embeddings in $\CC$, as any embedding of $K$ in $\CC$ factors through $\tilde{K}$. Define $\tilde{F}\subset \RR$ similarly. Let $\T'$ be the set of unital homomorphisms of $F$-algebras
\[
T:K\rightarrow F\otimes_\QQ\tilde{K},
\]
where the inclusion of $F$ in $F\otimes_\QQ\tilde{K}$ is via $f\mapsto f\otimes 1$. 
We have the decomposition $F\otimes_\QQ\tilde{F}\cong \oplus_{\hat{\sigma}\in\hat{\Sigma}}\hat{\sigma}$, where we denote by $\hat{\sigma}$ both a field embedding $\hat{\sigma}:F\rightarrow \tilde{F}$ and a $1$-dimensional  vector space over $\tilde{F}$ endowed with an $F$ vector space structure, such that for all $f\in F$ and $v\in \hat{\sigma}$ the equality $fv=\hat{\sigma}(f)v$ holds.
We have the isomorphisms
\[
F\otimes_\QQ\tilde{K}\cong
(F\otimes_\QQ\tilde{F})\otimes_{\tilde{F}}\tilde{K}\cong
\oplus_{\hat{\sigma}\in\hat{\Sigma}}\hat{\sigma}\otimes_{\tilde{F}}\tilde{K}.
\]
Note that the direct summand $\hat{\sigma}\otimes_{\tilde{F}}\tilde{K}$ is naturally isomorphic to $\tilde{K}$.
Composing an element $T:K\rightarrow F\otimes_\QQ\tilde{K}$ of $\T'$ with the projection onto the direct summand $\hat{\sigma}\otimes_{\tilde{F}}\tilde{K}$ we get an embedding $\sigma:K\rightarrow \tilde{K}$, which restricts to $F$ as $\hat{\sigma}$. Hence, the set $\T'$ is in natural bijection with the set $\T_K$ of CM-types.

The group $Gal(\tilde{K}/\QQ)$ acts on $F\otimes_\QQ\tilde{K}$ via its action on the right tensor factor. Hence $Gal(\tilde{K}/\QQ)$ acts on $\T'$ via 
$g(T):=(id_F\otimes g)\circ T$, for all $g\in Gal(\tilde{K}/\QQ)$ and $T\in \T'$. Let us describe the corresponding action on $\T_K$.
Let $r:\Sigma\rightarrow \hat{\Sigma}$ send $\sigma:K\rightarrow \tilde{K}$ to its restriction to $F$.
Then $\T_K=\{f:\hat{\Sigma}\rightarrow \Sigma \ : \ r\circ f=id_{\hat{\Sigma}}\}$. 
Given $g\in Gal(\tilde{K}/\QQ)$, define\footnote{The map $\fm$ from $Gal(\tilde{K}/\QQ)$ to the permutation group of $\Sigma$ is unrelated to the spin representation $m$.} 
$\fm_g:\Sigma\rightarrow\Sigma$ 
by $\fm_g(\sigma):=g\circ\sigma$. Denote by $\hat{g}\in Gal(\tilde{F}/\QQ)$ the restriction of $g$ to $\tilde{F}$. Define
$\fm_{\hat{g}}:\hat{\Sigma}\rightarrow\hat{\Sigma}$ by $\fm_{\hat{g}}(\hat{\sigma})=\hat{g}\circ\hat{\sigma}$. Then $Gal(\tilde{K}/\QQ)$ acts on $\T_K$ sending
$f\in \T_K$ to
\[
g(f) := \fm_g\circ f\circ \fm_{\hat{g}^{-1}}.
\]

Given $T\in \T'$,
we get the embedding
\begin{equation}
\label{eq-e-T-over-tilde-K}
e_T:V_\QQ\otimes_F K\RightArrowOf{id_V\otimes T} V_\QQ\otimes_F(F\otimes_\QQ\tilde{K})\cong V_\QQ\otimes_\QQ\tilde{K}.
\end{equation}
The embedding in (\ref{eq-embedding-of-V-otimes-FK-in-V-CC}) naturally factors through $e_T$ above, hence the use of the same notation. 
Given an $d$-dimensional $K$-subspace $W$ of $V_\QQ\otimes_F K$, 
let 
\[
W_T
\] 
be the subspace of $V_\QQ\otimes_\QQ\tilde{K}$ spanned over $\tilde{K}$ by $e_T(W)$. 
The set 
\[
\{
(id_V\otimes g)(W_T) \ : \ g \in Gal(\tilde{K}/\QQ)
\}
\]
determines a reduced subscheme of $Gr(2n,V_\QQ\otimes_\QQ\tilde{K})$ defined over $\QQ$. 
The set
\[
\{W_T \ : \ T\in\T'\}
\]
is the union of $Gal(\tilde{K}/\QQ)$-orbits, hence it corresponds to a reduced subscheme of $Gr(2n,V_\QQ\otimes_\QQ\tilde{K})$ defined over $\QQ$. 
Assume that $W_T$ is an even maximal isotropic subspace of $V_\QQ\otimes_\QQ\tilde{K}$, for all $T\in \T'$.
Let $\ell_T\in \PP(S^+_{\tilde{K}})$ be the pure spinor of $W_T$. 
The isomorphism between the spinorial variety in $\PP(S^+_{\tilde{K}})$ and the even component of the grassmannian of maximal isotropic subspaces of $V_\QQ\otimes_\QQ\tilde{K}$ is defined over $\QQ$. Hence, the set
$\{\ell_T \ : \ T\in\T'\}$ corresponds to a reduced subscheme defined over $\QQ$ and its span in $S^+_{\tilde{K}}$ is of the form $B\otimes_\QQ\tilde{K}$,
for some subspace $B$ of $S^+_\QQ$.

\begin{rem}
\label{rem-diagram-of-spin-groups-with-Spin-V-otimes-K}
The embedding $e_T:V_{\hat{\eta}}\otimes_FK\rightarrow V_\QQ\otimes_\QQ\tilde{K}$ extends to an injective algebra homomorphism 
$\tilde{e}_T: C(V_{\hat{\eta}}\otimes_FK)\rightarrow C(V_\QQ\otimes_\QQ\tilde{K})$ and yields an analogue of Diagram
(\ref{eq-diagram-of-spin-groups}) with $\Spin(V_{\hat{\eta}}\otimes_FK)$ replacing $\Spin(V_{\hat{\eta}})$ and $\Spin(V_\QQ\otimes_\QQ\tilde{K})$ replacing $\Spin(V_\RR)$. Note, however, that the algebras embedding
$\tilde{e}_T$ depends on the choice of a CM-type $T$.
\end{rem}
%****************************************************************
% 
%****************************************************************
\section{A linear space $B$ secant to the spinorial variety}
\label{sec-secant-space-B}
Let $W\subset V\otimes_F K$ be the maximal isotropic subspace (\ref{eq-W}). 
We have two $K$-actions on $V_\QQ\otimes_FK$, one via $\eta$ on the left factor $V_\QQ$, and the obvious one on the right factor. Each of $W$ and $\bar{W}$ is invariant with respect to both. The two actions coincide on $W$ and are $\iota$-conjugate on $\bar{W}$.
The $2$-dimensional  vector $K$-subspace $\wedge^d_KW\oplus\wedge^d_K\bar{W}$ of $\wedge^d_K(V_\QQ\otimes_FK)\cong(\wedge^d_F V_\QQ)\otimes_FK$
is $\iota$-invariant and is hence of the form $HW\otimes_FK$ for a $2$-dimensional $F$-subspace $HW$ of
$\wedge^d_FV_\QQ\subset H^d(X\times\hat{X},\QQ)$. Now $HW\otimes_FK$ is also a subspace of $\wedge^d_{\eta(K)}[V_\QQ\otimes_FK]$, and so 
$HW$ is a $1$-dimensional $\eta(K)$-subspace of $\wedge^d_{\eta(K)}V_\QQ$. The latter is one-dimensional over $K$ and so we get the equality  
\begin{equation}
\label{eq-HW-CM-field-case}
HW=\wedge^d_{\eta(K)}H^1(X\times\hat{X},\QQ).
\end{equation}
%The $\hat{\eta}(F)$-subspace $HW$ is independent of the CM-type $T$, when viewed as a subspace of $\wedge^d_{\eta(K)}V_{\hat{\eta},K}$.
% (see Remark \ref{rem-notation-T}). 
%but $T$ determines its embedding in $\wedge^dH^1(X\times\hat{X},\CC)$.
%the choice of the $\left(\prod_{\sigma\in T}\sigma\right)$-eigenspace $\wedge^d_FW_T$ in $HW$.
Given an embedding $\sigma\in\Sigma$ of $K$, denote by $V_{\sigma,\CC}$ the subspace of $V_\QQ\otimes_\QQ\CC$ on which the $\eta$ action of $K$ acts via the character $\sigma:K\rightarrow\CC$. 
We have the equality
\begin{eqnarray}
%HW\otimes_\QQ\CC&=&\oplus_{\hat{\sigma}\in\hat{\Sigma}} \left(
%\wedge^d_\CC[W\otimes_\QQ\RR]_{\hat{\sigma}} \oplus \wedge^d_\CC[\bar{W}\otimes_\QQ\RR]_{\hat{\sigma}} 
%\right),
%\\
\label{eq-direct-sum-decomposition-of-HW}
HW\otimes_\QQ\CC&=&
\oplus_{\sigma\in\Sigma}\wedge^d_\CC V_{\sigma,\CC}.
%\oplus_{\sigma\in T} \left(
%\wedge^d_\CC[W_{T,\CC}]_\sigma \oplus \wedge^d_\CC[\bar{W}_{T,\CC}]_{\bar{\sigma}}
%\right),
\end{eqnarray}
%for every CM-type $T$, where the characters are of the $\eta$ action of $K$
%and $W_{T,\CC}$ is given in (\ref{eq-W-T}).
Our goal in this section is to show that if $B$ is spanned by Hodge classes, then
$HW$ is spanned by Hodge classes (Lemma \ref{lemma-if_B-is-spanned-by-Hodge-classes-so-does-HW}). In Section \ref{sec-pure-spinors-CM-case} we show that $W_T$ is a maximal isotropic subspace of $V_\CC$, for each CM-type $T$.
In Section \ref{sec-isotypic-decomposition} we introduce the groups $\Spin(V_\QQ)_\eta$ and $\Spin(V_\QQ)_{\eta,B}$, the latter being  
%playing the roles of the Mumford-Tate group and 
the special Mumford-Tate group for a generic deformation of $(X\times\hat{X},\eta)$. We relate that characters $\ell_T$ and
$\wedge^d V_\sigma$ of $\Spin(V_\QQ)_\eta$ and use it to conclude that the lines $\{\ell_T \ : T\in\T_K\}$ are linearly independent.
In Section \ref{sec-Hermitian-form-H-t} we construct a $\Spin(V_\QQ)_\eta$-invariant hermitian form. In Section \ref{sec-if-B-is-Hodge-so-in-HW}
we prove Lemma \ref{lemma-if_B-is-spanned-by-Hodge-classes-so-does-HW}.

%Compute the intersection points of $\PP(P_{T,\CC})$ with the even spinorial variety. It contains 
%$\ell_{W_T}$ and $\ell_{\bar{W}_T}$ 
%but should contain a point for other CM-types for $K$.

%*************
% Hide
%*************
\hide{
%****************************************************************
% 
%****************************************************************
\section{New strategy - pure $F$-spinors}
Let $P\subset S^+_\QQ$ be the $e$-dimensional $\QQ$-subspace, such that $P_\RR$ is spanned by the pure spinors $\ell_{\hat{\sigma}}$ of 
$W_{\hat{\sigma}}$, $\hat{\sigma}\in\hat{\Sigma}$, and their conjugates $\bar{\ell}_{\hat{\sigma}}$. Note that $\ell_{\hat{\sigma}_0}$ is the pure spinor of the maximal isotropic subspace $W_{\ell_{\hat{\sigma}_0}}\oplus_{\ell_{\hat{\sigma}}\neq \ell_{\hat{\sigma}_0}}H^1_{\ell_{\hat{\sigma}}}(\hat{X})$.
Denote $\ell_{\hat{\sigma}}$ by $\ell_\sigma$, if $\sigma\in T$, and set $\ell_\sigma:=\bar{\ell}_{\hat{\sigma}}$, if $\sigma\in \iota(T)$.
%Let $P\subset S^+_{\hat{\eta}}$ be the $2$-dimensional $F$-subspace 
%corresponding to the $K$-subspace spanned by the two pure spinors $m_{g_0}(1)$ and its conjugate $\iota(m_{g_0}(1))$
%in $S^+_{\hat{\eta}}\otimes_F K$.
%Set $\ell_{W_T}:=\span_F\{m_{g_0}(1)\}$, considered as a point in the projective space $\PP(S^+_{\hat{\eta}}\otimes_F K)$ over $K$ and 
%let $\ell_{\bar{W}_T}$ be its $\iota$-conjugate.

Use Remark \ref{rem-image-via-varphi-of-top-wedge-over-F-of-a-maximal-isotropic-F-subspace} to relate the image via Orlov's equivalence of the subspace $P\otimes_\QQ P$ of $S^+_\QQ \otimes_\QQ S^+_\QQ$ (namely, of the subspace spanned by $\ell_{\sigma_1}\otimes\ell_{\sigma_2}$, $\sigma_1,\sigma_2\in\Sigma$) to 
the subspace $HW$ of $\wedge^d_FH^1(X\times\hat{X},\QQ)$. Now choose $k$-secant sheaves $F_i$ with $ch(F_i)$ in $P$.
(??? not helpful ???)

%{\bf Problem:} $P$ is contained in the subalgebra $S^+_{\hat{\eta}}$, which is contained in the subspace $\oplus_{k=0}^{2d}\wedge^k V_\QQ$, so 
%$ch_k(F_i)$ vanishes for $2d<k\leq 4n=de$. We get a constraint, when $e>2$.
%*************
% End Hide
%*************
}

%****************************************************************
% 
%****************************************************************
\subsection{Pure spinors}
\label{sec-pure-spinors-CM-case}
%These subspaces do not depend on the initial choice of CM-type $T$. 
%Define $V_{\hat{\sigma}}\subset V_\QQ\otimes_\QQ \RR$ similarly via $\hat{\eta}$.
Let
\begin{equation}
\label{eq-V-sigma}
V_\sigma
\end{equation}
be the direct summand of $V_\QQ\otimes_\QQ\tilde{K}$ on which $\eta(K)$ acts via the character $\sigma$, where the algebraic closure $\tilde{K}$ of $K$ in $\CC$ is introduced in section \ref{sec-Galois-action-on-set-of-CM-types}.
If $\sigma$ belongs to a CM-type $T$, then 
\[V_{\sigma,\CC}=W_{T,\hat{\sigma}}
%W_{T,\CC}\cap (V_{\hat{\sigma}}\otimes_\RR\CC) 
\ \  \mbox{and} \ \ V_{\bar{\sigma},\CC}=\bar{W}_{T,\hat{\sigma}},
%\bar{W}_{T,\CC}\cap (V_{\hat{\sigma}}\otimes_\RR\CC).
\]
where $W_{T,\hat{\sigma}}$ is given in (\ref{W-T-hat-sigma}).
We have the direct sum decomposition $V_\CC=\oplus_{\sigma\in\Sigma}V_{\sigma,\CC}$
and
\begin{equation}
\label{eq-decomposition-of-W-T-CC}
W_{T,\CC}=\oplus_{\{\hat{\sigma}\in\hat{\Sigma}\}}V_{T(\hat{\sigma}),\CC}.
\end{equation} 
%Set
%\begin{equation}
%\label{eq-V-sigma-CC}
%V_\sigma:=V_{\sigma,\CC}\cap [V_\QQ\otimes_\QQ \sigma(K)]\subset V_\QQ\otimes_\QQ\CC.
%\end{equation}
%*************
% Hide
%*************
\hide{
The composition $V_\QQ\rightarrow V_\CC\rightarrow V_{\sigma,\CC}$, of the natural embedding with the projection on the direct summand, 
maps $V_\QQ$ isomorphically onto $V_\sigma$. The direct sum $\oplus_{\sigma\in\Sigma}V_\sigma$, considered as a $\QQ$-subspace of $V_\QQ\otimes_\QQ(K\otimes_\QQ\CC)$,
is isomorphic to the image of $V_\QQ\otimes_\QQ K$ via the analogue of (\ref{eq-embedding-of-V-tensor-F}).
%*************
% End Hide
%*************
}

\begin{rem}
If $K$ is a Galois extension of $\QQ$, then $Gal(K/\QQ)$ acts transitively on the set of characters via $\sigma\mapsto \sigma\circ g^{-1}$, for all $g\in Gal(K/\QQ)$.
In that case the subfield $\sigma(K)\subset \CC$ is independent of $\sigma$ and we may identify $K$ with this subfield by choosing an embedding $\sigma_0$.
The subspace $V_\sigma$ is then the subspace of $V_\QQ\otimes_\QQ \sigma_0(K)$, where $\eta(a)$ acts via $id_{V_\QQ}\otimes \sigma(a)$, for all $a\in K$. Similarly, the subspaces $V_{\hat{\sigma}}$ are the analogous subspaces of $V_\QQ\otimes_\QQ \hat{\sigma}_0(F)$.
\end{rem}

%*************
% Hide
%*************
\hide{
For each CM-type $T'$ we define the subspace $\tilde{W}_{T'}$ of $V_\QQ\otimes_\QQ K$ (not $V_{\hat{\eta},K}$) as the
%the image via $g$  in (\ref{eq-g})  of the 
direct sum  
\[
\tilde{W}_{T'}:= \oplus_{\sigma\in T'}V_\sigma.
\] 
We define the subspace $W_{T'}$ of $V_{\hat{\eta},K}$ as the image of $\tilde{W}_{T'}$ via the natural homomorphism 
\begin{equation}
\label{eq-mapping-tilde-W-T-to-W-T}
V_\QQ\otimes_\QQ K
%\RightArrowOf{id_{V_\QQ}\otimes tr_{K/\QQ}} 
\rightarrow V_\QQ\otimes_F K=:V_{\hat{\eta},K}
\end{equation}
mapping $v\otimes(f_1+f_2\sqrt{-q})$ to $\hat{\eta}_{f_1}(v)+\hat{\eta}_{f_2}(v)\otimes\sqrt{-q}$, for all $v\in V_\QQ$ and $f_1, f_2\in F$.
(??? Check that for the originally chosen CM-type $T$ we do get the previously defined $W_T$. Note that $\dim_\QQ(V_\sigma)=\dim_\QQ(V_\QQ)=\dim_\QQ(W_T)$. We may need to change the definition and consider instead the image of $V_\QQ$ via its diagonal embedding into $\tilde{W}_{T'}$ and map it to $V_{\hat{\eta},K}$ via (\ref{eq-mapping-tilde-W-T-to-W-T})???)
(??? 
%When $T'=T$ we do not get $W_T$ so the notation is confusing. Use $\tilde{W}_{T'}$ instead. 
%The subspace $W_T$ is the image of $\tilde{W}_T$ via
%\begin{equation}
%\label{eq-tr-mapping-tilde-W-T-to-W-T}
%V_\QQ\otimes_\QQ K\cong V_\QQ\otimes_\QQ F\otimes_F K\LongRightArrowOf{(id_{V_\QQ}\otimes tr_{F/\QQ})\otimes id_K} V_\QQ\otimes_F K.
%\end{equation}
%What is the relationship between $W_T$ and $\tilde{W}_T$? 
Is $\tilde{W}_T$ isotropic? If indeed it is, what is its pure spinor and how is it related to that of $W_T$?   ???)
%What is the analogue $W_{T'}$ of $W_T$ for another CM-type $T'$? It can 
%*************
% End Hide
%*************
}

For each CM-type $T$ 
consider the subspace $W_{T,\CC}$ of
%Can $W_{T}$ be defined in 
$V_\QQ\otimes_\QQ\CC$.
%as

\begin{lem}
\label{lemma-W-T-is-maximal-isotropic}
The subspace $W_{T,\CC}$ of $V_\CC$ is maximal isotropic, for every CM-type $T$.
\end{lem}

\begin{proof}
If $\hat{\sigma}_1\neq \hat{\sigma}_2$, then $V_{\sigma_1,\CC}$ and $V_{\sigma_2,\CC}$ are orthogonal, since $V_{\hat{\sigma}_1}$ and $V_{\hat{\sigma}_2}$ are, by Lemma \ref{lemma-orthogonal-direct-sum-decomposition-V-hat-sigma}.
Ditto for $V_{\bar{\sigma}_1,\CC}$ and $V_{\sigma_2,\CC}$ and for  $V_{\bar{\sigma}_1,\CC}$ and $V_{\bar{\sigma}_2,\CC}$. The subspace 
$V_{\sigma,\CC}$ is isotropic, since $W_{T,\CC}$ is, and $V_{\bar{\sigma},\CC}$ is isotropic, since $\bar{W}_{T,\CC}$ is.
Hence, $W_{T,\CC}$ is  a maximal isotropic subspace of 
$V_\QQ\otimes_\QQ\CC$
%$V_{\hat{\eta},K}:=V_\QQ\otimes_F K$, 
for every CM-type $T$.  
\end{proof}

%**********
% Hide
%**********
\hide{
We set (??? these lines are not needed ???)
\[
P_{T'}:=\ell_{W_{T',\CC}}+\ell_{\bar{W}_{T',\CC}}.
\]
If $T'\neq T$, then $P_{T'}$ is a $2$-dimensional $\RR$-subspace of $H^{ev}(X,\RR)$. For our initially chosen CM-type $T$, $P_T$ is a $2$-dimensional $\QQ$-subspace (??? NO, it need not be defined over $\QQ$ ???) of $H^{ev}(X,\QQ)$. If $W_{T,\CC}$ is the image of $H^1(\hat{X},\QQ)\otimes_F K$ via the isometry $g_0$ in Example \ref{example-isometry-g-0}, then 
\[
P_T=\span_\QQ\left\{
\left(1-\hat{\eta}_q(\Theta\wedge_F\Theta)+\dots\right),
\left(\Theta-\hat{\eta}_q(\Theta\wedge_F\Theta\wedge_F\Theta)+\dots\right)
%\left(
%1-\hat{\eta}_q(\Theta^2/2)+\hat{\eta}_q^2(\Theta^4/4!)+\dots
%\right),
%\left(
%\Theta-\hat{\eta}_q(\Theta^3/3!)+\hat{\eta}_q^2(\Theta^5/5!)+\dots
%\right)
\right\}
\]
(see Equation (\ref{eq-pure-spinor-rho-g-0-of-1})).

%**********
% End Hide
%**********
}

Denote by $\T_K$ the set of $2^{e/2}$ complex multiplication types, where $e=[K:\QQ]$. 
The element $\iota$ of $Gal(K/\QQ)$ induces an involution of $\T_K$. Let $\ell_T\subset H^{ev}(X,\CC)$ be the one-dimensional subspace spanned by a pure spinor corresponding to the subspace $W_{T,\CC}$, $T\in\T_K$.

\begin{lem}
\label{lemma-B-is-rational}
The linear subspace of $H^{ev}(X,\CC)$ spanned by $\{\ell_T \ : \T\in\T_K\}$, is of the form $B\otimes_\QQ\CC$, where $B$ is a rational subspace 
\[
B\subset H^{ev}(X,\QQ).
\] 
\end{lem}

\begin{proof}
The statement follows from Lemma \ref{lemma-W-T-is-maximal-isotropic} and the discussion in Section \ref{sec-Galois-action-on-set-of-CM-types}.
\end{proof}

Let $\ell_\sigma\in S^+_{\hat{\sigma}}\otimes_\RR\CC$ be the pure spinor of the maximal isotropic subspace $V_\sigma$ of $V_{\hat{\sigma}}$ (see Section \ref{sec-embedding-of-isotropic-grassmannian}).
Under the tensor decomposition (\ref{eq-tensor-product-decomposition-of-pure-spinor}) of pure spinors of $F$-invariant maximal isotropic subspaces of $V_\CC$,
we have
\begin{equation}
\label{eq-tensor-factorization-of-ell-T}
\ell_T=\otimes_{\hat{\sigma}\in\hat{\Sigma}}\ell_{T(\hat{\sigma})},
\end{equation}
by the equality (\ref{eq-decomposition-of-W-T-CC}).
Denote by $P_{\hat{\sigma}}\subset S^+_{\hat{\sigma}}$ the secant plane spanned by $\ell_{\sigma}$ and $\ell_{\bar{\sigma}}$.
Then 
\begin{equation}
\label{eq-tensor-factorization-of-B}
B_\RR=\otimes_{\hat{\sigma}\in\hat{\Sigma}}P_{\hat{\sigma}}.
\end{equation}

%**********
% Hide
%**********
\hide{
\begin{rem}
Let $\tilde{K}$ be the (???)  Galois closure of $K$ over $\QQ$ and $H\subset Gal(\tilde{K}/\QQ)$ the subgroup leaving all elements of $K$ invariant. 
Let $\tilde{\Sigma}$ be the set of complex embeddings of $\tilde{K}$, so that $\Sigma=\tilde{\sigma}/H$. Then $Gal(\tilde{K}/\QQ)$ acts transitively on $\Sigma$, by
$\sigma\mapsto \sigma\circ g^{-1}$, $g\in Gal(\tilde{K}/\QQ)$. The stabilizer of each element of $\Sigma$ is $H$. 
The subfield $\tilde{K}':=\tilde{\sigma}(\tilde{K})$ of $\CC$ is independent of $\tilde{\sigma}\in\tilde{\Sigma}$, since $Gal(\tilde{K}/\QQ)$ acts transitively on $\tilde{\Sigma}$.
%Let $\tilde{\sigma}\in \tilde{\Sigma}$ and $g\in Gal(\tilde{K}/\QQ)$. 
The subspace $V_{\sigma,\CC}$ is defined over $\sigma(K)$ (and hence over $\tilde{K}'$), i.e., is associated to the eigenspace $[V_\QQ\otimes_\QQ \sigma(K)]_{\sigma(t)}$ of $\eta(t)$.
Hence, the subspace $W_{T',\CC}$ is defined over $\tilde{K}'$, for all $T'\in \T_K$. It follows that each line
$\wedge^{2n}_\CC W_{T',\CC}$ is defined over $\tilde{K}'$. The group $Gal(\tilde{K}/\QQ)$ acts on $V_\QQ\otimes_\QQ \tilde{K}'$ and
on $(\wedge^{2n}_\QQ V_\QQ)\otimes \tilde{K}'\cong \wedge^{2n}_{\tilde{K}'}[V_\QQ\otimes_\QQ\tilde{K}']$. 

The subspace $\wedge^{2n}W_{T',\CC}$ of $\wedge^{2n}V_\CC$ is characterized by the property that for $a\in K$, $\wedge^{2n}\eta(a)$ acts on $\wedge^{2n}W_{T',\CC}$ via $\prod_{\sigma\in T'}(\sigma(a))^d$. The set of characters
$\{\prod_{\sigma\in T'}(\sigma)^d \ : \ T'\in \T_K\}$ is $Gal(\tilde{K}/\QQ)$-invariant. Hence, the subspace
\[
\oplus_{T'\in\T_K}\wedge^{2n}W_{T',\CC}
\]
is defined over $\QQ$, as is the $0$-subscheme of $\PP(\wedge^{2n}V_\CC)$ of points
\[
\{\PP(\wedge^{2n}W_{T',\CC}) \ : \ T'\in\T_K\}.
\]
It follows that the $0$-dimensional subscheme $\{\ell_{W_{T',\CC}}\otimes\ell_{W_{T',\CC}}  \ : \ T'\in\T_K\}$
of $H^{ev}(X,\CC)\otimes H^{ev}(X,\CC)$, consisting of tensor squares of pure spinors, is defined over $\QQ$. 

We need to work with $K$-valued characters of the subgroup $\Spin(V_{\hat{\eta},K})_{\ell_\bullet}$ of $\Spin(V_{\hat{\eta},K})$ fixing every point the subset $\{\ell_{W_{T',\CC}}  \ : \ T'\in\T_K\}$. (???)
%It suffices to show that $\ell_{W_{T',\CC}}$ is the character $\prod_{\sigma\in T'}(\sigma)^{d/2}$. (???)

Let $t\in K$ be a primitive element satisfying $t\iota(t)=1$, as in Remark \ref{rem-isometry-as-primitive-element}.
There exists an element $\tilde{t}\in \Spin(V_{\hat{\eta},K})_\QQ$, such that $\rho_{\tilde{t}}=\eta(t)$, by Lemma \ref{lemma-units-are-isometries}.
The subspaces $V_{\sigma,\CC}$, $\sigma\in\Sigma$, of $V_\CC$, given in (\ref{eq-V-sigma-CC}), are the the eigenspaces of $\rho_{\tilde{t}}$ with eigenvalues 
$\sigma(t)$, since $t$ is a primitive element of $K$.
(???)
\EndProof
\end{rem}
%**********
% End Hide
%**********
}

%*******************************************************************************************************************
%
%*******************************************************************************************************************
\subsection{An isotypic decomposition of the $\Spin(V_\QQ)_\eta$-invariant space $B\otimes B$}
\label{sec-isotypic-decomposition}
Let 
\begin{equation}
\label{eq-Spin-V-eta}
\Spin(V_\QQ)_\eta
\end{equation}
be the subgroup of the group $\Spin(V_\QQ)_{\hat{\eta}}$, given in (\ref{eq-Spin-V-QQ-hat-eta}),
leaving invariant the subspace $V_\sigma$, for all $\sigma\in\Sigma$.
%, commuting with all elements of the subspace $\eta(K)$ of $\End_\QQ(V_\QQ)$. 
Let $\Spin(V)_\eta$ be the analogous subgroup of $\Spin(V)_{\hat{\eta}}$. We will see that $\Spin(V)_\eta/\{\pm 1\}$ is an arithmetic subgroup of a product of $e/2$ copies of the unitary group $U(d,d)$ (see Equation 
(\ref{eq-factorization-of-Spin-V-tilde-K-eta-mod-pm1}) and Lemma \ref{lemma-H-t-is-hermitian} below).

Set $V_{\tilde{K}}:=V_\QQ\otimes_\QQ\tilde{K}$. Let  $\Spin(V_{\tilde{K}})_{\{g_{\hat{\sigma}}\}}$ be the subgroup of $\Spin(V_{\tilde{K}})$ consisting of elements $g$ commuting with $g_{\hat{\sigma}}$, given in (\ref{eq-g-hat-sigma}), for all $\hat{\sigma}\in\hat{\Sigma}$. Let 
$\Spin(V_{\tilde{K}})_{\hat{\eta}}$ be the commutator subgroup of $\Spin(V_{\tilde{K}})_{\{g_{\hat{\sigma}}\}}$.
We have the short exact sequence
\[
0\rightarrow \{\pm 1\}\rightarrow \Spin(V_{\tilde{K}})_{\hat{\eta}}\rightarrow \prod_{\hat{\sigma}\in\hat{\Sigma}}SO_+(V_{\hat{\sigma}})\rightarrow 0.
\]
The group $\Spin(V_{\tilde{K}})_{\hat{\eta}}$ is equal to its commutator subgroup, since that property holds for each $\Spin(V_{\hat{\sigma}}\otimes_{\tilde{F}}\tilde{K})$ and the images of the latter via the natural inclusion $C(V_{\hat{\sigma}}\otimes_{\tilde{F}}\tilde{K})\subset C(V_{\tilde{K}})$ generate $\Spin(V_{\tilde{K}})_{\hat{\eta}}$.

Let  
\begin{equation}
\label{eq-Spin-V-tilde-K-eta}
\Spin(V_{\tilde{K}})_\eta
\end{equation}
be the subgroup of $\Spin(V_{\tilde{K}})_{\hat{\eta}}$ consisting of elements $g$ 
%be the subgroup of $\Spin(V_{\tilde{K}})_{\hat{\eta}}$ consisting of elements $g$ 
%commuting with $g_{\hat{\sigma}}$, given in (\ref{eq-g-hat-sigma}), for all $\hat{\sigma}\in\hat{\Sigma}$, such that $\rho_g$ 
%restricts to $V_{\sigma}\oplus V_{\bar{\sigma}}$ as an element of $SO_+(V_{\sigma}\oplus V_{\bar{\sigma}})$, for all $\sigma\in\Sigma$, and 
such that $V_\sigma$ is $\rho_g$-invariant, for all $\sigma\in\Sigma$. Note that $\Spin(V_\QQ)_\eta$ is contained in the $Gal(\tilde{K}/\QQ)$-invariant subgroup of $\Spin(V_{\tilde{K}})_\eta$, where the $Gal(\tilde{K}/\QQ)$-action is induced by that on $C(V_{\tilde{K}})\cong C(V_\QQ)\otimes_\QQ\tilde{K}$ via the action on the second tensor factor.

Let 
\begin{equation}
\label{eq-Spin-V-eta-B}
\Spin(V)_{\eta,B}
\end{equation}
be the subgroup of $\Spin(V)_{\eta}$ fixing every point of the subspace $B$ of $H^{ev}(X,\QQ)$. Denote by
$\Spin(V_{\tilde{K}})_{\eta,B}$ the analogous subgroup of $\Spin(V_{\tilde{K}})_\eta$.
When $F=\QQ$, the group $\Spin(V_{\tilde{K}})_\eta$ specializes to the group $\Spin(V_K)_{\ell_1,\ell_2}$ in \cite[Eq. (2.2.1)]{markman-sixfolds}
and the group $\Spin(V)_{\eta,B}$ specializes to $\Spin(V)_P$, given in \cite[Eq. (2.2.2)]{markman-sixfolds}.

The group $\Spin(V_{\tilde{K}})_\eta$ leaves invariant each of the subspaces $W_T$, $T\in \T_K$, by Equation (\ref{eq-decomposition-of-W-T-CC}).
Hence, each of the lines $\ell_T$ in $H^{ev}(X,\tilde{K})$ is $\Spin(V_{\tilde{K}})_\eta$-invariant. Thus, the commutator subgroup of $\Spin(V_{\tilde{K}})_\eta$
is contained in $\Spin(V_{\tilde{K}})_{\eta,B}$. 

\[
\xymatrix{
\Spin(V)_{\hat{\eta}} \ar[r]^\subset & \Spin(V_\QQ)_{\hat{\eta}} &
\Spin(V_{\hat{\eta}})\ar[l]_{(\ref{eq-homomorphism-from-Spin-V-hat-eta-to-commutator-in-Spin-V-QQ})} \ar[r]^-\subset&
\Spin(V_{\hat{\eta},K})_\QQ \ar[r]^-\subset & \Spin(V_{\hat{\eta}}\otimes_FK)
\\
\Spin(V)_\eta \ar[r]^\subset \ar[u]^\cup & \Spin(V_\QQ)_\eta \ar[r]^\subset \ar[u]^\cup & \Spin(V_{\tilde{K}})_\eta 
\\
\Spin(V)_{\eta,B}\ar[r]_\subset \ar[u]^\cup & 
\Spin(V_\QQ)_{\eta,B}\ar[r]_\subset\ar[u]^\cup &
\Spin(V_{\tilde{K}})_{\eta,B}\ar[u]^\cup.
}
\]
%\begin{question}
%Formulate the analogue of Lemma \ref{lemma-Spin-V-hat-eta} relating $\Spin(V_{\hat{\eta},K})_\QQ$, given in (\ref{eq-Spin-V-Q-hat-eta-K}), 
%to $\Spin(V_\QQ)_{\eta,B}$. The latter should preserve also a $K$-valued hermitian form. We know that $\Spin(V_{\hat{\eta},K})_\QQ$ 
%commutes with the scalar multiplication action of $K$ on $V_\QQ\otimes_F K$, but this action is different from the $\eta(K)$ action. 
%The $\QQ$-isomorphism from $W$ to $V_\QQ$, sending $w\in W$ to $w+\iota(w)$, is equivariant with respect to the $\eta(K)$-action on both.
%\end{question}
Note that $W_T\cap W_{T'}=V_\sigma$, if $\sigma$ is the only common value of $T$ and $T'$. Hence,
$\Spin(V)_{\eta,B}$ is equal to the subgroup of $\Spin(V)_{\hat{\eta}}$ fixing every point of $B$.

For each $\sigma\in\Sigma$ let $\Spin(V_\sigma\oplus V_{\bar{\sigma}})_\eta$ be the subgroup of $\Spin(V_\sigma\oplus V_{\bar{\sigma}})$
leaving each of $V_\sigma$ and $V_{\bar{\sigma}}$ invariant.
We have a natural homomorphism 
%from the product
%$
%\prod_{\hat{\sigma}\in\hat{\Sigma}}\Spin(V_\sigma\oplus V_{\bar{\sigma}})_{\eta,P_{\hat{\sigma}}}
%$
%to $\Spin(V_{\tilde{K}})_{\eta,B}$. 
\begin{equation}
\label{eq-factorization-of-Spin-V-tilde-K-B}
\prod_{\hat{\sigma}\in\hat{\Sigma}}\Spin(V_\sigma\oplus V_{\bar{\sigma}})_{\eta,P_{\hat{\sigma}}}
\rightarrow \Spin(V_{\tilde{K}})_{\eta,B}
\end{equation}
induced by the isomorphism
\[
C(V_{\tilde{K}})\cong\otimes_{\hat{\sigma}\in\hat{\Sigma}}C(V_\sigma\oplus V_{\bar{\sigma}}),
\] 
where the tensor product is in the category of $\ZZ_2$-graded algebras, as in (\ref{eq-diagonal-embedding-of-Clifford-algebra-of-F-vector-space-V}). Above $\sigma$ restricts to $F$ as $\hat{\sigma}$. 
Now, $\Spin(V_\sigma\oplus V_{\bar{\sigma}})_{\eta,P_{\hat{\sigma}}}$ is isomorphic to the group $SL(V_\sigma)$ of linear transformations of $V_\sigma$, as a $d$-dimensional vector space over $\tilde{K}$, of determinant $1$. An elements $g$ of $SL(V_\sigma)$ acts on $V_{\bar{\sigma}}$ by $(g^*)^{-1}$, under the identification of $V_{\bar{\sigma}}$ with $V_\sigma^*$ via the bilinear pairing $(\bullet,\bullet)_V$. The proof is identical to that of 
\cite[Lemma 2.2.2]{markman-sixfolds}.
%\ref{lemma-Spin-V-K-is-SL-n-K}.

Similarly, we have the homomorphism  
\[
\prod_{\hat{\sigma}\in\hat{\Sigma}}\Spin(V_\sigma\oplus V_{\bar{\sigma}})_\eta
\rightarrow
\Spin(V_{\tilde{K}})_\eta ,
\]
and $\Spin(V_\sigma\oplus V_{\bar{\sigma}})_\eta/\{\pm 1\}$ is isomorphic to $GL(V_\sigma)$, by the proof of \cite[Lemma 1]{igusa}. An elements $g$ of $GL(V_\sigma)$ acts on $V_{\bar{\sigma}}$ by $(g^*)^{-1}$, under the identification of $V_{\bar{\sigma}}$ with $V_\sigma^*$ via the bilinear pairing $(\bullet,\bullet)_V$.

\begin{equation}
\label{eq-diagram-of-Spin-groups-over-K-tilde}
\xymatrix{
\prod_{\hat{\sigma}\in\hat{\Sigma}}C(V_\sigma\oplus V_{\bar{\sigma}})^{even,\times} \ar[r] &
\left[\otimes_{\hat{\sigma}\in\hat{\Sigma}}(V_\sigma\oplus V_{\bar{\sigma}})\right]^{even,\times}
\ar[r]^-\cong &
C(V_{\tilde{K}})^{even,\times}
\\
\prod_{\hat{\sigma}\in\hat{\Sigma}}\Spin(V_\sigma\oplus V_{\bar{\sigma}})_\eta 
\ar[rr] \ar[u]^-\cup
& &
\Spin(V_{\tilde{K}})_\eta
\ar[u]^{\cup}.
}
\end{equation}

%************
% Hide
%************
\hide{
We have the left exact sequence
\[
1\rightarrow \Spin(V_{\tilde{K}})_\eta \rightarrow \prod_{\hat{\sigma}\in\hat{\Sigma}}G(V_\sigma\oplus V_{\bar{\sigma}})_\eta\RightArrowOf{j} \tilde{K}^\times\times \ZZ/2\ZZ,
\]
where 
\[
j((g_{\hat{\sigma}})_{\hat{\sigma}\in\hat{\Sigma}})=\left(\prod_{\hat{\sigma}\in\hat{\Sigma}}N(g_{\hat{\sigma}}),\sum_{\hat{\sigma}\in\hat{\Sigma}}p(g_{\hat{\sigma}})\right),
\]
$N:G(V_\sigma\oplus V_{\bar{\sigma}})\rightarrow\tilde{K}^\times$ is the norm\footnote{
The norm character $N:G(V_{\tilde{K}})\rightarrow \tilde{K}$ is defined by 
$N(g)=g\tau(g)$ as an element of the center $\tilde{K}$ of the Clifford algebra $C(V_{\tilde{K}})$. 
If $g=v_1\cdots v_k$, for $v_i\in V_{\tilde{K}}$, then $N(g)=\prod_{i=1}^k\frac{(v_i,v_i)}{2}$. See \cite[Sec. II.2.3]{chevalley}.
}
character,
and $p:G(V_\sigma\oplus V_{\bar{\sigma}})\rightarrow \ZZ/2\ZZ$ is the parity character. Denote by 
\[
N_{\hat{\sigma}}:\Spin(V_{\tilde{K}})_\eta\rightarrow \tilde{K}^\times \ \mbox{and} \ 
p_{\hat{\sigma}}:\Spin(V_{\tilde{K}})_\eta\rightarrow \ZZ/2\ZZ
\] 
the restriction of the norm and parity characters of the factor $G(V_\sigma\oplus V_{\bar{\sigma}})_\eta$. 
Let the character $\det_\sigma$ be the composition  
\[
\Spin(V_{\tilde{K}})_\eta \rightarrow G(V_\sigma\oplus V_{\bar{\sigma}})_\eta\rightarrow O(V_\sigma\oplus V_{\bar{\sigma}})_\eta\rightarrow GL(V_\sigma)\RightArrowOf{\det} \tilde{K}^\times.
\]
%************
% End Hide
%************
}
The bottom horizontal homomorphism in (\ref{eq-diagram-of-Spin-groups-over-K-tilde}) induces the isomorphism
\[
%\begin{equation}
%\label{eq-factorization-of-Spin-V-tilde-K-eta-mod-pm1}
\prod_{\hat{\sigma}\in\hat{\Sigma}}\left[\Spin(V_\sigma\oplus V_{\bar{\sigma}})_\eta /\{\pm 1\}\right]
\rightarrow
\Spin(V_{\tilde{K}})_\eta/\{\pm 1\},
\]
%\end{equation}
by an argument similar to the proof of Lemma \ref{lemma-Spin-V-hat-eta}.
We get the isomorphism
\begin{equation}
\label{eq-factorization-of-Spin-V-tilde-K-eta-mod-pm1}
\Spin(V_{\tilde{K}})_\eta/\{\pm 1\}\cong \prod_{\hat{\sigma}\in\hat{\Sigma}}GL(V_\sigma),
\end{equation}
where each $\sigma$ is chosen to restrict to $F$ as $\hat{\sigma}$.

Given $\sigma_0\in\Sigma$, 
let the character $\det_{\sigma_0}$ be the composition  
\begin{equation}
\label{eq-det-sigma}
\Spin(V_{\tilde{K}})_\eta \rightarrow \Spin(V_{\tilde{K}})_\eta/\{\pm 1\} 
%\RightArrowOf{(\ref{eq-factorization-of-Spin-V-tilde-K-eta-mod-pm1})^{-1}} 
%\prod_{\hat{\sigma}\in\hat{\Sigma}}\left[\Spin(V_\sigma\oplus V_{\bar{\sigma}})_\eta /\{\pm 1\}\right]
\IsomRightArrow
\prod_{\hat{\sigma}\in\hat{\Sigma}}GL(V_\sigma)
\rightarrow GL(V_{\sigma_0}) \RightArrowOf{\det} \tilde{K}^\times,
\end{equation}
where 
%the two direct products above involve a choice of $\sigma\in\Sigma$ restricting to $\hat{\sigma}$ and 
the third arrow is the projection onto the factor associated to $\sigma_0$.
The character $\ell_T$ of $\Spin(V_{\tilde{K}})_\eta$ satisfies 
\[
\ell_T\otimes\ell_T\cong \bigotimes_{\hat{\sigma}\in\hat{\Sigma}}\mbox{det}_{T(\hat{\sigma})}
\]
More generally, we have

\begin{lem}
\label{lemma-tensor-product-of-two-ells}
${\displaystyle \ell_{T_1}\otimes\ell_{T_2}\cong \bigotimes_{\{\hat{\sigma} \ : \ T_1(\hat{\sigma})=T_2(\hat{\sigma})\}}\det\!{}_{T_1(\hat{\sigma})}.
}$
\end{lem}

\begin{proof}
The isomorphism $\ell_\sigma\otimes\ell_\sigma\cong \det_\sigma$ of the $\Spin(V_\sigma\oplus V_{\bar{\sigma}})_\eta/{\pm 1}\cong GL(V_{\sigma})$ characters follows from \cite[III.3.2]{chevalley}. Hence, $\ell_\sigma\otimes\ell_{\bar{\sigma}}$ is the trivial characxter, as its square is $\det_{\bar{\sigma}}\otimes\det_\sigma$, which is the trivial character. 
The statement now follows from the factorization (\ref{eq-tensor-factorization-of-ell-T}).
\end{proof}

\begin{cor}
\label{cor-ell-T-are-linearly-independent}
The lines $\ell_T$, $T\in \T_K$, in $S^+_{\tilde{K}}$ are linearly independent over $\tilde{K}$.
\end{cor}

\begin{proof}
The group $\Spin(V_{\tilde{K}})_\eta$ leaves $B\otimes_\QQ\tilde{K}$ invariant and acts on the lines $\{\ell_T \ : \ T\in\T_K\}$ as characters. These characters are distinct, by 
Lemma \ref{lemma-tensor-product-of-two-ells}. 
\end{proof}

%\begin{question}
%Relate $N_{\hat{\sigma}}$ and $p_{\hat{\sigma}}$ to the pullback of the determinant character $\det_\sigma$.
%\end{question}

The image $m_g$, of $g\in \Spin(V)_{\eta,B}$ via the even half spin representation, leaves invariant each pure spinor $\ell_{T}$.
Given an element $g$ of $\Spin(V)_{\eta,B}$, its image $\rho_g$ via the vector representation
%even half spin representation 
leaves each of $W_T$, $T\in\T_K$,  invariant and acts on each with trivial determinant. Furthermore, $\rho_g$ commutes with $\hat{\eta}(F)$ and thus leaves invariant each of the subspaces $V_{\hat{\sigma},\RR}$ of $V_\RR$, for all $\hat{\sigma}\in\hat{\Sigma}$. Hence, $\rho_g$ leaves invariant each of the subspaces $V_{\sigma,\CC}$, for all $\sigma\in \Sigma$. 
%Consequently, $\rho_g$ leaves invariant each $W_{T'}$, $T'\in\T_K$. 
Denote by $\rho_{g,\sigma}$ the restriction of $\rho_g$ to $V_{\sigma,\CC}$.  Note that $\det_\sigma(g)=\det(\rho_{g,\sigma})$.

\begin{lem}
\label{lemma-det-rho-g-sigma=1}
For $g\in\Spin(V_\QQ)_{\eta,B}$ we have $\det(\rho_{g,\sigma})=1,$ for all $\sigma\in\Sigma$.
\end{lem}

\begin{proof}
Fix $\sigma_0\in\Sigma$ restricting to $F$ as $\hat{\sigma}_0\in\hat{\Sigma}$. 
Choose CM-types $T_1$ and $T_2$, such that $T_1(\hat{\sigma}_0)=T_2(\hat{\sigma}_0)$ and $T_1(\hat{\sigma})\neq T_2(\hat{\sigma})$, for $\hat{\sigma}\neq\hat{\sigma}_0$. 
Then $\ell_{T_1}\otimes\ell_{T_2}$ is the character $\det_{\sigma_0}$, by Lemma \ref{lemma-tensor-product-of-two-ells}. Thus, $\det_{\sigma_0}(g)=1$, for $g\in \Spin(V_\QQ)_{\eta,B}$, 
since $m_g$ fixes every point in $B$ and so $m_g\otimes m_g$ fixes every point in $B\otimes B$.
\end{proof}

%The following is an analogue of Lemma \ref{lemma-centralizer-of-rho-Spin-V-P}. 

%\begin{lem}
%The centralizer of $\rho(\Spin(V_\QQ)_{\eta,B})$ in $\tilde{O}(V_\QQ)$ is $\eta(K^\times)$.
%\end{lem}

\begin{lem}
\label{lemma-Spin-V-eta-tilde-K-is-product-of-SL-V-sigma}
The composition $\Spin(V_{\tilde{K}})_{\eta,B}\rightarrow \Spin(V_{\tilde{K}})_\eta/\{\pm 1\} \RightArrowOf{(\ref{eq-factorization-of-Spin-V-tilde-K-eta-mod-pm1})} \prod_{\hat{\sigma}\in\hat{\Sigma}}GL(V_{\sigma})$ maps $\Spin(V_{\tilde{K}})_{\eta,B}$ isomorphically onto 
$\prod_{\hat{\sigma}\in\hat{\Sigma}}SL(V_{\sigma})$, where each $\sigma$ is chosen to restrict to $F$ as $\hat{\sigma}$.
\end{lem}

\begin{proof}
The composition is injective, since (\ref{eq-factorization-of-Spin-V-tilde-K-eta-mod-pm1}) is an isomorphism and $-1$ does not belong to $\Spin(V_{\tilde{K}})_{\eta,B}$. The image of the composition is contained in $\prod_{\hat{\sigma}\in\hat{\Sigma}}SL(V_{\sigma}),$
by Lemma \ref{lemma-det-rho-g-sigma=1}. The image of the composition is the whole of $\prod_{\hat{\sigma}\in\hat{\Sigma}}SL(V_{\sigma})$,
since the subgroup of $\Spin(V_{\tilde{K}})_{\eta,B}$, generated by the images of $\Spin(V_{\sigma}\oplus V_{\bar{\sigma}})_{\eta,P_{\hat{\sigma}}}$, surjects onto 
$\prod_{\hat{\sigma}\in\hat{\Sigma}}SL(V_{\sigma})$.
\end{proof}

\begin{cor}
The subgroup $\Spin(V_{\tilde{K}})_B$ fixing all points of $B$ is equal to the subgroup $\Spin(V_{\tilde{K}})_{\eta,B}$.
\end{cor}

\begin{proof}
The homomorphism $\prod_{\hat{\sigma}\in\hat{\Sigma}}\Spin(V_\sigma\oplus V_{\bar{\sigma}})_{\eta,P_{\hat{\sigma}}}\rightarrow \Spin(V_{\tilde{K}})_{\eta,B}$, given in (\ref{eq-factorization-of-Spin-V-tilde-K-B}), is an isomorphism, by 
Lemma \ref{lemma-Spin-V-eta-tilde-K-is-product-of-SL-V-sigma}. The subgroup $\Spin(V_\sigma\oplus V_{\bar{\sigma}})_{\eta,P_{\hat{\sigma}}}$ is equal to 
the subgroup $\Spin(V_\sigma\oplus V_{\bar{\sigma}})_{P_{\hat{\sigma}}}$ fixing every point of $P_{\hat{\sigma}}$, by 
\cite[Lemma 2.2.2]{markman-sixfolds}.
%Lemma \ref{lemma-Spin-V-K-is-SL-n-K}.
Finally, $\Spin(V_{\tilde{K}})_B$ is generated by the images of $\Spin(V_\sigma\oplus V_{\bar{\sigma}})_{P_{\hat{\sigma}}}$, $\hat{\sigma}\in\hat{\Sigma}$.
\end{proof}

%*******************************************************************************************************************
%
%*******************************************************************************************************************
\subsection{$\Spin(V)_{\eta,B}$-invariant hermitian forms}
\label{sec-Hermitian-form-H-t}
Let $t$ be a non-zero element of the $(-1)$-eigenspace $K_-$ of $\iota:K\rightarrow K$. 
%************
% Hide
%************
\hide{
Choose a CM-type $T$. 
Consider the $F\otimes_\QQ\tilde{K}$-valued form $H_{t,T}:V_\QQ\times V_\QQ\rightarrow F\otimes_\QQ\tilde{K}\cong \prod_{\hat{\sigma}\in\hat{\Sigma}}\tilde{K}$ whose $\hat{\sigma}$ component is given by
\[
H_{t,T}(x_{\hat{\sigma}},y_{\hat{\sigma}}):= \hat{\sigma}(t\iota(t))(x_{\hat{\sigma}},y_{\hat{\sigma}})_V+\sigma(t)((\eta(t)x)_{\hat{\sigma}},y_{\hat{\sigma}})_V,
\]
where $\sigma=T(\hat{\sigma}):K\rightarrow \tilde{K}$, 
%restricts to $F$ as $\hat{\sigma}$, 
$x=\sum_{\hat{\sigma}\in\hat{\Sigma}}x_{\hat{\sigma}}$ with $x_{\hat{\sigma}}\in V_{\hat{\sigma}}$, and $y_{\hat{\sigma}}$ is defined similarly.
We have $H_{t,T}(x,y)=\iota(H_{t,T}(y,x))$, by Corollary \ref{cor-Neron-Severi-group}. 
%The choices of $\sigma$ above amount to a choice of a CM-type on which $H_t$ depends.

\begin{lem}
The form $H_{t,T}$ is hermitian,  i.e., 
\begin{equation}
\label{eq-K-linearity-in-second-variable}
H_{t,T}(x_{\hat{\sigma}},\eta(\lambda)y_{\hat{\sigma}})=\sigma(\lambda)H_{t,T}(x_{\hat{\sigma}},y_{\hat{\sigma}}), 
\end{equation}
for $\lambda\in K$, where $\sigma=T(\hat{\sigma})$. Furthermore, $H_{t,T}$ is $\Spin(V_\QQ)_{\eta,B}$-invariant.
\end{lem}

\begin{proof}
The form $H_{t,T}$ is $\Spin(V_\QQ)_{\eta,B}$-invariant, since $(\bullet,\bullet)_V$ is and $\eta(t)\rho_g(x)=\rho_g(\eta(t)x)$, by definition of $\Spin(V_\QQ)_{\eta,B}$.
We prove the $K$-linearity (\ref{eq-K-linearity-in-second-variable}) next.
Write $\lambda=a+bt$ with $a, b\in F$. The form $H_{t,T}$ on $V_{\hat{\sigma}}$ is $F$-bilinear, i.e., Equation (\ref{eq-K-linearity-in-second-variable}) holds for $\lambda\in F$. Hence, it suffices to prove the case $\lambda=t$. We have
\begin{eqnarray*}
H_{t,T}(x_{\hat{\sigma}},\eta(t)y_{\hat{\sigma}})&=&-\hat{\sigma}(t^2)(x_{\hat{\sigma}},\eta(t)y_{\hat{\sigma}})+
\sigma(t)(\eta(t)x_{\hat{\sigma}},\eta(t)y_{\hat{\sigma}})
\\
&\stackrel{\mbox{Cor. \ref{cor-Neron-Severi-group}}}{=}&
-\sigma(t)^2(x_{\hat{\sigma}},\eta(t)y_{\hat{\sigma}})+\sigma(t)(-\hat{\sigma}(t^2)x_{\hat{\sigma}},y_{\hat{\sigma}})
\\
&=&\sigma(t)\left[
\hat{\sigma}(-t^2)(x_{\hat{\sigma}},y_{\hat{\sigma}})-\sigma(t)(x_{\hat{\sigma}},\eta(t)y_{\hat{\sigma}})
\right]
\\
&\stackrel{\mbox{Cor. \ref{cor-Neron-Severi-group}}}{=}& \sigma(t)\left[
\hat{\sigma}(-t^2)(x_{\hat{\sigma}},y_{\hat{\sigma}})+\sigma(t)(\eta(t)x_{\hat{\sigma}},y_{\hat{\sigma}})
\right]
=\sigma(t)H_{t,T}(x_{\hat{\sigma}},y_{\hat{\sigma}}).
\end{eqnarray*}
%The image $\rho_g$, of  $g\in \Spin(V)_{\eta,B}$, restricts to $V_\sigma\oplus V_{\bar{\sigma}}$ as the image of an $\iota$-invariant element of 
%$\Spin(V_\sigma\oplus V_{\bar{\sigma}})_{\eta,P_{\hat{\sigma}}}.$ The latter preserves the restriction of $H_t$ to $V_\sigma\oplus V_{\bar{\sigma}}$,
%by the proof of Lemma \ref{lemma-su-3-3}.
\end{proof}

We saw in Section \ref{sec-Galois-action-on-set-of-CM-types} that the CM-type $T$ determines and embedding $T:K\rightarrow F\otimes_\QQ\tilde{K}$.
Equation (\ref{eq-K-linearity-in-second-variable}) becomes
\[
H_{t,T}(x,\eta(\lambda)y)_V=T(\lambda)H_{t,T}(x,y)_V.
\]
%************
% End Hide
%************
}

Let $H_t:V_\QQ\times V_\QQ\rightarrow K$ be given by
\begin{equation}
\label{eq-H-t}
H_t(x,y) :=
(-t^2)(x,y)_{V_{\hat{\eta}}}+t(\eta(t)x,y)_{V_{\hat{\eta}}}.
\end{equation}

\begin{lem}
\label{lemma-H-t-is-hermitian}
The form $H_t$ is hermitian, i.e., it satisfies $H_t(x,y)=\iota(H_t(y,x))$ and $H_t(x,\eta(\lambda)y)=\lambda H_t(x,y)$. Furthermore, $H_t$ is $\Spin(V_\QQ)_\eta$-invariant.
%There exists a unique $K$-valued hermitian form $H_t:V_\QQ\times V_\QQ\rightarrow K$, such that $H_{t,T}=T\circ H_t$.
%Furthermore, $H_t(x,y)=\iota(H_t(y,x))$ and $H_t(x,\eta(\lambda)y)=\lambda H_t(x,y)$.
\end{lem}

\begin{proof}
The equality $H_t(x,y)=\iota(H_t(y,x))$ follows from Lemma \ref{lemma-units-are-isometries}. The form $H_t$ is $F$-bilinear. Hence, the equality 
$H_t(x,\eta(\lambda)y)=\lambda H_t(x,y)$, for all $\lambda\in K$, would follow from the case $\lambda=t$. We have
\begin{eqnarray*}
H_t(x,\eta(t)y)&=& -t^2(x,\eta(t)y)_{V_{\hat{\eta}}}+t(\eta(t)x,\eta(t)y)_{V_{\hat{\eta}}}
\\
&\stackrel{\mbox{Lemma \ref{lemma-units-are-isometries}}}{=}&
-t^2(x,\eta(t)y)_{V_{\hat{\eta}}}+t(\hat{\eta}(-t^2)x,y)_{V_{\hat{\eta}}}
\\
&=&-t^2(x,\eta(t)y)_{V_{\hat{\eta}}}+-t^3(x,y)_{V_{\hat{\eta}}}
\\
&\stackrel{\mbox{Lemma \ref{lemma-units-are-isometries}}}{=}&
t\left[
-t^2(x,y)_{V_{\hat{\eta}}}+t(\eta(t)x,y)_{V_{\hat{\eta}}}
\right]=tH_t(x,y).
\end{eqnarray*}

The pairing $(\bullet,\bullet)_{V_{\hat{\eta}}}$ is $\Spin(V_{\hat{\eta}})$-invariant. It follows that it is also $\Spin(V_\QQ)_{\hat{\eta}}$-invariant, by Lemma \ref{lemma-Spin-V-hat-eta}, which implies that the two groups have the same image in $GL(V_\QQ)$. Now $\Spin(V_\QQ)_\eta$ is a subgroup of $\Spin(V_\QQ)_{\hat{\eta}}$, by definition.
The form $H_t$ is $\Spin(V_\QQ)_\eta$-invariant, since $(\bullet,\bullet)_{V_{\hat{\eta}}}$ is and $\eta(t)\rho_g(x)=\rho_g(\eta(t)x)$, by definition of $\Spin(V_\QQ)_\eta$.
\end{proof}

%************
% Hide
%************
\hide{
\begin{proof}
It suffices to prove the equality $H_{t,T}=T\circ H_t$.
The embedding $T:K\rightarrow F\otimes_\QQ\tilde{K}$ restricts to $F$ as the standard embedding of $F$ in $F\otimes_\QQ\tilde{F}$. 
%In particular, given $x, y\in V_\QQ$, 
%$(x,y)_{V_{\hat{\eta}}}$ belongs to $T(F)$ and so does $(\eta(t)x,y)_{V_{\hat{\eta}}}$. 
Now,
\begin{eqnarray*}
H_{t,T}(x,y)&=&T(-t^2)\left[T\left((x,y)_{V_{\hat{\eta}}}\right)+T(t)T\left((\eta(t)x,y)_{V_{\hat{\eta}}}\right)
\right]
\\
&=&T\left[
(-t^2)(x,y)_{V_{\hat{\eta}}}+t(\eta(t)x,y)_{V_{\hat{\eta}}}
\right]=T(H_t(x,y)).
\end{eqnarray*}
%and so it has values in $T(K)$.
\end{proof}
%************
% End Hide
%************
}

Let $H:V_\eta\times V_\eta\rightarrow K$ be a non-degenerate hermitian form. 
Choose an orthogonal basis $\B=\{v_1, \dots, v_d\}$ of $V_\eta$ as a vector space over $K$ and let $[H]_\B$ be the hermitian matrix of $H$ in the basis $\B$. 
Given an embedding $\hat{\sigma}:F\rightarrow \RR$, let $a_{\hat{\sigma}}$ be the number of basis elements $v_i$, such that 
$\hat{\sigma}(H(v_i,v_i))>0$ and let $b_{\hat{\sigma}}$ be the number of basis elements $v_i$, such that 
$\hat{\sigma}(H(v_i,v_i))<0$. The non-degeneracy of $H$ implies that $a_{\hat{\sigma}}+b_{\hat{\sigma}}=d$.

\begin{defi}
\label{def-discriminant}(\cite[Sec. 4]{deligne-milne})
  The {\em discriminant} $\disc(H)$ is the coset of $\det([H]_\B)$ in $F^\times/Nm_{K/F}(K^\times)$ and $\disc(H)$ is independent of the choice of $\B$.
\end{defi}

%Given an embedding $\hat{\sigma}:F\rightarrow \RR$, the hermitian form $H_t$ induces a $K\otimes_F\RR$-valued hermitian form 
%on $V_{\hat{\sigma},\RR}$ and we let $a_{\hat{\sigma}}$ be the maximal dimension of a subspace on which it is positive definite and $b_{\hat{\sigma}}$ 
%the maximal dimension of a subspace on which it is negative definite. 

\begin{lem} 
\label{lemma-invariants-of-hermitian-forms}
\begin{enumerate}
\item
(\cite{lan})
The signatures $(a_{\hat{\sigma}},d-a_{\hat{\sigma}})$, $\hat{\sigma}\in\hat{\Sigma}$, and the discriminant $\disc(H_t)$ determine the isomorphism class of 
$(V_\eta, H)$. Given integers $0\leq a_{\hat{\sigma}}\leq d$ and an element $\delta\in F^\times/Nm_{K/F}(K^\times)$, there exists a $K$-valued hermitian form $H$ on $V_\eta$ with these invariants.
\item
\label{lemma-item-two-characterizations-of-split-type}
(\cite[Cor. 4.2]{deligne-milne})
The following are equivalent:
\begin{itemize}
\item
$a_{\hat{\sigma}}=b_{\hat{\sigma}}$, for all $\hat{\sigma}\in\hat{\Sigma}$, and $\disc(H)=(-1)^{d/2}$;
\item
there exists an isotropic subspace of $(V_\eta,H)$ of dimension $d/2$.
\end{itemize}
\end{enumerate}
\end{lem}

\begin{defi}
\label{def-split-type}
The pair $(V_\eta,H)$ is said to be of {\em split type}, if it satisfies the conditions of
Lemma \ref{lemma-invariants-of-hermitian-forms} (\ref{lemma-item-two-characterizations-of-split-type}).
\end{defi}

%*******************************************************************************************************************
%
%*******************************************************************************************************************
\subsection{If $B$ is spanned by Hodge classes then $HW$ is as well}
\label{sec-if-B-is-Hodge-so-in-HW}
%The special Mumford-Tate group of the generic deformation of $(X\times\hat{X},\eta)$  contains the Zariski closure in $SO(V_\QQ)$
%of the image via $\rho$ of subgroup 
%$\Spin(V_\QQ)_{\eta,B}$ 
%of $\Spin(V_\QQ)$. This follows from the following Lemma.
In this section we show that if the secant space $B$ is spanned by Hodge classes, then condition 
(\ref{eq-condition-for-HW-to-consist-of-Hodge-classes}) is satisfied for the $K$-action $\eta$ on 
$H^1(X\times\hat{X},\QQ)$, and so the subspace $HW$ of $H^d(X\times\hat{X},\QQ)$ consists of Hodge classes.

\begin{lem}
\label{lemma-if_B-is-spanned-by-Hodge-classes-so-does-HW}
Assume that $B$ is spanned by Hodge classes. Then the complex structure $I_{X\times\hat{X}}:H^1(X\times\hat{X},\RR)\rightarrow H^1(X\times\hat{X},\RR)$
of $X\times\hat{X}$ commutes with the action of $\eta(K)$. 
Furthermore, 
%the subspace $B$ is spanned by Hodge classes, then 
\begin{equation}
\label{eq-dim-V-sigma-0-1=dim-V-sigma-1-0}
\dim(V_{\sigma,\CC}^{0,1})=\dim(V_{\sigma,\CC}^{1,0}),
\end{equation}
for all $\sigma\in\Sigma$, and consequently $HW$ is spanned by Hodge classes.
\end{lem}

\begin{proof}
Set $I:=I_{X\times\hat{X}}$. Let $T\in\T_K$ be a CM-type. We get the two complex conjugate lines $\ell_T$ and $\ell_{\iota(T)}$ in $B$, each spanned by a Hodge class. 
The proof of \cite[Lemma 2.2.6]{markman-sixfolds}
%Lemma \ref{lemma-decomposition-into-4-direct-summands} 
establishes that each of $W_{T,\CC}$ and $\bar{W}_{T,\CC}=W_{\iota(T),\CC}$ is $I$-invariant and the dimensions of the $I$-eigenspaces in each are equal,
\begin{eqnarray}
\label{eq-W-1-0=W-0-1}
\dim(W^{1,0}_{T,\CC})&=&\dim(W^{0,1}_{T,\CC})=2n
\\
\dim(\bar{W}^{1,0}_{T,\CC})&=&\dim(\bar{W}^{0,1}_{T,\CC})=2n.
\nonumber
\end{eqnarray} 
The endomorphisms in $\eta(F)$ commute with $I$, by assumption. Hence, the subspace $V_{\hat{\sigma},\RR}:=H^1(X\times\hat{X},\RR)_{\hat{\sigma}}$ is $I$-invariant, for all $\hat{\sigma}$ in $\hat{\Sigma}$. We get the equality
\begin{equation}
\label{eq-dim-V-hat-sigma-1-0-is-d}
\dim(V_{\hat{\sigma},\CC}^{1,0})=\dim(V_{\hat{\sigma},\CC}^{0,1})=d,
\end{equation}
for all $\hat{\sigma}$ in $\hat{\Sigma}$.

The following equalities hold for all $\sigma\in T$
\begin{eqnarray*}
V_{\sigma,\CC} & = & [V_{\hat{\sigma},\RR}\otimes_\RR\CC]\cap W_{T,\CC},
\\
V_{\bar{\sigma},\CC} & = & [V_{\hat{\sigma},\RR}\otimes_\RR\CC]\cap \bar{W}_{T,\CC},
\end{eqnarray*}
by construction. Hence, the $\eta(K)$ eigenspaces $V_{\sigma,\CC}$ are $I$-invariant, for all $\sigma\in\Sigma$.

We know the equality
\begin{equation}
\label{eq-sum-of-dimensions}
\sum_{\sigma\in T}\dim(V_{\sigma,\CC}^{0,1})=\sum_{\sigma\in T}\dim(V_{\sigma,\CC}^{1,0})=2n,
\end{equation}
by equality (\ref{eq-W-1-0=W-0-1}).
%If $B$ is spanned by Hodge classes, then the subspace $P_{T'}$ of $H^{ev}(X,\RR)$ is contained in $H^{1,1}(X,\RR)$. 
%The proof of Lemma \ref{lemma-decomposition-into-4-direct-summands} implies equation (\ref{eq-W-1-0=W-0-1}) and hence 
%equation (\ref{eq-sum-of-dimensions}) for all $T'\in\T_K$. 
We prove Equation (\ref{eq-dim-V-sigma-0-1=dim-V-sigma-1-0}) by contradiction. 
Assume that $\dim(V^{1,0}_{\sigma_0,\CC})>\dim(V^{0,1}_{\sigma_0,\CC})$, for some $\sigma_0\in\Sigma$.  Then $\dim(V^{1,0}_{\bar{\sigma}_0,\CC})<\dim(V^{0,1}_{\bar{\sigma}_0,\CC})$, 
and so $\dim(V^{1,0}_{\sigma_0,\CC})>d/2>\dim(V^{1,0}_{\bar{\sigma}_0,\CC})$,
by (\ref{eq-dim-V-hat-sigma-1-0-is-d}). 
Thus, if $\sigma_0$ belongs to $T_1$ and we replace it by $\bar{\sigma}_0$ to obtain a CM-type $T_2$, then
%changing the value of a CM-type $T$ at $\sigma_0$ changes 
the second sum in (\ref{eq-sum-of-dimensions}) for $T_1$ will be different from the sum for $T_2$. This contradicts the fact that the sum is $2n$, for all CM-types $T$.
Equation (\ref{eq-dim-V-sigma-0-1=dim-V-sigma-1-0}) thus hold for all $\sigma\in\Sigma$.
The equalities (\ref{eq-dim-V-sigma-0-1=dim-V-sigma-1-0}) hold for all $\sigma\in\Sigma$, if and only if $HW$ consists of rational Hodge classes, by \cite[Prop. 4.4]{deligne-milne}.
\end{proof}

%*******************************************************************************************************************
%
%*******************************************************************************************************************
\section{The $\Spin(V)_{\eta,B}$-invariant subalgebra of $H^*(X\times\hat{X},\QQ)$}
\label{sec-Mumford-Tate-groups}
Let $W$ be the maximal isotropic subspace (\ref{eq-W}). The associated complex multiplication 
$\eta:K\rightarrow \End_\QQ[H^1(X,\QQ)\oplus H^1(\hat{X},\QQ)]=\End_\QQ(V_\QQ)$ is given in (\ref{eq-eta}).
Let $B\subset H^{ev}(X,\QQ)$ be the associated secant linear subspace in Lemma \ref{lemma-B-is-rational}.
Let $\A:=(\wedge^*V_\QQ)^{\Spin(V)_{\eta,B}}$ be the subalgebra of $\wedge^*V_\QQ$
of $\Spin(V)_{\eta,B}$-invariant classes with respect to  the natural extension of the representation $\rho:\Spin(V_\QQ)\rightarrow SO_+(V_\QQ)$ to the exterior algebra. 
%given in (\ref{eq-rho-extended-to-exterior-algebra}). 

By $\wedge^*V_{\tilde{F}}$ we will refer to $\wedge^*_{\tilde{F}}V_{\tilde{F}}$, where $V_{\tilde{F}}=V_\QQ\otimes_\QQ\tilde{F}$.

\begin{prop}
\label{prop-generator-for-the-invariant-subalgebra}
\begin{enumerate}
\item
\label{lemma-item-Spin-V-eta-B-invariant-2-forms}
The subspace $(\wedge^2V_\QQ)^{\Spin(V)_{\eta,B}}$ is equal to the image of the injective homomorphism $\Xi:K_-\rightarrow  \wedge^2V_\QQ$
in Equation (\ref{eq-Xi}). 
%of Corollary (\ref{cor-Neron-Severi-group}). 
This subspace is thus $(e/2)$-dimensional over $\QQ$ and is naturally a $1$-dimensional  vector space over $F$, and it
consists of $\Spin(V)_\eta$-invariant classes. 
\item
\label{lemma-item-generator-for-the-invariant-subalgebra}
The subalgebra $\A$ is commutative and is generated by the $(e/2)$-dimensional $\QQ$-vector space 
$(\wedge^2V_\QQ)^{\Spin(V)_{\eta,B}}$
and the $e$-dimensional $\QQ$-vector space 
$
HW=\wedge^d_{\eta(K)}V_\QQ.
$
\item
\label{lemma-item-decomposition-of-HW-into-characters}
The subspace  $HW\otimes_\QQ\tilde{K}$ decomposes as a direct sum of the characters $\det_\sigma$, $\sigma\in\Sigma$, of $\Spin(V_{\tilde{K}})_\eta$, each with multiplicity $1$.
\item
\label{lemma-item-generatyors-for-the-invariant-subalgebra-R}
The subalgebra $\R:=(\wedge^*V_\QQ)^{\Spin(V)_\eta}$ is generated by $(\wedge^2V_\QQ)^{\Spin(V)_{\eta,B}}$.
\end{enumerate}
\end{prop}

%Note that the algebra $(\wedge^*V_\QQ)^{\Spin(V)_{\eta,B}}$ need not be generated by $(\wedge^2V_\QQ)^{\Spin(V)_{\eta,B}}$ and $HW$, 
%since the $Gal(\tilde{F}/\QQ)$-invariant subalgebra of $(\wedge^*V_{\tilde{F}})^{\Spin(V)_{\eta,B}}$ may be larger than the subalgebra 
%generated by the $Gal(\tilde{F}/\QQ)$-invariant subspaces of the generators 
%in part (\ref{lemma-item-generator-for-the-invariant-subalgebra}) above.

\begin{proof}
(\ref{lemma-item-generator-for-the-invariant-subalgebra}) \underline{Step 1:}
We prove first that the subalgebra $(\wedge^*V_{\tilde{F}})^{\Spin(V)_{\eta,B}}$ is commutative and is generated by the $(e/2)$-dimensional $\tilde{F}$-vector space 
$(\wedge^2V_\QQ)^{\Spin(V)_{\eta,B}}\otimes_\QQ{\tilde{F}}$
and the $e$-dimensional $\tilde{F}$-vector space 
$
HW\otimes_\QQ{\tilde{F}}.
%=(\wedge^dV_{\tilde{F}})^{\Spin(V)_{\eta,B}}.
$

We have the isomorphism
\[
\wedge^*V_{\tilde{F}}\cong \bigotimes_{\hat{\sigma}\in\hat{\Sigma}}
\wedge^*V_{\hat{\sigma}}.
\]
The $\Spin(V)_{\eta,B}$-invariant subalgebra is equal to the $\Spin(V_{\tilde{F}})_{\eta,B}$-invariant subalgebra. The factorization 
(\ref{eq-factorization-of-Spin-V-tilde-K-B}) is defined over $\tilde{F}$ and yields an analogous factorization of  $\Spin(V_{\tilde{F}})_{\eta,B}$. Using the latter factorization we get the isomorphism
\begin{equation}
\label{eq-factorization-of-Spin-V-B-invariant-subalgebra-over-tilde-K}
(\wedge^*V_{\tilde{F}})^{\Spin(V)_{\eta,B}}\cong
\bigotimes_{\hat{\sigma}\in\hat{\Sigma}}
(\wedge^*V_{\hat{\sigma}})^{\Spin(V_{\hat{\sigma}})_{P_{\hat{\sigma}}}}.
\end{equation}
The argument of \cite[Lemma 2.2.7]{markman-sixfolds}
%Lemma \ref{lemma-Spin-V-P-invariant-classes-are-Hodge} 
applies, replacing $\QQ$ by $\tilde{F}$ and $K$ by $\tilde{K}$, and yields that 
\[
\dim_{\tilde{F}}\left((\wedge^kV_{\hat{\sigma}})^{\Spin(V_{\hat{\sigma}})_{P_{\hat{\sigma}}}}\right)=
\left\{
\begin{array}{ccl}
0&\mbox{if}&k \ \mbox{is odd}
\\
1&\mbox{if} & k \ \mbox{is even and} \ k\neq d
\\
3&\mbox{if}& k=d.
\end{array}
\right.
\]
Furthermore, $(\wedge^dV_{\hat{\sigma}})^{\Spin(V_{\hat{\sigma}})_{P_{\hat{\sigma}}}}$ decomposes as a direct sum of three one-dimensional   
$\Spin(V_{\hat{\sigma}}\otimes_{\tilde{F}}\tilde{K})_\eta$-invariant subspaces corresponding to the trivial character and the characters $\det_\sigma$ and $\det_{\bar{\sigma}}$, where $\sigma$ restricts to $\hat{\sigma}$. The direct sum of the latter two is of the form $HW_{\hat{\sigma}}\otimes_{\tilde{F}}\tilde{K}$
of a $2$-dimensional $\tilde{F}$ subspace $HW_{\hat{\sigma}}$. 
Then $HW\otimes_\QQ\tilde{K}=\oplus_{\hat{\sigma}\in\hat{\Sigma}}HW_{\hat{\sigma}}$.
Choose an element $\Xi_{\hat{\sigma}}$ spanning $(\wedge^2V_{\hat{\sigma}})^{\Spin(V_{\hat{\sigma}})_{P_{\hat{\sigma}}}}$. Then $\Xi_{\hat{\sigma}}^{d/2}$ spans that trivial character in $(\wedge^dV_{\hat{\sigma}})^{\Spin(V_{\hat{\sigma}})_{P_{\hat{\sigma}}}}$ and
the symmetric square of the subspace $HW_{\hat{\sigma}}$ is mapped onto $\wedge^{2d}V_{\hat{\sigma}}$, which is spanned by $\Xi_{\hat{\sigma}}^d$.
The tensor factors on the righthand side of (\ref{eq-factorization-of-Spin-V-B-invariant-subalgebra-over-tilde-K}) are all even, hence they commute with respect to wedge-product. We conclude that $\A_{\tilde{F}}:=(\wedge^*V_{\tilde{F}})^{\Spin(V)_{\eta,B}}$ is a commutative graded subalgebra of $\wedge^*V_{\tilde{F}}$ and it is generated by $HW\otimes_\QQ\tilde{F}$ and its graded summand $\A_{\tilde{F}}^2$ in degree $2$.

We conclude that $(\wedge^2V_{\tilde{F}})^{\Spin(V)_{\eta,B}}$ is an $(e/2)$-dimensional $\tilde{F}$-vector space with basis 
$\{\Xi_{\hat{\sigma}} \ : \ \hat{\sigma}\in\hat{\Sigma}\}$. Taking $Gal(\tilde{F}/\QQ)$ invariants we get that
$(\wedge^2V_\QQ)^{\Spin(V)_{\eta,B}}$ is an $(e/2)$-dimensional $\QQ$-vector space. Furthermore,
\[
(\wedge^2V_\QQ)^{\Spin(V)_{\eta,B}}=\{\Xi_t \ : \ t\in K_-\},
\]
 where $\Xi_t$ is given in (\ref{eq-Xi-t}).

\underline{Step 2:} Choose a basis $\{\beta_1, \dots, \beta_{e/2}\}$ of $\A^2_\QQ:=(\wedge^2V_\QQ)^{\Spin(V)_{\eta,B}}$ and a basis 
$\{\delta_1, \dots, \delta_e\}$ of $HW$. Then every element of $(\wedge^*V_{\tilde{F}})^{\Spin(V)_{\eta,B}}$ is the sum 
$\gamma:=\sum_i c_if_i$, where $c_i\in\tilde{F}$ and $f_i$ is a polynomial in the basis elements $\beta_j$'s and $\delta_k$'s. 
These two bases consist of rational classes. 
We may  further assume that the $f_i$'s are linearly independent over $\QQ$. Then they are also linearly independent over $\tilde{F}$,
since the homomorphism $\wedge^*_\QQ V_\QQ\rightarrow (\wedge^*_\QQ V_\QQ)\otimes_\QQ\tilde{F}$, given by $f\mapsto f\otimes 1$, sends linearly independent subsets over $\QQ$ to linearly independent subsets over $\tilde{F}$.
Hence, if $\gamma$ is $Gal(\tilde{F}/\QQ)$-invariant, then so is each of the coefficients $c_i$. 
Consequently, the algebra $(\wedge^*V_\QQ)^{\Spin(V)_{\eta,B}}$ is generated by $\A^2_\QQ$ and $HW$.

(\ref{lemma-item-Spin-V-eta-B-invariant-2-forms})
Let us verify the $\Spin(V)_\eta$-invariance of $(\wedge^2V_\QQ)^{\Spin(V)_{\eta,B}}$.  For $g\in\Spin(V)_\eta$ and $x,y\in V_\QQ$ we have
\[
\Xi_t(\rho_g(x),\rho_g(y)):=(\eta(t)\rho_g(x),\rho_g(y))_V=(\rho_g(\eta(t)x),\rho_g(y))_V=(\eta(t)x,y)_V=\Xi_t(x,y).
\]

(\ref{lemma-item-decomposition-of-HW-into-characters}) Follows from the direct sum decomposition (\ref{eq-direct-sum-decomposition-of-HW}).

(\ref{lemma-item-generatyors-for-the-invariant-subalgebra-R}) 
Assume that $x\in\wedge^*_\QQ V_\QQ$ is $\Spin(V)_\eta$-invariant. Then each of its factors $x_{\hat{\sigma}}$ in $(\wedge^*V_{\hat{\sigma}})^{\Spin(V)_{\hat{\sigma}}}$, with respect to the factorization (\ref{eq-factorization-of-Spin-V-B-invariant-subalgebra-over-tilde-K}), is $\Spin(V)_\eta$-invariant, as different factors contribute linearly independent characters over $\ZZ$. Hence, $x_{\hat{\sigma}}$ belongs to the subalgebra generated by the one-dimensional subspace $\Xi_{\hat{\sigma}}$ of $\wedge^2V_{\hat{\sigma}}$. Thus, $x$ is a rational class of the subalgebra generated by 
$(\wedge^2V_{\tilde{F}})^{\Spin(V)_{\eta,B}}$. Arguing as in Step 2 of the proof of part (\ref{lemma-item-generator-for-the-invariant-subalgebra}) we conclude that $x$ belongs to $\R$.
\end{proof}

Let $\Sigma'\subset\Sigma$ be a subset, such that $\Sigma'\cap\iota(\Sigma')=\emptyset$. 
Set $\det_{\Sigma'}:=\otimes_{\{\sigma\in\Sigma'\}}\det_\sigma.$
Let $k$ be the cardinality of $\Sigma'$.
Denote by $(\wedge^{dj}V_{\tilde{K}})_{\Sigma'}$ the isotypic subspace of $\wedge^{dj}V_{\tilde{K}}$ corresponding to the character $\det_{\Sigma'}$ of
$\Spin(V_{\tilde{K}})_\eta$. Note that $(\wedge^{dj}V_{\tilde{K}})_{\Sigma'}$ vanishes, if $j<k$ or $j>e-k$, and it is one-dimensional if $j=k$ or if $j=e-k$.
Let $(\wedge^{dj}V_{\tilde{K}})_k$ be the direct sum of all $(\wedge^{dj}V_{\tilde{K}})_{\Sigma'}$, for subsets $\Sigma'\subset\Sigma$ as above of cardinality $k$.
The dimension of $(\wedge^{dk}V_{\tilde{K}})_k$ is $\Choose{e/2}{k}2^k$.
Note that $(\wedge^{dk}V_{\tilde{K}})_k$ is contained in the subalgebra generated by $HW$ and 
$(\wedge^{d}V_{\tilde{K}})_d=HW\otimes_\QQ\tilde{K}.$

\begin{lem}
If $d>2$, then $HW$ is contained in the primitive cohomology with respect to every $\Spin(V)_{\eta,B}$-invariant ample class. 
%the cup product induces the zero homomorphism $(\wedge^2V_\QQ)^{\Spin(V)_{\eta,B}}\otimes HW\rightarrow \wedge^{d+2}V_\QQ)^{\Spin(V)_{\eta,B}}$.
More generally, the same is true for $(\wedge^{dk}V_{\tilde{K}})_k$.
\end{lem}

\begin{proof}
Let $h$ be a class in $(\wedge^2V_\QQ)^{\Spin(V)_{\eta,B}}$. Then $h$ is $\Spin(V_\QQ)_\eta$-invariant, 
by Proposition \ref{prop-generator-for-the-invariant-subalgebra}(\ref{lemma-item-Spin-V-eta-B-invariant-2-forms}).
%Note the equality $n-d=d(e-1)$ and so the classes in $(\wedge^{d(e-1)+2}V_\QQ)^{\Spin(V)_{\eta,B}}$ are all $\Spin(V_\QQ)_\eta$-invariant,
%by Lemma \ref{prop-generator-for-the-invariant-subalgebra}, as so are the classes in $(\wedge^{d-2}V_\QQ)^{\Spin(V)_{\eta,B}}$.
%Then $h^{d(e-2)/2+1}HW$ consists of $\Spin(V_\QQ)_\eta$-invariant classes. On the other hand, as $HW\otimes_\QQ\CC$ is a 
%direct sum of non-trivial $\Spin(V_\QQ)_\eta$-characters, and $h$ is $\Spin(V_\QQ)_\eta$-invariant, 
%then if non-zero, $h^{1+d(e-2)/2}HW$ is a direct sum of non-trivial characters.
%
%The case of $(\wedge^{dk}V_{\tilde{K}})_k$ is similar. 
Thus $h^{1+d(e-2k)/2}(\wedge^{dk}V_{\tilde{K}})_k$ is contained in 
$(\wedge^{d(e-k)+2}V_{\tilde{K}})_k$, which vanishes.
\end{proof}

%Let $\R$ be the subalgebra of $\wedge^*V_\QQ$ generated by $(\wedge^2V_\QQ)^{\Spin(V)_{\eta,B}}$.
%Let $\R'$ be the $\left((\wedge^2V_{\tilde{F}})^{\Spin(V)_{\eta,B}}\right)^{Gal(\tilde{F}/\QQ)}$. Note that $\R$ is a subalgebra of $\R'$ and the two 
%coincide if the $Gal(\tilde{F}/\QQ)$-action on $\hat{\Sigma}$ induces a surjection of $Gal(\tilde{F}/\QQ)$ onto the full permutation group of $\hat{\Sigma}$.

%\begin{lem}
%\label{lemma-generators-of-the-module-of-Spin-V-B-invariant-classes}
%The $\R$-subalgebra $(\wedge^*V_\QQ)^{\Spin(V)_{\eta,B}}$ of $\wedge^*V_\QQ$ is generated, as an $\R$-module,
%by the subspaces $(\wedge^{dk}V_{\tilde{K}})_k$, $0\leq k\leq e/2$.
%\end{lem}

%*******************************************************************************************************************
%
%*******************************************************************************************************************
\section{An adjoint orbit in $\Spin(V_\RR)_B$ as a period domain of abelian varieties of Weil type}
\label{sec-adjoint-orbit-n-Spin-V-RR-eta-B}

%*******************************************************************************************************************
%
%*******************************************************************************************************************
\subsection{The group $\Spin(V_\RR)_B$}
\label{sec-Spin-V-R-B}

Let $\Spin(V_\RR)_B$ be the subgroup of $\Spin(V_\RR)$ consisting of elements $g$, such that $m_g$ leaves every point of $B$ (and hence also of $B_\RR:=B\otimes_\QQ\RR$) invariant. 
Choose a non-zero element $t$ of $K_-$. Given $\sigma\in\Sigma$ restricting to $F$ as $\hat{\sigma}$, let 
\begin{equation}
\label{eq-SO+V-RR-hat-sigma-eta-t}
SO_+(V_{\hat{\sigma},\RR})_{\eta_t}
\end{equation} 
be the subgroup of 
$SO_+(V_{\hat{\sigma},\RR})$ of elements $g$, which commute with the restriction of $\eta_t$ to $V_{\hat{\sigma},\RR}$ and which restriction to each of $V_{\sigma,\CC}$ and $V_{\bar{\sigma},\CC}$ has determinant $1$.

%*********
% Hide
%*********
\hide{
and let $H_t$ be the corresponding $K$-valued hermitian form on $V_\QQ$ given in (\ref{eq-H-t}).
Denote the associated $K\otimes_\QQ\CC\cong\prod_{\sigma\in\Sigma} \CC$-valued form on $V_\RR$ and $V_\CC$ by $H_t$ as well.
Denote by $H_{t,\sigma}$ the $\CC$-valued form,
which is the restriction of $\sigma\circ H_t$ to $V_{\hat{\sigma},\RR}$, where $\sigma$ restricts to $F$ as $\hat{\sigma}$.
% picking up the factors corresponding to the two elements of $\Sigma$ restricting to $F$ as $\hat{\sigma}$. 
 Let $SU(V_{\hat{\sigma},\RR})$ be the subgroup of
$SO_+(V_{\hat{\sigma},\RR})$ consisting of elements leaving $H_{t,\sigma}$ invariant. Note that $SU(V_{\hat{\sigma},\RR})$ is independent of the choice of $\sigma$ and of $t\in K_-\setminus\{0\},$ as different choices lead to rescaling the imaginary part of $H_{t,\sigma}$ by a real factor.
%*********
% End Hide
%*********
}
\begin{lem}
\label{lemma-factorization-of-rho-of-Spin-V-B}
The representation $\rho:\Spin(V_\RR)\rightarrow SO_+(V_\RR)$ maps $\Spin(V_\RR)_B$ isomorphically onto 
$\prod_{\hat{\sigma}\in\hat{\Sigma}}SO_+(V_{\hat{\sigma},\RR})_{\eta_t}$.
%SU(V_{\hat{\sigma},\RR})$.
\end{lem}

\begin{proof}
The representation $\rho$ maps $\Spin(V_\RR)_B$ injectively into $SO_+(V_\RR)$, since $-1\in\Spin(V_\RR)$ does not belong to $\Spin(V_\RR)_B$. 
We prove first the inclusion of $\rho(\Spin(V_\RR)_B)\subset \prod_{\hat{\sigma}\in\hat{\Sigma}}SO_+(V_{\hat{\sigma},\RR})_{\eta_t}$.
We have the inclusion
\begin{equation}
\label{eq-rho-of-Spin-V-B}
\rho(\Spin(V_\RR)_B)\subset 
\{g\in SO_+(V_\RR) \ : \ g(V_\sigma)=V_\sigma \ \mbox{and} \ \det(\restricted{g}{V_\sigma})=1, \forall \ \sigma\in\Sigma \}.
\end{equation}
Indeed, fix $\sigma_0\in \Sigma$ restricting to $\hat{\sigma}_0\in\hat{\Sigma}$ and let $T_1$ and $T_2$ be two CM-types, such that $T_1(\hat{\sigma}_0)=T_2(\hat{\sigma}_0)=\sigma_0$ and $T_1(\hat{\sigma})\neq T_2(\hat{\sigma}),$ for all $\hat{\sigma}\neq\hat{\sigma}_0$. Then $g\in\rho(\Spin(V_\RR)_B)$
maps $\ell_{T_i}$ to itself, for $i=1,2$, hence it maps $W_{T_i}$ to itself, for $i=1,2$, and thus it maps $V_\sigma=W_{T_1}\cap W_{T_2}$ to itself.
The equality $\det(\restricted{g}{V_\sigma})=1$ is proved by the same argument proving Lemma \ref{lemma-det-rho-g-sigma=1}.
We conclude that $\rho(\Spin(V_\RR)_B)$ commutes with $\eta(K)$. The right hand side of (\ref{eq-rho-of-Spin-V-B}) is equal to 
$\prod_{\hat{\sigma}\in\hat{\Sigma}}SO_+(V_{\hat{\sigma},\RR})_{\eta_t}$.
%*********
% Hide
%*********
\hide{
The $\prod_{\sigma\in\Sigma} \CC$-valued hermitian form $H_t$
is invariant by elements in the right hand side of (\ref{eq-rho-of-Spin-V-B}), by its definition. Indeed, given $x,y\in V_\RR$, write $x=\sum_{\hat{\sigma}\in\hat{\Sigma}}x_{\hat{\sigma}}$ and define $y_{\hat{\sigma}}$ similarly. Then
\[
H_t(x,y)=\left(H_{t,\sigma}(x_{\hat{\sigma}},y_{\hat{\sigma}})
\right)_{\sigma\in\Sigma}
\]
and elements in the right hand side of (\ref{eq-rho-of-Spin-V-B}) leave $H_{t,\sigma}$ invariant, as they restrict to $V_{\hat{\sigma},\RR}$ as elements of 
$\rho(\Spin(V_{\hat{\sigma},\RR})_{P_{\hat{\sigma}}}$) (see Lemmas \ref{lemma-stabilizer-is-isomorphic-to-so-f} and \ref{lemma-su-3-3}). We conclude that $\rho((\Spin(V_\RR)_B)$
is contained in $\prod_{\hat{\sigma}\in\hat{\Sigma}}SU(V_{\hat{\sigma},\RR})$.
%*********
% End Hide
%*********
}

Conversely, let $g=\prod_{\hat{\sigma}\in\hat{\Sigma}}g_{\hat{\sigma}}$ be an element of $\prod_{\hat{\sigma}\in\hat{\Sigma}}SO_+(V_{\hat{\sigma},\RR})_{\eta_t}$. Then each $g_{\hat{\sigma}}$ lifts to an element of $\Spin(V_{\hat{\sigma},\RR})_{P_{\hat{\sigma}}}$, by \cite[Lemma 3.1.1]{markman-sixfolds}. 
%\ref{lemma-stabilizer-is-isomorphic-to-so-f}. 
Hence, $g$ lifts to $\prod_{\hat{\sigma}\in\hat{\Sigma}}\Spin(V_{\hat{\sigma},\RR})_{P_{\hat{\sigma}}}$, which maps to
$\Spin(V_\RR)_B$ via the right bottom arrow in (\ref{eq-diagram-of-spin-groups}).
\end{proof}

\begin{rem}
\label{rem-SU-V-hat-sigma-RR}
An element of $SO_+(V_{\hat{\sigma},\RR})$ leaves invariant the restriction of the $\CC$-valued hermitian form $\sigma\circ H_t$ to $V_{\hat{\sigma},\RR}$
if and only if it commutes with $\eta_t$. Hence, $SO_+(V_{\hat{\sigma},\RR})_{\eta_t}$ is equal to the subgroup of $SO_+(V_{\hat{\sigma},\RR})$ of elements leaving $\sigma\circ H_t$ invariant and whose restrictions to each of $V_{\sigma,\CC}$ and $V_{\bar{\sigma},\CC}$ have determinant $1$. 
Let $U_+(V_{\hat{\sigma},\RR})$ be the subgroup 
 of $SO_+(V_{\hat{\sigma},\RR})$ leaving the restriction of $\sigma\circ H_t$ invariant.
The determinant of the restriction of elements of $U_+(V_{\hat{\sigma},\RR})$ to $V_{\sigma,\RR}$ is a non-trivial character of $U_+(V_{\hat{\sigma},\RR})$.  
The subgroup $SO_+(V_{\hat{\sigma},\RR})_{\eta_t}$ thus deserves to be denoted $SU(V_{\hat{\sigma},\RR})$, since the imaginary parts of $\sigma\circ H_t$ and of $\bar{\sigma}\circ H_t$ defer only by a sign and their restrictions to $V_{\hat{\sigma},\RR}$ depend on $t$ only up to a real factor, and so $SO_+(V_{\hat{\sigma},\RR})_{\eta_t}$ is independent of $\sigma$ and of $t$.  
\end{rem}
%*******************************************************************************************************************
%
%*******************************************************************************************************************
\subsection{A period domain of abelian varieties of Weil type}
\label{sec-period-domain-for-abelian-varieties-of-Weil-type-CM-case}

%Recall that $V_{\sigma,\CC}=W_{T,\CC}\cap V_{\hat{\sigma},\CC}$, for any CM-type $T$ satisfying $T(\hat{\sigma})=\sigma$.
Fix a non-zero $t\in K_-$. 
Set $T_t:=\{\sigma\in \Sigma \ : \ Im(\sigma(t))>0\}$. $T_t$ is a CM-type and 
all CM-types arise\footnote{\label{footnote-on-CM-types}
 Choose $t_0\in K_-\setminus\{0\}$. It suffices to show that 
 $
 K_-\ni t\mapsto (\mbox{sign}(\hat{\sigma}(t/t_0))_{\hat{\sigma}\in\hat{\Sigma}}\in \{\pm 1\}^{e/2}
 $
 is surjective. This follows from the fact that the map $F\rightarrow \RR^{\hat{\Sigma}}$, given by $x\mapsto (\hat{\sigma}(x))_{\hat{\sigma}\in\hat{\Sigma}}$,
 contains a full lattice, for example the image of the ring of algebraic integers in $F$.
 } this way.
Given $I\in \rho(\Spin(V_\RR)_B)$, satisfying $I^2=-id_V$, set 
\begin{equation}
\label{eq-symmetric-bilinear-form-g_i-CM-case}
g_I(x,y) := \Xi_t(x,I(y))\stackrel{(\ref{eq-Xi-t})}{=}(\eta_t(x),I(y))_V.
\end{equation}
The automorphisms $I$ and $\eta_t$ of $V_\RR$ commute and both are anti-self-dual with respect to the symmetric bilinear pairing $(\bullet,\bullet)_V$ on $V_\RR$. Hence, $g_I$ is symmetric.
Let
\begin{equation}
\label{eq-omega-B-t}
\Omega_{B,t}\subset \rho(\Spin(V_\RR)_B)
\end{equation}
be the subset of elements $I$, such that $I$ is a complex structure on $V_\RR$, the eigenspaces $V_{\sigma,\CC}^{1,0}$ and $V_{\sigma,\CC}^{0,1}$ of the restriction of $I$ to $V_{\sigma,\CC}$ are both $d/2$-dimensional, for all $\sigma\in\Sigma$, and the symmetric bilinear form (\ref{eq-symmetric-bilinear-form-g_i-CM-case})
is positive definite. The subset $\Omega_{B,t}$ depends only on the CM-type $T_t$ associated to $t\in K_-\setminus\{0\}$, but the polarization $\Xi_t$ depends on $t\in K_-$. Denote $\Omega_{B,t}$ by $\Omega_{B,T_t}$ as well and set
\begin{equation}
\label{eq-omega-B}
\Omega_B:=\cup_{T\in\T_K}\Omega_{B,T}.
\end{equation}
Given two CM-types $T_1$ and $T_2$, define $\delta_{T_1,T_2}:V_\RR\rightarrow V_\RR$ as the isometry acting on $V_{\hat{\sigma},\RR}$ by multiplication by $-1$, if $T_1(\hat{\sigma})\neq T_2(\hat{\sigma})$ and by $1$ if $T_1(\hat{\sigma})= T_2(\hat{\sigma})$. The isometry $\delta_{T_1,T_2}$ belongs to $\Spin(V_\RR)_B$, by Lemma \ref{lemma-factorization-of-rho-of-Spin-V-B}, since $-id_{V_{\hat{\sigma},\RR}}$ belongs to $SO_+(V_{\hat{\sigma},\RR})_{\eta_t}$, as $d$ is even. Hence, given $I_1\in\Omega_{B,T_1}$, the complex structure $I_2:=\delta_{T_1,T_2}\circ I_1$ belongs to $\Omega_{B,T_2}$. 

\begin{rem}
The secant subspace $B$ does not determine the embedding $\eta:K\rightarrow \End_{Hdg}(H^1(X\times\hat{X},\QQ))$. Given $g\in Gal(K/\QQ)$, the embedding $\eta\circ g$ has the same secant subspace $B$. If we incorporate the embedding $\eta$ in the notation and denote $\Omega_{B,T}$ by $\Omega_{B,\eta,T}$, then $\Omega_{B,\eta\circ g,T}=\Omega_{B,\eta,g^{-1}(T)}$. Consider for example the involution $\iota\in Gal(K/F)\subset Gal(K/\QQ)$. 
If $I_{X\times\hat{X}}$ belongs to $\Omega_{B,\eta,T}$, then it belongs also to $\Omega_{B,\eta\circ\iota,\bar{T}}$ and 
$-I_{X\times\hat{X}}$ belongs to $\Omega_{B,\eta,\bar{T}}$. 
\end{rem}

\begin{lem}
The subset $\Omega_{B,t}$ of $\rho(\Spin(V_\RR)_B)$ is the image via $\rho$ of an adjoint orbit of $\Spin(V_\RR)_B$.
\end{lem}

\begin{proof}
%Recall that $\rho$ maps $\Spin(V_\RR)_B$ injectively into $SO_+(V_\RR)$. 
Step 1: (Reduction to $V_{\hat{\sigma},\RR}$). Set $T:=T_t$. On
each direct summand $V_{\hat{\sigma},\RR}$ we have the equalities $\rho(\Spin(V_{\hat{\sigma},\RR})_{P_{\hat{\sigma}}})=SO_+(V_{\hat{\sigma},\RR})_{\eta_t}$  and
 \[
\rho(\Spin(V_\RR)_B)= \prod_{\hat{\sigma}\in\hat{\Sigma}}\rho(\Spin(V_{\hat{\sigma},\RR})_{P_{\hat{\sigma}}}),
\]
established in the proof of Lemma \ref{lemma-factorization-of-rho-of-Spin-V-B}. Hence, 
$\Omega_{B,t}=\prod_{\hat{\sigma}\in\hat{\Sigma}}\Omega_{B,t,\hat{\sigma}}$, where 
$\Omega_{B,t,\hat{\sigma}}$ is the subset of $\rho(\Spin(V_{\hat{\sigma},\RR})_{P_{\hat{\sigma}}})$
consisting of $I$, such that $I^2=-id$, the eigenspaces of its restrictions to $V^{1,0}_{T(\hat{\sigma}),\CC}$ 
and $V^{1,0}_{\bar{T}(\hat{\sigma}),\CC}$ are $d/2$-dimensional\footnote{
Note that $V^{0,1}_{\bar{T}(\hat{\sigma}),\CC}$ is the complex conjugate of $V^{1,0}_{T(\hat{\sigma}),\CC}$ and so their dimensions are equal, for all CM-type and all $\hat{\sigma}\in\hat{\Sigma}$. Hence, if $V^{1,0}_{T(\hat{\sigma}),\CC}$ 
and $V^{1,0}_{\bar{T}(\hat{\sigma}),\CC}$ are $d/2$-dimensional, so are $V^{0,1}_{\bar{T}(\hat{\sigma}),\CC}$ and $V^{0,1}_{T(\hat{\sigma}),\CC}$.
}, and $(\eta_t(\bullet),I(\bullet))_{V_{\hat{\sigma}}}$
is positive definite. It suffices to show that $\Omega_{B,t,\hat{\sigma}}$ is an adjoint orbit in 
$\rho(\Spin(V_{\hat{\sigma},\RR})_{P_{\hat{\sigma}}})$.

Step 2: 
Let $I$ be an element in $\rho(\Spin(V_{\hat{\sigma},\RR})_{P_{\hat{\sigma}}})$ with $I^2=-id$, such that 
the eigenspaces  $V^{1,0}_{T(\hat{\sigma}),\CC}$ 
and $V^{1,0}_{\bar{T}(\hat{\sigma}),\CC}$ of its restrictions to $V_{T(\hat{\sigma})}$ are $d/2$-dimensional. We claim that 
$(\eta_t(\bullet),I(\bullet))_{V_{\hat{\sigma}}}$
is positive definite, if and only if the subspace $U^+_\CC:=V^{1,0}_{T(\hat{\sigma})}\oplus V^{0,1}_{\bar{T}(\hat{\sigma})}$ is positive definite with respect to 
the restriction of $(\bullet,\bullet)_{V_{\hat{\sigma}}}$ and 
$U^-_\CC:=V^{1,0}_{\bar{T}(\hat{\sigma})}\oplus V^{0,1}_{T(\hat{\sigma})}$ is negative definite.
Indeed, both subspaces are self conjugate, corresponding to subspaces $U^+$ and $U^-$ of $V_{\hat{\sigma},\RR}$,  and given $x,y\in V^{1,0}_{T(\hat{\sigma})}\oplus V^{0,1}_{\bar{T}(\hat{\sigma})}$, we have $\eta_tI(y)=T(\hat{\sigma})(t)iy$ and
\[
\Xi_t(x,I(y))=(\eta_t(x),I(y))_{V}=(x,-\eta_tI(y))_V=-T(\hat{\sigma})(t)i(x,y)=|T(\hat{\sigma})(t)|(x,y).
\] 
Given $x,y\in V^{1,0}_{\bar{T}(\hat{\sigma})}\oplus V^{0,1}_{T(\hat{\sigma})}$,  we have $\eta_tI(y)=-T(\hat{\sigma})(t)iy$ and
\[
\Xi_t(x,I(y))=(\eta_t(x),I(y))_{V}=(x,-\eta_tI(y))_V=T(\hat{\sigma})(t)i(x,y)=-|T(\hat{\sigma})(t)|(x,y).
\] 

Step 3: Let $I_1, I_2$ be elements of $\Omega_{B,t,\hat{\sigma}}$. Let $V_{T(\hat{\sigma})}^{(1,0),I_j}$ and $V_{T(\hat{\sigma})}^{(0,1),I_j}$
be the eigenspaces of the restriction of $I_j$ to $V_{T(\hat{\sigma})}$, $j=1,2$. 
Set $U^{+,I_j}_\CC:=V^{(1,0),I_j}_{T(\hat{\sigma})}\oplus V^{(0,1),I_j}_{\bar{T}(\hat{\sigma})}$ and define $U^{-,I_j}_\CC$ analogously and let $U^{+,I_j}$ and
$U^{-,I_j}$ be the corresponding real subspaces of $V_{\hat{\sigma},\RR}$, $j=1,2$.
We claim that there exists an element $g\in \rho(\Spin(V_{\hat{\sigma},\RR})_{P_{\hat{\sigma}}})$ satisfying
\begin{eqnarray}
\label{eq-g-maps-U-I-1-to-U-I-2}
g(U^{+,I_1})&=& U^{+,I_2},
%g(V_{T(\hat{\sigma})}^{(1,0),I_1}\oplus V_{\bar{T}(\hat{\sigma})}^{(0,1),I_1})&=& V_{T(\hat{\sigma})}^{(1,0),I_2}\oplus V_{\bar{T}(\hat{\sigma})}^{(0,1),I_2},
\\
\nonumber
g(U^{-,I_1})&=& U^{-,I_2}
%g(V_{\bar{T}(\hat{\sigma})}^{(1,0),I_1}\oplus V_{T(\hat{\sigma})}^{(0,1),I_1})&=& V_{\bar{T}(\hat{\sigma})}^{(1,0),I_2}\oplus V_{T(\hat{\sigma})}^{(0,1),I_2}.
\end{eqnarray}
$H_t$ restricts to $V_{\hat{\sigma},\RR}$ as a $K\otimes_{F,\hat{\sigma}}\RR$-valued hermitian form and its signature is the same whether we identify $K\otimes_{F,\hat{\sigma}}\RR$ with $\CC$ via $T(\hat{\sigma})$ or via $\bar{T}(\hat{\sigma})$.
The four subspaces $U^{+,I_j}$, $U^{-,I_j}$, $j=1,2$, 
%$V_{T(\hat{\sigma})}^{(1,0),I_j}\oplus V_{\bar{T}(\hat{\sigma})}^{(0,1),I_j}$ and $V_{\bar{T}(\hat{\sigma})}^{(1,0),I_j}\oplus V_{T(\hat{\sigma})}^{(0,1),I_j}$ 
are $K\otimes_{F,\hat{\sigma}}\RR$ subspaces of $V_{\hat{\sigma},\RR}$ and there exist isomorphisms of hermitian spaces $g_+:U^{+,I_1}\rightarrow U^{+,I_2}$
and $g_-:U^{-,I_1}\rightarrow U^{-,I_2}$, as the signatures of the restriction of $H_t$ are the same, by Step 2, hence we get a unitary automorphism $g$ of $(V_{\hat{\sigma},\RR},H_t)$ satisfying (\ref{eq-g-maps-U-I-1-to-U-I-2}).
The embedding $\eta:K\rightarrow \End(V_{\hat{\sigma},\RR})$ extends to an embedding of $K\otimes_{F,\hat{\sigma}}\RR$ and, possibly after composing $g$ with $\eta_s$
for a suitable $s\in K\otimes_{F,\hat{\sigma}}\RR$ of absolute value $1$, we may assume that $g$ belongs to $SU(V_{\hat{\sigma},\RR})=\Spin(V_{\hat{\sigma},\RR})_{P_{\hat{\sigma}}}$ (see Remark \ref{rem-SU-V-hat-sigma-RR}).

The definition of $H_t$ implies that $g$ is an isometry of $V_{\hat{\sigma},\RR}$ with respect to $(\bullet,\bullet)_{V_{\hat{\sigma}}}$
and it commutes with the restriction of $\eta_t$. By construction, we have
%$(\eta_t(x),I_1(y))=(\eta_t(g(x)),I_2(g(y)))$, for all $x,y\in V_{\hat{\sigma},\RR}$. Hence, 
$g(\eta_tI_1)g^{-1}=\eta_tI_2$. Combining with $g\eta_t=\eta_tg$ we get that $gI_1g^{-1}=I_2$. 

Step 4: We have seen that $\Omega_{B,t,\hat{\sigma}}$ is contained in a single adjoint orbit in $\Spin(V_{\hat{\sigma},\RR})_{P_{\hat{\sigma}}}.$
It remains to show the reverse inclusion, i.e., that $\Omega_{B,t,\hat{\sigma}}$ is invariant under the conjugation action of $\Spin(V_{\hat{\sigma},\RR})_{P_{\hat{\sigma}}}.$ The check is straightforward. 
\end{proof}

The connected component $\Omega_{B,t}$ may be regarded as 
a $\Spin(V_\RR)_B$-adjoint orbit, since $\rho$ maps $\Spin(V_\RR)_B$ injectively into $SO_+(V_\RR)$.
%The proof is identical to that of \cite[Lemma 4.0.2]{markman-sixfolds} 
%on each direct summand $V_{\hat{\sigma},\RR}$ using the equalities $\rho(\Spin(V_{\hat{\sigma},\RR})_{P_{\hat{\sigma}}})=SO_+(V_{\hat{\sigma},\RR})_{\eta_t}$ 
%and
 %\[
%\rho(\Spin(V_\RR)_B)= \prod_{\hat{\sigma}\in\hat{\Sigma}}\rho(\Spin(V_{\hat{\sigma},\RR})_{P_{\hat{\sigma}}}),
%\]
%established in the proof of Lemma \ref{lemma-factorization-of-rho-of-Spin-V-B}.
Let $A$ be the differentiable manifold underlying $X\times\hat{X}$. Given $I\in \Omega_{B,t}$ denote by $(A,I)$ the abelian variety with complex structure $I$.
The period domain $\Omega_{B,t}$ parametrizes polarized abelian  varieties of Weil type $((A,I),\eta,\Xi_t)$ with complex multiplication by $K$.
In Section \ref{sec-example-isometry-g-0} we will consider examples in which the complex structure $I_{X\times\hat{X}}$ belongs to $\Omega_{B,t}$ and
$\Omega_{B,t}$ parametrizes abelian varieties of split Weil type, i.e., for which $(V_\eta,H_t)$ is of split type 
 (see Definition \ref{def-split-type} and 
Lemma \ref{lemma-I-X-hat-X-belongs-to-Omega-B}).

\begin{lem}
\label{lem-if-Omega-B-t-is-of-split-type-then-all-triples-are-represented}
Assume that $(X\times\hat{X},\eta,\Xi_t)$ is of split Weil type.
Every $2n$-dimensional polarized abelian variety of split Weil type $(A',\eta',h)$, $\eta':K\rightarrow \End_{Hdg}(H^1(A',\QQ))$, is isogenous to $((A,I),\eta,\Xi_{t'})$, for some $I\in\Omega_{B,t'}$ and some non-zero $t'\in K_-$.
%where $T_{t'}=T_t$.
\end{lem}

\begin{proof}
There exists a unique hermitian form $H:H_1(A',\QQ)\times H_1(A',\QQ)\rightarrow K$ and a non-zero element $s\in K_-$, such that 
the polarization $h$ of   $(A',\eta',h)$ satisfies $h(x,y)=tr_{K/\QQ}(sH(x,y))$, for all $x,y\in H_1(A',\QQ)$, by 
\cite[Lem. 4.6]{deligne-milne}. 
There exists an isomorphism
$f:(H_1(A',\QQ),H) \rightarrow (V_\QQ,H_t)$, by Lemma \ref{lemma-invariants-of-hermitian-forms}, as both hermitian spaces have the same dimension over $K$ and both are of split type. 
Hence, $f$ conjugates $h(x,y)$ to 
\[
tr_{K/\QQ}(sH_t(x,y))=tr_{K/\QQ}(st(\eta_t(x),y)_{V_{\hat{\eta}}})=tr_{K/\QQ}((\eta_{st^2}(x),y)_{V_{\hat{\eta}}})=2\Xi_{st^2}(x,y).
\]
Let $I':H_1(A',\RR)\rightarrow H_1(A',\RR)$ be the complex structure of $A'$. Then the conjugate $I:=f\circ I'\circ f^{-1}:V_\RR\rightarrow V_\RR$ 
is such that $tr_{K/\QQ}((\eta_{st^2}(x),I(y))_{V_{\hat{\eta}}})$ is positive definite.
Hence, so is $tr_{K/\QQ}((\eta_{-st^2}(x),-I(y))_{V_{\hat{\eta}}})$. Now, the complex structure on $H_1(X\times\hat{X},\QQ)$ is related to minus the complex structure of $H^1(X\times\hat{X},\QQ)$ via the isomorphism $V_\RR\ni x\mapsto (x,\bullet)_V\in V_\RR^*$. 
So $(A',\eta',h)$ is isomorphic to $((A,-I),\eta,\Xi_{-st^2})$ and $-I$ belongs to $\Omega_{B,t'}$, where $t'=-st^2$.
\end{proof}

\begin{rem}
\label{rem-complete-family}
The statement and proof of Lemma \ref{lem-if-Omega-B-t-is-of-split-type-then-all-triples-are-represented} hold after replacing the assumption that 
$(X\times\hat{X},\eta,\Xi_t)$ and $(A',\eta',h)$ are of split Weil type by the assumption that their hermition forms have the same signatures and discriminant invariants (see Lemma \ref{lemma-invariants-of-hermitian-forms}).  In particular, $\Omega_{B,t}$ parametrizes a complete family of abelian varieties of Weil type, even if those are not of split type. 
\end{rem}

The connected component $\Omega_{B,t}$ is isomorphic to an open subset of 
$\prod_{\hat{\sigma}\in\hat{\Sigma}}Gr(d/2,V_{T_t(\hat{\sigma}),\CC})$
in the analytic topology, 
%where $\sigma\in\Sigma$ is a choice of one of the two characters which restricts to $F$ as $\hat{\sigma}$,
by Lemma \cite[4.0.1]{markman-sixfolds}
%\ref{lemma-coadjoint-orbit-embedds-as-open-subset-of-Grassmannian} 
applied to each factor separately.
 
 \begin{rem} 
 The period domains of abelian varieties of split Weil type are described in \cite{deligne-milne} as symmetric hermitian domains (see also \cite[Lemma 11.5.25]{charles-schnell}).
 The construction in \cite{deligne-milne} of the period domain of abelian varieties of split Weil type starts with the polarization $\Xi_t$ (called there the Riemann form) rather than the pairing $(\bullet,\bullet)_{V_{\hat{\eta}}}$. Note that we have the equalities $\Xi_t=tr_{F/\QQ}\circ\tilde{\Xi}_t$,
 \begin{eqnarray*}
 (x,y)_{V_{\hat{\eta}}} & = & \tilde{\Xi}_t(\eta_{t^{-1}}(x),y),
 \\
 H_t(x,y) &=& (-t^2)\tilde{\Xi}_t(\eta_{t^{-1}}(x),y)+t\tilde{\Xi}_t(x,y).
 \end{eqnarray*}
 Strictly speaking, $\tilde{\Xi}_t$ is an element of $\wedge_F^2\Hom_F(V_{\hat{\eta}},F)$ which is dual to $\wedge^2_FH^1(X\times\hat{X},\QQ)$.
 We regard $\Xi_t:=tr_{F/\QQ}\circ\tilde{\Xi}_t$ as a class in $\wedge^2_FH^1(X\times\hat{X},\QQ)\subset H^2(X\times\hat{X},\QQ)$ by identifying $V_{\hat{\eta}}$ with its dual via $(\bullet,\bullet)_{V_{\hat{\eta}}}$. The diffeomorphism induced by $(\bullet,\bullet)_V$, of the compact complex tori $X\times\hat{X}$ and its dual $(V_\RR/V,I)$,
 is an isomorphism of abelian varieties, for all complex structures $I$ in $\Omega_{B,t}$, since $\rho_I$ is an isometry of $V_\RR$, so the induced isomorphism
 $V_\RR\rightarrow V_\RR^*$ satisfies $I(v)\mapsto (Iv,\bullet)=(v,-I(\bullet))=I_{X\times\hat{X}}(v,\bullet)$.
 \end{rem}

%***************************************************************************
% 
%***************************************************************************
\section{Hodge-Weil classes, pure spinors, and Orlov's equivalence}
\label{sec-HW-classes-pure-spinors-and-Orlov-equivalence}
In section \ref{sec-grading-of-BB} we recall Orlov's derived equivalence $\Phi:D^b(X\times X)\rightarrow D^b(X\times\hat{X})$ and its $\Spin(V)$-equivariance properties. 
Given two coherent sheaves $F_i$, $i=1,2$, on $X$ with $ch(F_i)\in B$ and satisfying $\rank(\Phi(F_1\boxtimes F_2^\vee))\neq 0$, the class $\kappa(\Phi(F_1\boxtimes F_2^\vee))$ is $\Spin(V)_{\eta,B}$-invariant. Hence, $\kappa_{d/2}(\Phi(F_1\boxtimes F_2^\vee))$
belongs to $\A^d=\R^d\oplus HW(X\times\hat{X},\eta)$.
We introduce a grading on $B\otimes B$ in section \ref{sec-grading-of-BB}
and use it in section \ref{sec-HW-classes-from-pure-spinors} to give a sufficient condition for the class $\kappa_{d/2}(\Phi(F_1\boxtimes F_2^\vee))$ to project to a nonzero class in the direct summand $HW(X\times\hat{X},\eta)$ of $\A^d$. 

%***************************************************************************
% 
%***************************************************************************
\subsection{A grading of $B\otimes B$}
\label{sec-grading-of-BB}
We have the decomposition $B\otimes_\QQ\tilde{K}=\oplus_{T\in\T_K}\ell_T$, by Corollary \ref{cor-ell-T-are-linearly-independent}.
Given $T,T'\in\T_K$, set $T\cap T':=\{\hat{\sigma}\in\hat{\Sigma} \ : \ T(\hat{\sigma})=T'(\hat{\sigma})\}$.
Let $|T\cap T'|$ denote the cardinality of $T\cap T'$. 
Then $\dim_{\tilde{K}}(W_T\cap W_{T'})=\sum_{\hat{\sigma}\in T\cap T'}\dim_{\tilde{K}}(V_{T(\hat{\sigma})})=d|T\cap T'|$.

The subspace
\[
BB_{i,\tilde{K}} := \bigoplus_{\{(T,T')\in \T_K\times\T_K \ : \ |T\cap T'|=i\}}\ell_T\otimes\ell_{T'}
\]
of $S^+_\QQ\otimes_\QQ S^+_\QQ\otimes_\QQ \tilde{K}$ is $Gal(\tilde{K}/\QQ)$-invariant, hence of the form $BB_i\otimes_\QQ\tilde{K}$, for a subspace $BB_i$ of $S^+_\QQ\otimes_\QQ S^+_\QQ$.
The vector space $B\otimes_\QQ B$ decomposes as the direct sum
\[
B\otimes_\QQ B = \oplus_{i=0}^{e/2} BB_i.
\]
Note that the characters $\ell_{T_1}\otimes\ell_{T'_1}$ and $\ell_{T_2}\otimes\ell_{T'_2}$ of $\Spin(V_\QQ)_\eta$ are isomorphic, if and only if 
$T_1\cap T_1'=T_2\cap T_2'$, by Lemma \ref{lemma-tensor-product-of-two-ells}. Consequently, the characters $\ell_{T_1}\otimes\ell_{T'_1}$ and $\ell_{T_2}\otimes\ell_{T'_2}$ are different, if $|T_1\cap T'_1|\neq |T_2\cap T'_2|.$

The subspace $HW$ of $\wedge^dV_\QQ$ is the intersection $\wedge^dV_\QQ\cap[\oplus_{\sigma\in\Sigma}\wedge^dV_\sigma]$ in $\wedge^d(V\otimes_\ZZ\tilde{K})$. Denote by $\langle HW\rangle$ the intersection of $\wedge^*V_\QQ$ with the subalgebra of $\wedge^*(V\otimes_\ZZ\tilde{K})$
generated by $HW\otimes_\QQ\tilde{K}=\oplus_{\sigma\in\Sigma}\wedge^dV_\sigma$. Then $\langle HW\rangle$ is a graded subalgebra 
and
$\langle HW\rangle=\oplus_{k=0}^e\langle HW\rangle^{kd}$, where $\langle HW\rangle^{kd}$ consists of the rational points of the direct sum of all wedge products $\wedge_{i=1}^k\wedge^d V_{\sigma_i}$. Let $\tilde{\varphi}:S\otimes_\ZZ S\rightarrow \wedge^*V$ be the isomorphism given in (\ref{eq-tilde-varphi}) (and over $\ZZ$ in \cite[Eq. (2.3.2]{markman-sixfolds}).

Consider the increasing filtration $F^k(\wedge^*V):=\oplus_{i\leq k}\wedge^iV$ and the decreasing filtration $F_k(\wedge^*V):=\oplus_{i\geq k}\wedge^iV$.
Define $F^k(\wedge^*V_\QQ)$ and $F_k(\wedge^*V_\QQ)$ analogously. 

\begin{lem}
\label{lemma-on-lines-tensor-ell-T-ell-T-prime}
If $|T\cap T'|=k$, then the line $\tilde{\varphi}(\ell_T\otimes\ell_{T'})$ belongs to $F^{d(e-k)}(\wedge^*V_\CC)$ and it projects onto the line
$\wedge^{d(e-k)}[W_T+W_{T'}]$ in $\wedge^{d(e-k)}V_\CC$. The image of the composition 
\begin{equation}
\label{eq-composition-from-BB-1-to-wedge-de-d}
BB_1\RightArrowOf{\tilde{\varphi}}F^{d(e-1)}(\wedge^*V_\QQ)\rightarrow \wedge^{d(e-1)}V_\QQ
\end{equation}
is the $e$-dimensional subspace $\langle HW\rangle^{d(e-1)}$.
%which is the graded summand of degree $d(e-1)$ of the subalgebra $\langle HW\rangle$ of $\wedge^*V_\QQ$ generated by $HW$.
\end{lem}

Note that $\dim_\QQ(B)=2^{\dim_\QQ(F)}=2^{e/2}$ and so $\dim_\QQ(B\otimes_\QQ B)=2^e$. 
The character $\det_\sigma$ of $\Spin(V_{\tilde{K}})_\eta$ appears in $BB_k\otimes_\QQ\tilde{K}$, 
if and only if $k=1$ and its multiplicity in $BB_1\otimes_\QQ\tilde{K}$ is $2^{(\frac{e}{2}-1)}$. 
We see that $\dim_\QQ(BB_1)=e2^{(\frac{e}{2}-1)}$
%We have $\dim_\QQ(BB_1)=\frac{e}{2}2^{e/2}$ 
and so the composition (\ref{eq-composition-from-BB-1-to-wedge-de-d}) is injective, if and only if $e=2$.

\begin{proof}
If $d/2$ is even, then the line $\ell_\sigma\ell_{\bar{\sigma}}$ in the subspace $\Sym^2P_{\hat{\sigma}}$ of $S^+_{\hat{\sigma}}\otimes S^+_{\hat{\sigma}}$ is mapped via $\tilde{\varphi}_{\hat{\sigma}}$, given in (\ref{eq-tilde-varphi-hat-sigma}), to a line in $\wedge^*V_{\hat{\sigma}}$ which 
projects onto the top graded summand $\wedge^{2d}V_{\hat{\sigma}}$, by \cite[Lemma 2.3.2]{markman-sixfolds}.
%\ref{lemma-symmetric-or-alternating-product-of-two-pure-spinors-has-weight-2}. 
In that case the image $\tilde{\varphi}_{\hat{\sigma}}(\ell_\sigma\wedge \ell_{\bar{\sigma}})$, of the line $\ell_\sigma\wedge \ell_{\bar{\sigma}}$ in the subspace $\wedge^2P_{\hat{\sigma}}$ of $S^+_{\hat{\sigma}}\otimes S^+_{\hat{\sigma}}$, 
is contained in $F^{2d-2}(\wedge^*V_{\hat{\sigma}})$, by \cite[Lemma 2.3.2]{markman-sixfolds}.
%Lemma \ref{lemma-symmetric-or-alternating-product-of-two-pure-spinors-has-weight-2}. 

If $d/2$ is odd, then the line $\ell_\sigma\wedge \ell_{\bar{\sigma}}$ in the subspace $\wedge^2P_{\hat{\sigma}}$ of $S^+_{\hat{\sigma}}\otimes S^+_{\hat{\sigma}}$ is mapped via $\tilde{\varphi}_{\hat{\sigma}}$ to a line in $\wedge^*V_{\hat{\sigma}}$ which 
projects onto the top graded summand $\wedge^{2d}V_{\hat{\sigma}}$, by \cite[Lemma 2.3.2]{markman-sixfolds}.
%Lemma \ref{lemma-symmetric-or-alternating-product-of-two-pure-spinors-has-weight-2}. 
In that case the line $\tilde{\varphi}_{\hat{\sigma}}(\ell_\sigma \ell_{\bar{\sigma}})$ is contained in $F^{2d-2}(\wedge^*V_{\hat{\sigma}})$.

The line $\ell_\sigma\otimes\ell_\sigma$  is mapped via $\tilde{\varphi}_{\hat{\sigma}}$ to a line in $\wedge^*V_{\hat{\sigma}}$ which belongs to $F^d(\wedge^*V_{\hat{\sigma}})$, but not to $F^{d-1}(\wedge^*V_{\hat{\sigma}})$, and 
projects onto the line $\wedge^dV_\sigma$ in $\wedge^dV_{\hat{\sigma}}$, by  \cite[III.3.2]{chevalley}. 
Assume that $T,T'\in \T_K$ satisfy $|T\cap T'|=1$. Say $T(\hat{\sigma}_0)=T'(\hat{\sigma}_0)=\sigma_0$. Then
\[
\ell_T\otimes\ell_{T'}=(\ell_{\sigma_0}\otimes\ell_{\sigma_0})\otimes\bigotimes_{\{\hat{\sigma}\in\hat{\Sigma} \ : \ \hat{\sigma}\neq\hat{\sigma}_0\}}(\ell_{T(\hat{\sigma})}\otimes \ell_{T'(\hat{\sigma})}),
\]
by (\ref{eq-tensor-factorization-of-ell-T}).
The line $\tilde{\varphi}(\ell_T\otimes\ell_{T'})$ thus belongs to $F^{d(e-1)}(\wedge^*V_\CC)$ and projects onto the line\footnote{
It follows that the isomorphism $\phi:S^+_\CC\otimes S^+_\CC\rightarrow \wedge^*V_\CC$, given in (\ref{eq-phi}),
%of Lemma \ref{lemma-orlov-isomorphism-is-chevalley}, 
maps the line $\ell_T\otimes\ell_{T'}$
into a line in $F_d(\wedge^*V_\CC)$, which projects onto the line $\wedge^d(W_T\cap W_{T'})$ in $HW\otimes_\QQ\CC.$ 
}
$\wedge^{d(e-1)}[W_T+W_{T'}]$ in $\wedge^{d(e-1)}V_\CC$. The last statement follows also from the proof of \cite[III.3.3]{chevalley}.
We conclude that 
$
\tilde{\varphi}(BB_1)\otimes_\QQ\CC
$
is contained in $F^{d(e-1)}(\wedge^*V_\CC)$ and it projects onto
$
\oplus_{\sigma\in\Sigma}\wedge^{d(e-1)}[(V_\sigma)^\perp].
$
In particular, the image of the composition (\ref{eq-composition-from-BB-1-to-wedge-de-d})
%\begin{equation}
%\label{eq-composition-from-BB-1-to-wedge-de-d}
%BB_1\RightArrowOf{\tilde{\varphi}}F^{d(e-1)}(\wedge^*V_\QQ)\rightarrow \wedge^{d(e-1)}V_\QQ
%\end{equation}
is $e$-dimensional. 

If $|T\cap T'|=k$, then 
\[
\ell_T\otimes\ell_{T'}=\bigotimes_{\{\hat{\sigma}\in\hat{\Sigma} \ : \ T(\hat{\sigma})=T'(\hat{\sigma})\}}(\ell_{T(\hat{\sigma})}\otimes \ell_{T(\hat{\sigma})})
\otimes\bigotimes_{\{\hat{\sigma}\in\hat{\Sigma} \ : \ T(\hat{\sigma})\neq T'(\hat{\sigma})\}}(\ell_{T(\hat{\sigma})}\otimes \ell_{T'(\hat{\sigma})}),
\]
and the above argument applies verbatim to show that $\tilde{\varphi}(\ell_T\otimes\ell_{T'})$ belongs to $F^{d(e-k)}(\wedge^*V_\CC)$ and in projects onto the line $\wedge^{d(e-k)}[W_T+W_{T'}]$ in $\wedge^{d(e-k)}V_\CC$.
\end{proof}

\begin{rem}
Let $KB_1$ be the kernel of the composition (\ref{eq-composition-from-BB-1-to-wedge-de-d}).
Let us describe $KB_1$ explicitly.
It suffices to describe $KB_1\otimes_\QQ\tilde{K}$. Given $\sigma\in\Sigma$, let
$(\T_K\times \T_K)_{\sigma}$ be the set of ordered pairs $(T,T')\in \T_K\times \T_K$, such that $T\cap T'=\{\sigma\}$.
The equality $W_T+W_{T'}=\wedge^{d(e-1)}[\oplus_{\{\sigma'\in\Sigma \ : \sigma'\neq\bar{\sigma}\}}V_{\sigma'}]=(V_{\sigma})^\perp$ holds, 
for each pair $(T,T')\in(\T_K\times \T_K)_{\sigma}$. The line $\ell_T\otimes\ell_{T'}$ is canonically isomorphic to the line 
$\wedge^{d(e-1)}(V_{\sigma}^\perp)$, for each pair $(T,T')\in(\T_K\times \T_K)_{\sigma}$, by
Lemma \ref{lemma-on-lines-tensor-ell-T-ell-T-prime}. In particular, the lines $\ell_{T_1}\otimes\ell_{T'_1}$
and $\ell_{T_2}\otimes\ell_{T'_2}$ are canonically isomorphic,
for two pairs $(T_1,T'_1)$ and $(T_2,T'_2)$ in $(\T_K\times \T_K)_{\sigma}$. Choosing a non-zero vector $t_{\sigma}$ in $\wedge^{d(e-1)}[(V_{\sigma})^\perp]$, for each
$\sigma\in\Sigma$, 
we get a non-zero vector $\lambda_{(T,T')}$ in $\ell_T\otimes\ell_{T'}$ mapping to $t_{\sigma}$ via (\ref{eq-composition-from-BB-1-to-wedge-de-d}), for each 
$(T,T')$ with $|T\cap T'|=1$. The subspace $KB_1\otimes_\QQ\tilde{K}$ is explicitly given by
\begin{equation}
\label{eq-explicit-equations-for-KB-1}
KB_1\otimes_\QQ\tilde{K} \ = \ 
\left\{
\sum_{\{(T,T') \ : \ |T\cap T'|=1\}}c_{(T,T')}\lambda_{(T,T')} \ : \
\sum_{(T,T')\in(\T_K\times \T_K)_{\sigma}}c_{(T,T')}=0, \ \ \forall \sigma\in\Sigma
\right\},
\end{equation}
where $c_{(T,T')}\in \tilde{K}$. 
\EndProof
\end{rem}

%*************
% Hide
%*************
\hide{
\begin{question}
Let $m_{g_0}(1)\in[\wedge^{ev}_FH^1(X,\QQ)]\otimes_FK$ be the pure spinor of $W$.
Then $(id\otimes T(\hat{\sigma}))(m_{g_0}(1))$ is an element of 
$[\wedge^{ev}H^1_{\hat{\sigma}}(X)]\otimes_{\tilde{F}}\tilde{K}$ and 
\[
\otimes_{\hat{\sigma}\in\hat{\Sigma}}(id\otimes T(\hat{\sigma}))(m_{g_0}(1))
\]
is an explicit class $\lambda_T$ in $\ell_T$. 
%Hence, $\lambda_T\otimes\lambda_{T'}$ is an explicit class in $\ell_T\otimes\ell_{T'}$. 
Can we take $\lambda_{(T,T')}$ in the above remark to be $\lambda_T\otimes\lambda_{T'}$?
I.e., are the two cosets $\tilde{\varphi}(\lambda_{T_1}\otimes\lambda_{T'_1})+F^{d(e-1)}(\wedge^*V_\CC)$ and
$\tilde{\varphi}(\lambda_{T_2}\otimes\lambda_{T'_2})+F^{d(e-1)}(\wedge^*V_\CC)$
equal, for any two pairs $(T_i,T'_i)\in(\T_K\times\T_K)_\sigma$, $i=1,2$?
\end{question}
}
%*************
% End Hide
%*************

Let $\P$ be the normalized Poincar\'{e} line bundle over $\hat{X}\times X$. Let $\mu:X\times X\rightarrow X\times X$ be the automorphism $\mu(x,y)=(x+y,y)$.
Let $\Phi_\P:D^b(\hat{X})\rightarrow D^b(X)$ be the equivalence with Fourier-Mukai kernel $\P$. We get the equivalence
$id\times\Phi_\P:D^b(X\times\hat{X})\rightarrow D^b(X\times X)$ with Fourier-Mukai kernel $\pi_{1,3}^*(\StructureSheaf{\Delta})\otimes\pi_{2,4}^*\P$,
where $\pi_{i,j}$ is the projection from $X\times\hat{X}\times X\times X$ onto the product of the $i$-th and $j$-th factors and $\Delta$ is the diagonal in $X\times X$. Orlov's equivalence
\begin{equation}
\label{eq-Orlov-derived-equivalence}
\Phi:D^b(X\times X)\rightarrow D^b(X\times\hat{X})
\end{equation}
is the inverse of $\mu_*\circ (id\times \Phi_\P):D^b(X\times\hat{X})\rightarrow D^b(X\times X).$
Let 
\begin{equation}
\label{eq-phi}
\phi:H^*(X\times X,\ZZ)\rightarrow H^*(X\times\hat{X},\ZZ)
\end{equation} 
be the isomorphism induced by the cohomological action of $\Phi$. Set
\[
\tilde{\phi}:=\exp(-c_1(\P)/2)\cup\phi\circ(id\otimes\tau):H^*(X\times X,\QQ)\rightarrow H^*(X\times\hat{X},\QQ),
\]
where 
%$\phi$ is given in (\ref{eq-phi}) and 
\begin{equation}
\label{eq-tau}
\tau:H^*(X,\ZZ)\rightarrow H^*(X,\ZZ)
\end{equation} 
acts on $H^i(X,\ZZ)$ by multiplication by $(-1)^{i(i-1)/2}$.

\begin{prop}
\label{prop-tilde-phi-is-spin-equivariant}
\cite[Prop. 1.3.1]{markman-sixfolds}
The linear transformation $\tilde{\phi}$ is an isomorphism of two $\Spin(V_\QQ)$-representations, the tensor square $S_\QQ\otimes S_\QQ$ of the spin representation and $\wedge^*V_\QQ$
\end{prop}.

Note the inclusion $\tilde{\phi}(B\otimes B)\subset \A$, where $\A:=(\wedge^*V_\QQ)^{\Spin(V)_{\eta,B}}$.
% where $\tilde{\phi}$ is given in (\ref{eq-tilde-phi}).

\begin{lem}
\label{lemma-on-the-phi-image-of-lines-tensor-ell-T-ell-T-prime}
If $|T\cap T'|=k$, then the line $\phi(\ell_T\otimes\tau(\ell_{T'}))$ belongs
%\footnote{Should it be the line $(\phi\circ (id\otimes\tau))(\ell_T\otimes\ell_{T'})$?}
 to $F_{dk}(\wedge^*V_\CC)$ and it projects onto the line
$\wedge^{dk}[W_T\cap W_{T'}]$ in $\wedge^{dk}V_\CC$. The image of the composition 
\begin{equation}
\label{eq-composition-from-BB-1-to-wedge-d}
BB_1\RightArrowOf{\phi\circ(id\otimes\tau)}F_d(\wedge^*V_\QQ)\rightarrow \wedge^dV_\QQ
\end{equation}
is $HW$.
%an $e$-dimensional subspace, which is the graded summand of degree $d$ of $\langle HW\rangle$.
\end{lem}

\begin{proof}
Let $\phi_\P:H^*(X,\ZZ)\rightarrow H^*(\hat{X},\ZZ)$ be the isomorphism induced by $ch(\P)$.
The statement follows from Lemma \ref{lemma-on-lines-tensor-ell-T-ell-T-prime}, by the equality $\phi\circ(id\otimes\tau)=(\phi_\P\otimes\phi_\P^{-1})\circ\tilde{\varphi}$ \cite[Lemma 6.3.1]{markman-sixfolds} and the fact that 
$\phi_\P\otimes\phi_\P^{-1}$ restricts to each graded summand of $\wedge^*V_\QQ$ 
as Poincar\'{e} duality, up to sign \cite[Lemma 6.3.2]{markman-sixfolds}.
%Lemma \ref{lemma-phi-P-psi-P-inverse-is-PD-up-to-sign}. 
We claim that Poincar\'{e} duality followed by the cohomological action of the automorphism interchanging the factors $X$ and $\hat{X}$ results in an automorphism of $\wedge^*V$ which
maps $\langle HW\rangle^{d(e-1)}$ to $\langle HW\rangle^d$.
It suffices to show that it maps $\wedge^{d(e-1)}[W_T+ W_{T'}]$ to $\wedge^d[W_T\cap W_{T'}]$. 
Indeed, Poincar\'{e} duality maps the top exterior power $\wedge^t_{\tilde{K}}Z$ of a $t$-dimensional subspace $Z$ of $V_{\tilde{K}}$ to the top exterior power 
$\wedge^{de-t}_{\tilde{K}}(\ann(Z))$, where $\ann(Z)$ is the annihilator of $Z$ in $V_{\tilde{K}}^*=\Hom(V_{\tilde{K}},\tilde{K})$ with respect to the evaluation pairing. Interchanging the factors of $X$ and $\hat{X}$ induces on the level of first cohomologies the isomorphism $V\cong V^*$ induced by the pairing $(\bullet,\bullet)_V$, which maps $\ann(Z)$ to the orthogonal complement $Z^\perp$ with respect to $(\bullet,\bullet)_V$.
Now note that 
$W_T\cap W_{T'}$ is the orthogonal complement of $W_T+ W_{T'}$ with respect to $(\bullet,\bullet)_V$ and 
given $\sigma\in\Sigma$ restricting to $\hat{\sigma}\in\hat{\Sigma}$, the subspace $V_\sigma$ is the orthogonal complement of $V_\sigma\oplus[\oplus_{\{\hat{\sigma}'\in\hat{\Sigma} \ : \ \hat{\sigma}'\neq\hat{\sigma}\}}V_{\hat{\sigma}'}]$.
The subspace $\langle HW\rangle^d$ is the direct sum of the top exterior powers of the former subspaces and 
$\langle HW\rangle^{d(e-1}$ is the direct sum of the top exterior powers of the latter subspaces.
\end{proof}

%**************
% Hide
%**************
\hide{
Consider the decreasing filtration $F_k(B\otimes_\QQ B):=\oplus_{i\geq k}BB_i$
and the increasing filtration $F^k\wedge^*V_\QQ:=\oplus_{i\leq k}\wedge^iV_\QQ$.
% and the increasing filtration $F^k(B\otimes_\QQ B):=\oplus_{i\leq k}BB_i$.
The isomorphism $\tilde{\varphi}:S_\QQ\otimes_\QQ S_\QQ\rightarrow \wedge^*V_\QQ$, given in (\ref{eq-tilde-varphi}),
maps $F_k(B\otimes_\QQ B)$ into $F^{d(e-k)}(\wedge^*V_\QQ)$, but the image of $F_k(B\otimes_\QQ B)$ is not contained in $F^{d(e-k)-1}(\wedge^*V_\QQ)$, by \cite[III.3.3]{chevalley}. 
We get the induced homomorphism $BB_k\rightarrow \wedge^{d(e-k)}V_\QQ$. 
The proof of \cite[III.3.3]{chevalley} shows furthermore that the image of $\ell_T\otimes\ell_{T'}$ in $\wedge^{d(e-k)}V_\QQ$ is the top exterior power of 
the subspace $W_T+ W_{T'}$ of $V_{\tilde{K}}$.
The associated homomorphism $BB_k\otimes_\QQ\tilde{K}\rightarrow \wedge^{d(e-k)}V_{\tilde{K}}$
is furthermore $\Spin(V_{\tilde{K}})_\eta$-equivariant. 
%Choose an ordering of the set $\T_K$. 
The image of the induced homomorphism 
$BB_1\otimes_\QQ\tilde{K}\rightarrow \wedge^{d(e-1)}V_{\tilde{K}}$ is thus 
\[
\bigoplus_{\sigma\in\Sigma}
\wedge^{d(e-1)}\left(\bigoplus_{\{\sigma' \ : \ \sigma'\neq \sigma\}}V_{\sigma'}
\right).
%=HW\otimes_\QQ\tilde{K}
\] 
%and $\tilde{\varphi}(BB_{e-1})=HW$. 
The isomorphism $\phi\circ (id\otimes\tau)$, $\phi$ as in Lemma \ref{lemma-orlov-isomorphism-is-chevalley}, maps 
$F_k(B\otimes_\QQ B)$ into $F_{dk}(\wedge^*V_\QQ)$ and the induced homomorphism on the graded summands maps 
$BB_1$ onto $HW$.
%since the characters $\ell_T\otimes \ell_{T'}$, . 

%**************
% End Hide
%**************
}

%Note that $\dim_\QQ(B)=2^{\dim_\QQ(F)}=2^{e/2}$ and so $\dim_\QQ(B\otimes_\QQ B)=2^e$. 
%The character $\det_\sigma$ of $\Spin(V_{\tilde{K}})_\eta$ appears in $BB_k\otimes_\QQ\tilde{K}$, 
%if and only if $k=1$ and its multiplicity in $BB_1\otimes_\QQ\tilde{K}$ is $2^{(\frac{e}{2}-1)}$. 
%We see that $\dim_\QQ(BB_1)=e2^{(\frac{e}{2}-1)}$.
%The character $\otimes_{\{\hat{\sigma} \ : \ \hat{\sigma}\neq \hat{\sigma}_0\}}\det_{T(\hat{\sigma})}$
%appears in $BB_k\otimes_\QQ\tilde{K}$, if and only if $k=\frac{e}{2}-1$ and it appears in $BB_{(\frac{e}{2}-1)}\otimes_\QQ\tilde{K}$ with multiplicity $2$.
%There are $(e/2)2^{(\frac{e}{2}-1)}$ characters of this type. 
%We see that $\dim_\QQ(BB_{(\frac{e}{2}-1)})=e2^{(\frac{e}{2}-1)}$.
The character $\det_{\sigma_1}\otimes \cdots \otimes \det_{\sigma_k}$, with $\hat{\sigma}_1, \dots, \hat{\sigma}_k$ distinct in $\hat{\Sigma}$,
appears only in $BB_k\otimes_\QQ\tilde{K}$ and it appears in $BB_k\otimes_\QQ\tilde{K}$ with multiplicity $2^{(\frac{e}{2}-k)}$.
Denote by 
\[
BB_{T\cap T'}
\] 
the subspace of $BB_k\otimes_\QQ\tilde{K}$, $k=|T\cap T'|$, such that $BB_{T\cap T'}$ is the direct sum of all 
characters of $\Spin(V_{\tilde{K}})_\eta$ isomorphic to $\otimes_{\{\hat{\sigma}\in T\cap T'\}}\det_{T(\hat{\sigma})}$. 
There are $2^k\Choose{e/2}{k}$ pairwise non-isomorphic such characters, so the dimension of $BB_k$ is $\Choose{e/2}{k}2^{e/2}$.
%***************************************************************************
% 
%***************************************************************************
\subsection{Hodge-Weil classes from pure spinors} 
\label{sec-HW-classes-from-pure-spinors}
The $\QQ$-subalgebra  $\langle HW\rangle$ of $\wedge^*V_\QQ$ in Lemma \ref{lemma-on-lines-tensor-ell-T-ell-T-prime}
%generated by the $e$-dimensional subspace $HW$ of $\wedge^dV_\QQ$. 
is a graded subalgebra and we set $\langle HW\rangle^k:=\langle HW\rangle\cap \wedge^kV_\QQ.$ 
The graded summand $\langle HW\rangle^k$ vanishes, if $d\not| \ k$, and
\[
\dim_\QQ \langle HW\rangle^{dk}=\Choose{e}{k},
\]
for $0\leq k\leq e$.
The multiplicity of $\det_\sigma$ in $\langle HW\rangle^{dk}\otimes_\QQ\tilde{K}$ is $0$, if $k$ is even, and it is $\Choose{(e/2)-1}{i}$, if $k=2i+1$, for $0\leq i\leq \frac{e}{2}-1$. Hence, the multiplicity of $\det_\sigma$ in $\langle HW\rangle\otimes_\QQ\tilde{K}$ is equal to its multiplicity in $B\otimes_\QQ B\otimes_\QQ\tilde{K}$. 
%Indeed, we have the following.

%**************
% Hide
%**************
\hide{
\begin{lem}
The isomorphism
\[
\tilde{\phi}:=\exp\left(\frac{1}{2}c_1(\P)\right)\cup \phi\circ (id\otimes \tau):H^*(X\times X,\QQ)\rightarrow H^*(X\times\hat{X},\QQ),
\]
which is  $\Spin(V)$-equivariant with respect to $(m\otimes m^\dagger,\rho)$ 
by Proposition \ref{prop-extension-class-of-decreasing-filtration-of-spin-V-representations}, 
%maps the $2^e$-dimensional subspace 
%$B\otimes_\QQ B$ of $S\otimes_\QQ S$ isomorphically onto 
%$\langle HW\rangle$. 
%Furthermore, 
satisfies
\begin{equation}
\label{equality-in-degree-de-over-2}
\tilde{\phi}(BB_{e/2})\subset \langle HW\rangle^{de/2},
\end{equation}
where the right hand side is the graded summand of $\langle HW\rangle$ in degree $de/2$. 
\end{lem}

\begin{proof}
The isomorphism $\tilde{\phi}$ extends to a $\Spin(V_{\tilde{K}})_\eta$-equivariant one from $S_\QQ\otimes_\QQ S_\QQ\otimes\tilde{K}$ onto
$\wedge^*V_{\tilde{K}}$. The subspace $BB_{e/2}\otimes_\QQ\tilde{K}$ is the sum of all 
characters of $\Spin(V_{\tilde{K}})_\eta$ isomorphic to $\ell_T\otimes\ell_T$, for some $T\in\hat{\Sigma}$.
The subspace $\langle HW\rangle^{de/2}\otimes_\QQ\tilde{K}$ of $\wedge^*V_{\tilde{K}}$ is the sum of all characters of $\Spin(V_{\tilde{K}})_\eta$ isomorphic to
tensor products $\otimes_{\sigma\in \Sigma'}\det_\sigma$, for some subset $\Sigma'\subset\Sigma$ of cardinality $e/2$, which projects onto (??? no ???) $\hat{\Sigma}$ under the restriction map $\Sigma\rightarrow\hat{\Sigma}$. The two sets of characters are equal, 
by Lemma \ref{lemma-tensor-product-of-two-ells}.
\end{proof}
%*********
% End Hide
%*********
}

%\begin{rem}
%Additional grading on $\wedge^*V_{\tilde{K}}$ is given by the action via $T\circ \eta$ of $K^\times$. 
%According to this grading $\langle HW\rangle\otimes_\QQ\tilde{K}$
%is the sum of all $d$-th powers of characters of $K^\times$.
%\end{rem}

Let $c$ be a class in $B\otimes_\QQ B$. 
Write $c=\sum_i c_i$, with $c_i$ in $BB_i$. Set $\check{\phi}:=\phi\circ(id\otimes\tau).$
Let $\check{\phi}(c)_i$ be the graded summand of $\check{\phi}(c)$ in $\wedge^{2i}V_\QQ$. 
Set $r:=\check{\phi}(c)_0$, considered as a rational number via the isomorphism $\wedge^0V_\QQ\cong \QQ$. 
Set $\kappa(\check{\phi}(c)):=\check{\phi}(c)\exp\left(-\check{\phi}(c)_1/r\right)$ whenever $r\neq 0$.

\begin{prop}
\label{prop-kappa-class-of-image-of-secant-class-yields-a-HW-class}
Assume that $c_1$ does not belong to $KB_1$, $r\neq 0$, and $d>2$. 
\begin{enumerate}
\item
\label{lemma-item-kappa-phi-c-is-invariant}
The class $\kappa(\check{\phi}(c))$ is $\Spin(V)_{\eta,B}$-invariant. 
\item
\label{lemma-item-kappa-phi-c-0-is-invariant}
The class $\kappa(\check{\phi}(c_0))$ belongs to $\R$. 
\item
\label{lemma-item-difference-projects-to-HW}
The equality $\check{\phi}(c)_0=\check{\phi}(c_0)_0$ holds. 
The difference
$\kappa(\check{\phi}(c))-\kappa(\check{\phi}(c_0))$ belongs to $F_d(\wedge^*V_\QQ)$ and its projection to $\wedge^dV_\QQ$ is a non-zero element of $HW$.
\end{enumerate}
\end{prop}

\begin{proof}
(\ref{lemma-item-kappa-phi-c-is-invariant}) 
%The $\Spin(V)_{\eta,B}$-invariance of the class $\kappa(\check{\phi}(c))$ follows from Proposition \ref{prop-tilde-phi-is-spin-equivariant}, 
%by the same argument proving \cite[Corollary 1.3.2]{markman-sixfolds}.
%We include the argument here for the convenience of the reader. 
The class 
$\check{\phi}(c)$ is $\Spin(V)_{\eta,B}$-invariant, by Proposition \ref{prop-tilde-phi-is-spin-equivariant}.
The group $\Spin(V_\QQ)$ acts on $\wedge^*V_\QQ$ via graded algebra automorphisms. Hence, the graded summands of $\check{\phi}(c)$ are $\Spin(V)_{\eta,B}$-invariant. Thus, $\check{\phi}(c)_1$ and consequently also $\exp\left(-\check{\phi}(c)_1/r\right)$ are $\Spin(V)_{\eta,B}$-invariant. The class $\kappa(\check{\phi}(c))$, being the product of two $\Spin(V)_{\eta,B}$-invariant factors, is thus $\Spin(V)_{\eta,B}$-invariant.

(\ref{lemma-item-kappa-phi-c-0-is-invariant})
The class $\check{\phi}(c_k)$ belongs to $F_{dk}(\wedge^*V_\QQ)$, by Lemma \ref{lemma-on-the-phi-image-of-lines-tensor-ell-T-ell-T-prime}, and $d>2$, by assumption. Hence, $\check{\phi}(c)_0=\check{\phi}(c_0)_0$ and 
$\check{\phi}(c)_1=\check{\phi}(c_0)_1$. 
We thus have
\[
\kappa(\check{\phi}(c))=\kappa(\check{\phi}(c_0))+\check{\phi}(c-c_0)\exp\left(-\check{\phi}(c)_1/r\right).
\]
The class $\kappa(\check{\phi}(c_0))$ is $\Spin(V)_\eta$-invariant, by the same argument proving 
\cite[Corollary 1.3.2]{markman-sixfolds},
%Corollary \ref{cor-kappa-class-is-Spin-V-P-invariant},
and is hence in $\R$, by Proposition \ref{prop-generator-for-the-invariant-subalgebra}(\ref{lemma-item-generatyors-for-the-invariant-subalgebra-R}). 

(\ref{lemma-item-difference-projects-to-HW})
The class $\check{\phi}(c-c_0)\exp\left(-\check{\phi}(c)_1/r\right)$ belongs to $F_d(\wedge^*V_\QQ)$
and its projection to $\wedge^dV_\QQ$ is equal to that of $\check{\phi}(c-c_0)$, which in turn is equal to the projection of $\check{\phi}(c_1)$. 
Now, the projection of $\check{\phi}(c_1)$ to $\wedge^dV_\QQ$ is a non-zero element of $HW$, by Lemma \ref{lemma-on-the-phi-image-of-lines-tensor-ell-T-ell-T-prime}.
\end{proof}

\begin{rem}
There exist classes $c$ of the form $c=\alpha\otimes\beta$, with $\alpha,\beta\in B$, satisfying the assumptions of Proposition \ref{prop-kappa-class-of-image-of-secant-class-yields-a-HW-class}. Indeed, any subspace containing classes of the form $\alpha\otimes\beta$, $\alpha,\beta\in B$, 
 is the whole of $B\otimes B$ and so for generic choices of $\alpha$ and $\beta$ we have that $c_1\not\in KB_1$ and $\check{\phi}(c)_0\neq 0$. One would like to choose $\alpha$ and $\beta$ to be Chern classes of objects in $D^b(X)$. 
\end{rem}

The following Corollary is an immediate consequence of Proposition \ref{prop-kappa-class-of-image-of-secant-class-yields-a-HW-class}. It is the reduction of the algebraicity of the Weil classes on abelian varieties with complex multiplication to the variational Hodge conjecture. 
Let $(A,\eta',h')$ be a polarized abelian variety of Weil type parametrized by the same connected component of the period domain $\Omega_B$ corresponding to 
$(X\times\hat{X},\eta,h)$, where $h=\Xi_t$ for some $t\in K_-$, given in (\ref{eq-omega-B-t}).  

\begin{cor}
\label{corollary-VHC}
Assume that the class $c$ is algebraic. 
If the flat deformation of the algebraic class $\kappa_{d/2}(\check{\phi}(c))$ in $H^d(X\times\hat{X},\QQ)$ remains algebraic in $H^d(A,\QQ)$, then every class in $HW(A,\eta')$ is algebraic.
\end{cor}

%*************
% Hide
%*************
\hide{
Let $c$ be a class in $B\otimes_\QQ B$. 
Write
$c=\sum_i c_i$, with $c_i$ in $BB_i$ and assume that $c_1$ does not belong to $KB_1$.

\begin{lem}
There exists a class $\delta\in \R$, such that 
$\tilde{\varphi}(c)-\delta$ belongs to $F^{d(e-1)}(\wedge^*V_\QQ)$ and its projection to $\wedge^{d(e-1)}V_\QQ$ is a non-zero class in 
$\langle HW\rangle^{d(e-1)}$. 
\end{lem}

\begin{proof}
The summand $c_0$ spans a trivial $\Spin(V_\QQ)_\eta$-character.
Hence, there exists a class $\alpha\in \R$, such that the $\Spin(V_\QQ)_\eta$-invariant class 
$\tilde{\varphi}(c_0)-\alpha$ belongs to $F^{d(e-1)}(\wedge^*V_\QQ)$, by Proposition \ref{prop-generator-for-the-invariant-subalgebra} (??? the lemma is in terms of the $\rho$-action, but $\tilde{\varphi}$ is not equivariant with respect to this action ???).
Let $\beta$ be the projection of $\tilde{\varphi}(c_0)-\alpha$ to $\wedge^{d(e-1)}V_\QQ$.
The class $\beta$ is $\Spin(V_{\tilde{K}})_\eta$-invariant and so in $\R$, by Lemma \ref{lemma-generators-of-the-module-of-Spin-V-B-invariant-classes}. 
The summand $\tilde{\varphi}(c_i)$ belongs to $F^{d(e-i)}(\wedge^*V_\QQ)$. Hence, $\tilde{\varphi}(c-c_0)$ belongs to $F^{d(e-1)}(\wedge^*V_\QQ)$
and its image in $\wedge^{d(e-1)}V_\QQ$ is equal to the image of $\tilde{\varphi}(c_1)$, 
which is an element $\gamma$ of $\langle HW\rangle^{d(e-1)}$, by Lemma \ref{lemma-on-lines-tensor-ell-T-ell-T-prime}. The element $\gamma$ 
%of $\langle HW\rangle^{d(e-1)}$ 
is non-zero,  by assumption. The 
vector $\tilde{\varphi}(c)-\alpha=(\tilde{\varphi}(c_0)-\alpha)+\tilde{\varphi}(c-c_0))$ is thus in $F^{d(e-1)}(\wedge^*V_\QQ)$ and it projects to $\wedge^{d(e-1)}V_\QQ$ as 
the sum of a $\Spin(V_{\tilde{K}})_\eta$-invariant class $\beta$ and a non-zero class $\gamma$ in $\langle HW\rangle^{d(e-1)}$.
The statement follows by setting $\delta=\alpha+\beta$.
\end{proof}
%*************
% End Hide
%*************
}

%Let $T$ and $T'$ be CM-types, such that the set 
%$\{\hat{\sigma}\in \hat{\Sigma} \ : \ T(\hat{\sigma})= T'(\hat{\sigma})\}$ consists of a single character $\hat{\sigma}_{T,T'}$. 
%Denote by $\sigma_{T,T'}\in\Sigma$ the common value of $T$ and $T'$ at $\hat{\sigma}_{T,T'}$. 
%The intersection $W_{T,\CC}\cap W_{T',\CC}$ is $V_{\sigma_{T,T'}}$ and their sum
%$W_{T,\CC}+ W_{T',\CC}$ is the direct sum of all $V_\sigma$, for $\sigma\neq \iota(\sigma_{T,T'})$. Now use \cite[III.3.3]{chevalley} 
%to recover $\wedge^{top}V_{\iota(\sigma_{T,T'})}$ from the tensor product of the two pure spinors $\ell_T\otimes \ell_{T'}$.

%Show  that for a generic pair of points $b_1, b_2\in B$, the image of 
%$b_1\otimes b_2$ in $\wedge^*V_\QQ$ via Orlov's cohomological isomorphism, normalized as in the definition of the $\kappa$ class, 
%has graded component in the sum of $\wedge^{d}V_\QQ$, and the subspace spanned by polynomials in divisor class, 
%and the projection $\kappa(b_1\otimes b_2)$ to $HW$ is non-zero. Then choose coherent sheaves $F_i$ with $ch(F_i)=b_i$ and 
%deduce that the classes in $HW$ are algebraic. Then deform 
%the image of $F_1\boxtimes F_2$ via Orlov's derived equivalence, to all abelian $2n$-folds with complex multiplication by $K$ deformation equivalent 
%to $(X\times\hat{X},\eta)$

%*******************************************************************************************************************
%
%*******************************************************************************************************************
\section{Examples}
%\begin{example}
\label{sec-example-isometry-g-0}

Let $K$ be a CM-field and $F$ its maximal totally real subfield. In section \ref{sec-example-of-K-secant-for-CM-field}
we construct a $K$-secant linear subspace $B\subset S^+_\QQ:=H^{ev}(X,\QQ)$, spanned by Hodge classes, for every simple abelian variety $X$ with an embedding $\hat{\eta}:F\rightarrow\End_\QQ(X)$.
Let $\tilde{K}$ be the Galois closure of $K$ over $\QQ$ and $\tilde{F}$ that of $F$. 
We have the short exact sequence
\[
0\rightarrow (\ZZ/2\ZZ)^\nu\rightarrow \mbox{Gal}(\tilde{K}/\QQ)\rightarrow \mbox{Gal}(\tilde{F}/\QQ)\rightarrow 0,
\] 
where $1\leq \nu\leq e/2:=[F:\QQ]$, by \cite[Sec. 1.1]{dodson}. 
In section \ref{sec-example-of-a-class-c} we specialize to the case where $\nu=1$ and the above short exact sequence splits. 
This is the case when $K=F(\sqrt{-q})$, where $q$ is a positive rational number.
In this case we provide an example of a class $c=\alpha_0\otimes\beta$ in $B\otimes B$, which satisfies the assumptions of Proposition \ref{prop-kappa-class-of-image-of-secant-class-yields-a-HW-class} and is hence mapped to a non-zero class in $HW(X\times \hat{X},\eta)$ 
(see Lemma \ref{lemma-a-criterion-for-alpha-0-tensor-beta-to-map-via-orlov-equivalence-to-a-class-inducing-non-zero-HW-class}).
Here $\alpha_0$ is a class in the $\QQ(\sqrt{-q})$-secant line associated to the pure spinors $\exp(\sqrt{-q}\Theta)$ and its complex conjugate as in \cite[Equation (1.2.4)]{markman-sixfolds}.
%(\ref{eq-P-introduction}). 
If $X$ is a jacobian and $\Theta$ is its principal polarization, a choice of $\alpha_0=ch(E)$, where $E$ is a coherent sheaf, is discussed in \cite[Example 8.2.3]{markman-sixfolds}.
%\ref{example-secant-sheaf-on-jacobian-of-genus-4-curve}. 
The class $\beta\in B$ satisfies an explicit genericity assumption 
(see Lemma \ref{lemma-a-criterion-for-alpha-0-tensor-beta-to-map-via-orlov-equivalence-to-a-class-inducing-non-zero-HW-class}). In case $[F:\QQ]=2,$ we provide an example in which 
the class $\beta$ is the Chern character of a coherent sheaf (Lemma \ref{lemma-box-tensor-product-of-sheaves-maps-to-one-with-equivariant-kappa-class-case-of-degree-4-extension}).

%*******************************************************************************************************************
%
%*******************************************************************************************************************
\subsection{An example of a $K$-secant $B$ spanned by Hodge classes}
\label{sec-example-of-K-secant-for-CM-field}
The subspace $\wedge_F^kH^1(X,\QQ)$ consists of elements $\gamma$ of $\wedge_Q^kH^1(X,\QQ)\subset H^1(X,\QQ)^{\otimes k}$, such that the endomorphism 
 $1 \otimes \cdots \otimes \hat{\eta}_q \otimes \cdots \otimes 1$, with $\hat{\eta}_q$ in the $j$-th tensor factor, maps $\gamma$ to the same element $\hat{\eta}_q(\gamma)$, regardless of the choice of $j$. Wedge products over $\QQ$ of such classes need not be such.
 %If $\alpha$ and $\beta$ are elements of $H^1(X)_{\hat{\sigma}}\cap H^1(X,\hat{\sigma}(F))$ and $F$ is Galois over $\QQ$, and we let $\Theta$ be  
 %the sum of the Galois orbit of $\alpha\wedge\beta$, then $\Theta$ belongs to $\wedge_F^2H^1(X,\QQ)$.

\begin{lem}
\label{lemma-wedge-2-F-H-1-contains-an-ample-class}
If $X$ is simple and $\End_\QQ(X)=F$, then 
$H^{1,1}(X,\QQ)$ is contained in the subspace $\wedge^2_FH^1(X,\QQ)$ of $H^2(X,\QQ)$. 
\end{lem}

\begin{proof}
%This is seen as follows. 
We have 
$\dim_\QQ H^{1,1}(X,\QQ)=\dim_\QQ(F)$, by
\cite[Prop. 5.5.7]{BL}. In this case $NS(X)\otimes_\ZZ\QQ$ is isomorphic to $\End_\QQ(X)$, by 
\cite[Prop. 5.2.1]{BL}. Let $L_0$ be a polarization on $X$ and let $L$ be a line bundle.
Denote their hermitian forms by $H_0$ and $H$. Let $X=V/\Lambda$ and $\hat{X}=\Omega/\hat{\Lambda}$, where $\Omega:=\Hom_{\bar{\CC}}(V,\CC)$ and
$\hat{\Lambda}=\{\omega\in\Omega \ : \ Im(\omega(\lambda))\in\ZZ, \forall \lambda\in\Lambda\}$.
Let $\phi_{H_0}:V\rightarrow\Omega$ be given by $\phi_{H_0}(v)=H_0(v,\bullet)$ and define $\phi_H$ similarly. Then $\phi_H$ is the analytic representation of $\phi_L$.
Set $\theta:=Im(H)$ and $\theta_0:=Im(H_0)$. Let $f:=\phi_H^{-1}\phi_{H_0}\in\End_\QQ(X)$.
Then $\phi_{H_0}(x)=\phi_H(f(x))$ and $H_0(x,\bullet)=H(f(x),\bullet)$, for all $x\in V$. Thus,
\begin{eqnarray*}
\theta(f(x),y)&=&Im H(f(x),y)=Im H_0(x,y)=-Im H_0(y,x) = -Im H(f(y),x)
\\
&=&Im H(x,f(y))=\theta(x,f(y)).
\end{eqnarray*}
So $\theta$ belongs to $\wedge^2_F \Lambda_\QQ^*=\wedge^2_FH^1(X,\QQ)$, since every element of $F=\End_\QQ(X)$ is represented by such $f$.
%Requiring the existence of a Hodge class in $\wedge^2_FH^1(X,\QQ)$ is in analogy to the case $F=\QQ$, 
%in which the construction was carried out only for projective compact complex tori. Note, however, that $\wedge^2_FH^1(X,\QQ)$
%is spanned by Hodge classes if $\dim_FH^1(X,\QQ)=2$ and $\dim H^{1,0}_{\hat{\sigma}}(X)=\dim H^{0,1}_{\hat{\sigma}}(X)$, 
%for all $\hat{\sigma}\in\hat{\Sigma}$. 
\end{proof}

Assume that $X$ is simple and $\End_\QQ(X)=F$.
 Then the subspace $\wedge^2_FH^1(X,\QQ)$ of $H^2(X,\QQ)$ contains a Hodge class $\Theta$, by Lemma \ref{lemma-wedge-2-F-H-1-contains-an-ample-class}.
Contraction with $\Theta$ induces an $\hat{\eta}(F)$-equivariant homomorphism $\Theta\Contract: H^1(X,\QQ)^*\rightarrow H^1(X,\QQ)$.
The element $g$ of $\Spin(V_\QQ)$ corresponding to cup product with $\exp(\Theta)$ acts on $V_\QQ:=H^1(X,\QQ)\oplus H^1(X,\QQ)^*$ 
by $\rho_g(w,\xi)=(w-\Theta\Contract \xi,\xi)$. It follows that $\rho_g$ is $\hat{\eta}(F)$-equivariant
\[
f(\rho_g(w,\xi))=f(w-\Theta\Contract\xi,\xi)=(f(w)-f(\Theta\Contract\xi),\xi\circ f)=(f(w)-\Theta\Contract(\xi\circ f),\xi\circ f)=f(\rho_g(w,\xi)).
\] 
Hence, the induced action of $\rho_g$ on $\wedge^*V_\QQ$ commutes with the induced action of $\hat{\eta}(F)$, where $\hat{\eta}_f$ acts via $\wedge^k\hat{\eta}_f$ on $\wedge^kV_\QQ$, for all $f\in F$.
Note that the nilpotent element of square zero $\rho_g-id_{V_\QQ}$ depends linearly on $\Theta$ and $\rho_g=\exp(\rho_g-id_{V_\QQ})$.
We thus define the element $g_0\in\Spin(V_{\hat{\eta},K})$ to be the unique element, such that $m_{g_0}$ acts on the spin representation $\wedge^*_K[H^1(X,\QQ)\otimes_FK]$
via cup product with $\exp(\sqrt{-q}\Theta)$, where $q$ is an element of $F$ such that $K=F[\sqrt{-q}]$. 
%as a lift of the isometry of $V_\QQ\otimes_F K$ corresponding to $id_{V_\QQ\otimes_F K}+(\rho_g-id)\otimes\sqrt{-q}$, 
The isometry $\rho_{g_0}$ of $V_\QQ\otimes_F K$ corresponds to $id_{V_\QQ\otimes_F K}+(\rho_g-id)\otimes\sqrt{-q}$.
Explicitly, decompose elements of $V_\QQ\otimes_F K$ as  $(v_1,v_2\otimes\sqrt{-q})$ with $v_i\in V_\QQ$. Then 
\begin{eqnarray*}
(\rho_g-id)(v_1,v_2\otimes\sqrt{-q})&=& ((\rho_g-id)(v_1),(\rho_g-id)(v_2)\otimes\sqrt{-q})
\\
\sqrt{-q}(v_1,v_2\otimes\sqrt{-q})&=&(\hat{\eta}_{-q}(v_2),v_1\otimes\sqrt{-q}),
\\
(\sqrt{-q}\circ (\rho_g-id))(v_1,v_2\otimes\sqrt{-q})&=&(\hat{\eta}_{-q}((\rho_g-id)(v_2)),(\rho_g-id)(v_1)\otimes\sqrt{-q})
\end{eqnarray*}
Set 
\begin{eqnarray*}
\nu&:=&(\sqrt{-q}\circ (\rho_g-id)),
\\
\rho_{g_0}&:=&\exp(\nu).
\end{eqnarray*} 
In block matrix form we have:
\[
\nu=\left(
\begin{array}{cc}
0 & -\hat{\eta}_q(\rho_g-id))
\\
\rho_g-id & 0
\end{array}
\right), \ \mbox{and} \
\exp(\nu)=\left(
\begin{array}{cc}
1 & -\hat{\eta}_q(\rho_g-id))
\\
\rho_g-id & 1
\end{array}
\right).
\]
%The lift $g_0$ of $\rho_{g_0}$ is unique up to multiplication by $-\one_S$ (acting as multiplication by $-1$ on the spin representation), 
%and we choose $g_0$ so that $g_0-\one_S$ is nilpotent. 
The induced action of $\nu$ on $\wedge^*_{K}[V_\QQ\otimes_F K]$ is nilpotent, but no longer of square zero. The induced action of $\rho_{g_0}$ on $\wedge^*_{K}[V_\QQ\otimes_F K]$ is still $\exp(\nu)$. The maximal isotropic subspace $W=\rho_{g_0}(H^1(\hat{X},\QQ)\otimes_F K)$ is the graph of
$\restricted{\nu}{}:H^1(\hat{X},\QQ)\otimes_FK\rightarrow H^1(X,\QQ)\otimes_FK$ sending $\xi\in H^1(\hat{X},\QQ)$ to
$-(\Theta\rfloor\xi)\otimes\sqrt{-q}$. Clearly, $W\cap \iota(W)=(0)$ verifying Equation (\ref{eq-W-T-and-its-conjugate-are-complementary}).

%We get the induced action $m_{g_0}$ on $\wedge^*_K [H^1(X,\QQ)\otimes_F K]$. 
The pure spinor of the maximal isotropic subspace $H^1(\hat{X},\QQ)$ is spanned by $1\in H^{ev}(X,\QQ)$. 
The pure spinor of $W$ is spanned by $m_{g_0}(1)\in 
\wedge^*_K[H^1(X,\QQ)\otimes_FK]\subset H^*(X,\QQ)\otimes_F K$. The element $m_{g_0}(1)$
is given by
\begin{eqnarray}
\nonumber
m_{g_0}(1)&=&1+\Theta\otimes\sqrt{-q} -\hat{\eta}_q(\Theta\wedge_{K}\Theta)/2-\hat{\eta}_q(\Theta\wedge_{K}\Theta\wedge_{K}\Theta)\otimes\sqrt{-q}/3!+\dots
\\
\label{eq-pure-spinor-rho-g-0-of-1}
&=& \left(1-\hat{\eta}_q(\Theta\wedge_{K}\Theta)/2+\dots\right)+
\left(\Theta-\hat{\eta}_q(\Theta\wedge_{K}\Theta\wedge_{K}\Theta)/3!+\dots\right)\otimes\sqrt{-q}
\end{eqnarray}
The $j$-th wedge product $\Theta\wedge_{K}\cdots\wedge_{K} \Theta$ 
is equal to $\Theta\wedge_F\cdots\wedge_F\Theta$, where we consider $\wedge^{2j}_FH^1(X,\QQ)$
as an $F$-subspace of $\wedge^{2j}_K[H^1(X,\QQ)\otimes_FK]$ via the isomorphism 
\[
[\wedge^{2j}_FH^1(X,\QQ)]\otimes_F K\cong \wedge^{2j}_K[H^1(X,\QQ)\otimes_FK].
\]
The product $\Theta\wedge_F\cdots\wedge_F\Theta$ is the projection of $\Theta\wedge_\QQ\cdots\wedge_\QQ\Theta$ into the 
direct summand $\wedge^{2j}_FH^1(X,\QQ)$ of $\wedge^{2j}_\QQ H^1(X,\QQ)$. Write 
$\Theta=\sum_{\hat{\sigma}\in\hat{\Sigma}}\Theta_{\hat{\sigma}}$, where $\Theta_{\hat{\sigma}}$ is in
$\wedge^2H^1_{\hat{\sigma}}(X)$. Under the embedding of $\wedge^{2j}_\QQ H^1(X,\QQ)$ in $\wedge^{2j}_\RR H^1(X,\RR)$ we get the equality $\Theta\wedge_F\cdots\wedge_F\Theta=\sum_{\hat{\sigma}\in\hat{\Sigma}}\Theta_{\hat{\sigma}}\wedge_\RR\cdots\wedge_\RR\Theta_{\hat{\sigma}}$.

%********
% Hide
%********
\hide{
Note that for $\rho_{g_0}$ to act on $H^1(X\times\hat{X},\CC)$ we need to choose an embedding 
$H^*(X,\QQ)\otimes_F K\rightarrow H^1(X\times\hat{X},\CC)$ to determine the action of $\sqrt{-q}$. So whether $\rho_{g_0}$ is an automorphism of the Hodge structure would depend on this embedding. Our choice is such that $H^1(X\times\hat{X},\QQ)$ is invariant with respect to the $\rho_{g_0}$ action (??? NO ???).
This depends on equality (\ref{eq-W-T-and-its-conjugate-are-complementary}) which $\rho_{g_0}$ needs to satisfy.
We are thus able to induce an action on $H^*(X\times\hat{X},\QQ)\cong\wedge^*_\QQ V_\QQ$ and similarly on $H^*(X,\QQ)$. 
The construction is somewhat subtle. The element $\rho_{g_0}$ maps $H^1(\hat{X},\QQ)\otimes_F K$ to a maximal isotropic subspace $W$ of $V_\QQ\otimes_F K$ defined over $K$. We then use $W$ and a choice of CM-type $T$ to define the embedding of $H^1(X\times\hat{X},\QQ)\otimes_F K$ in $H^1(X\times\hat{X},\CC)$ so that $H^1(X\times\hat{X},\QQ)$ becomes a $K$ vector space. Only then we are able to lift $\rho_{g_0}$ to an element of $\Spin(V_\QQ)$ acting on $H^*(X\times\hat{X},\QQ)=\wedge^*_\QQ V_\QQ$ rather than an element of  $\Spin(V_{\hat{\eta},K})$ acting on $\wedge^*_K [V_\QQ\otimes_F K].$ 
%********
% End Hide
%********
}
%\hfill{\EndProof}
%********
% Hide
%********
\hide{
Let $g_0$ be cup product with $\exp(\sqrt{-q}\Theta)$, where $q$ is an element of $F$ such that $K=F[\sqrt{-q}]$. 
Note that cup product with $\Theta$ commutes 
with $\hat{\eta}_q$, by our choice of $\Theta$.
We have
\[
\exp(\sqrt{-q}\Theta)=
\left(
1-\hat{\eta}_q\circ\Theta^2/2+\hat{\eta}_q^2\circ\Theta^4/4!+\dots
\right)+
\sqrt{-q}\left(
\Theta-\hat{\eta}_q\circ\Theta^3/3!+\hat{\eta}_q^2\circ\Theta^5/5!+\dots
\right)
\]
Powers of $\Theta$ above denote the endomorphisms of $H^*(X,\QQ)$ of cup product with these powers and $\circ$ denotes composition of endomorphisms. 
(??? the above doesn't make sense, since the $\hat{\eta}_q$ action on $\wedge^kH^1(X,\QQ)$ does depend on the tensor factor in $H^1(X,\QQ)^{\otimes k}$ we place it. It acts only on the subspace $\wedge^k_FH^1(X,\QQ)$. On the other hand, the cup product with the exponential of a $2$-form in $\wedge^2_F(X,\QQ)$ induces an isometry of $H^1(X,\QQ)\oplus \Hom_F(H^1(X,\QQ),F)$, which induces an isometry of $V_\QQ$, which should correspond to an element of $\Spin(V_\QQ)$, which should then act on the spin representation $H^*(X,\QQ)$. Assume that $F$ is Galois over $\QQ$. An example of a class $\alpha$, which belongs to $\wedge^k_FH^1(X,\QQ)$ is the sum of a Galois orbit of classes $\alpha_{\hat{\theta}}$ in $\wedge^kH^1(X)_{\hat{\theta}}$. Then $\hat{\eta}_q(\alpha)=\sum_{\hat{\theta}\in\hat{\Sigma}}\hat{\theta}(q)\alpha_{\hat{\theta}}$. Wedge powers of $\alpha$ are not of this type. ???)
Note that for $\exp(\sqrt{-q}\Theta)$ to act on $H^{ev}(X,\CC)$ we need to choose an embedding 
$H^*(X,\QQ)\otimes_F K\rightarrow H^*(X,\CC)$ to determine the action of $\sqrt{-q}$. So whether $g_0$ is an automorphism of the Hodge structure would depend on this embedding.
%********
% End Hide
%********
}
%\end{example}

The secant space $B$ in this example is computed via the factorization (\ref{eq-tensor-factorization-of-B}). Write 
$\Theta=\sum_{\hat{\sigma}\in\hat{\Sigma}}\Theta_{\hat{\sigma}}$, where $\Theta_{\hat{\sigma}}$ is in
$\wedge^2H^1_{\hat{\sigma}}(X)$. 
$V_\CC=\oplus_{\hat{\sigma}\in\hat{\Sigma}}V_{\hat{\sigma},\RR}\otimes_\RR\CC$, where
$V_{\hat{\sigma},\RR}=H^1_{\hat{\sigma}}(X)\oplus H^1_{\hat{\sigma}}(\hat{X})$. 
Choose a CM-type $T$ and let 
$\sigma=T(\hat{\sigma})\in\Sigma$ be the embedding restricting to $F$ as $\hat{\sigma}$. 
The isometry $\rho_{g_0}$ is mapped by the homomorphism $\tilde{e}_T$ of Remark \ref{rem-diagram-of-spin-groups-with-Spin-V-otimes-K} to an isometry which
restricts to $V_{\hat{\sigma},\RR}\otimes_\RR\CC$ as $\exp(\sigma(\sqrt{-q})\Theta_{\hat{\sigma}})$. We get the pure spinor in
$S^+_{\hat{\sigma}}$ corresponding to the element $\exp(\sigma(\sqrt{-q})\Theta_{\hat{\sigma}})$. We get the complex numbers
$[T(\hat{\sigma})(\sqrt{-q})]$ and 
the product in $S^+_\CC$
\begin{equation}
\label{eq-set-of-pure-spinnors-spanning-B}
\bigwedge_{\hat{\sigma}\in \hat{\Sigma}}\exp\left([T(\hat{\sigma})(\sqrt{-q})]\Theta_{\hat{\sigma}}\right)=
\exp\left(\sum_{\hat{\sigma}\in \hat{\Sigma}}[T(\hat{\sigma})(\sqrt{-q})]\Theta_{\hat{\sigma}}
\right)
\end{equation}
is a pure spinor spanning $\ell_T$. The secant space $B\otimes_\QQ\CC$ is the span of all these pure spinors in $S^+_\CC$. Equivalently, $B\otimes_\QQ\RR$ is the product $\otimes_{\hat{\sigma}\in\hat{\Sigma}}P_{\hat{\sigma}}$ in $S^+_\RR$, where 
\[
P_{\hat{\sigma}}\otimes_\RR\CC:=\span_\CC\{\exp(\sigma(\sqrt{-q})\Theta_{\hat{\sigma}}), \exp(\bar{\sigma}(\sqrt{-q})\Theta_{\hat{\sigma}})\}
\]
is a plane in $S^+_{\hat{\sigma}}$. The equality $\tau(P_{\hat{\sigma}})=P_{\hat{\sigma}}$, where $\tau$ is the dualization involution  
given in (\ref{eq-tau}),
%as in (\ref{eq-Mukai-pairing}), 
implies the equality
\[
\tau(B)=B.
\]

We claim that $B$ is spanned by Hodge classes. It suffices to prove that 
the $2$-form $\Theta_{\hat{\sigma}}$ is a $(1,1)$-form, for every $\hat{\sigma}\in\hat{\Sigma}$. Now, the $F$ action on $H^1(X,\QQ)$ preserves the Hodge structure, by assumption.
Hence, $H^{1,1}(X,\CC)$ is $F^\times$-invariant, where the $F^\times$-action is via $\wedge^2\hat{\eta}$. The two form $\Theta$ is of type $(1,1)$, by assumption. Hence, $\Theta=\sum_{\{\hat{\sigma}_1,\hat{\sigma}_2\}\subset\hat{\Sigma}}\Theta_{\hat{\sigma}_1\otimes\hat{\sigma}_2}$, where 
the summand $\Theta_{\hat{\sigma}_1\otimes\hat{\sigma}_2}$ is in the isotypic subspace for the character $\hat{\sigma}_1\otimes\hat{\sigma}_2$
of $F^\times$. The isotypic subspace for $\hat{\sigma}_1\otimes\hat{\sigma}_2$ is invariant with respect to the complex structure $I_X$ of $X$ and so is $\Theta$. Hence, $\Theta_{\hat{\sigma}_1\otimes\hat{\sigma}_2}$ is $I_X$-invariant and so of type $(1,1)$. Our choice of $\Theta$ in $\wedge^2_FH^1(X,\CC)$ means that $\Theta_{\hat{\sigma}_1\otimes\hat{\sigma}_2}$ is non-zero, if and only if $\hat{\theta}_1=\hat{\theta}_2$, and the summand $\Theta_{\hat{\sigma}}$ in $\wedge^2H^1_{\hat{\sigma}}(X)$ is the summand 
$\Theta_{\hat{\sigma}\otimes\hat{\sigma}}$. Hence, indeed $\Theta_{\hat{\sigma}}$ is of type $(1,1)$.

\begin{lem}
\label{lemma-I-X-hat-X-belongs-to-Omega-B}
There exists a non-zero $t\in K_-$, such that 
the complex structure $I_{X\times\hat{X}}$ of $X\times\hat{X}$ belongs to the period domain $\Omega_{B,t}$ given in (\ref{eq-omega-B-t}).
Furthermore, the polarized abelian variety of Weil type $(X\times\hat{X},\eta,\Xi_t)$ is of split type (Def. \ref{def-split-type}).
\end{lem}

\begin{proof}
Set $I:=I_{X\times\hat{X}}$. $I$ belongs to $\rho(Spin(V_\RR)_B)$, since it is the image of the element $\tilde{I}$ of $\Spin(V_\RR)$, which acts in $S_\RR$ as $m(\tilde{I})=I_X$, where $I_X$ is the complex structure of $X$, and $\tilde{I}$ belongs to $\Spin(V_\RR)_B$, since $B$ is spanned by Hodge classes. 

We show next that 
the symmetric bilinear form $g_I(x,y):=(\eta_t(x),I(y))_V$, given in (\ref{eq-symmetric-bilinear-form-g_i-CM-case}), is positive definite for some $t\in K_-$. It suffices to prove it for the restriction of $g_I$ to 
$V_{\hat{\sigma},\RR}$, for each $\hat{\sigma}\in\hat{\Sigma}$. 
Now, $\eta_t$ acts on $V_{\hat{\sigma},\CC}=V_{\sigma,\CC}\oplus V_{\bar{\sigma},\CC}$ by acting on $V_{\sigma,\CC}$ via multiplication by $\sigma(t)$ and on 
$V_{\bar{\sigma},\CC}$ via multiplication by $\bar{\sigma}(t)$. The positive-definiteness of the restriction of $g_I$ to 
$V_{\hat{\sigma},\RR}$ now follows by the same argument used to prove 
\cite[Proposition 2.4.4]{markman-sixfolds}. The latter argument shows that the restriction of $g_I$ to $V_{\hat{\sigma},\RR}$ is definite, and the choice of $t$ is made so that it is positive definite. Such a choice of $t$ exists, as observed in footnote \ref{footnote-on-CM-types} in Section \ref{sec-period-domain-for-abelian-varieties-of-Weil-type-CM-case}.
%\ref{prop-polarized-abelian-variety-of-Weil-type}.

We show next that each of $V_{\sigma,\CC}^{1,0}$ and $V_{\sigma,\CC}^{0,1}$ is $d/2$-dimensional. 
The dimension of $V_{\hat{\sigma},\RR}$ is $2d$, for each $\hat{\sigma}\in\hat{\Sigma}$, and the dimension of 
$V_\sigma$ is $d$, for each $\sigma\in\Sigma$. We know that $V_{\sigma,\CC}$ is $I$-invariant, since $I$ commutes with $\eta(K)$.
The dimension of $V_{\hat{\sigma},\RR}^{1,0}$ is equal to that of $V_{\hat{\sigma},\RR}^{0,1}$, as they are complex conjugate. It remains to prove that the dimensions of $V_{\sigma,\CC}^{0,1}$ and $V_{\sigma,\CC}^{1,0}$ are equal, for all $\sigma\in\Sigma$.
The proof is identical to that of \cite[Lemma 3.2.1]{markman-sixfolds},
%\ref{lemma-3-dimensional-eigenvalues}, 
with $V_{\hat{\sigma},\RR}$ here playing the role of $V_\RR$ there,
$V_{\sigma,\CC}$ here playing the role of $W_{1,\CC}$ there,
and $V_{\bar{\sigma},\CC}$ here playing the role of $W_{2,\CC}$ there.

Finally, we show that $(X\times\hat{X},\eta,\Xi_t)$ is of split type. Choose a basis $\{y_1, \dots, y_d\}$ of $H^1(\hat{X},\QQ)\cong\Hom_F(H^1(X,\QQ),F)$ as a vector space over $F$. Recall that $\Theta$ is an element of $\wedge^2_FH^1(X,\QQ)$.  Let $\theta:H^1(\hat{X},\QQ)\rightarrow H^1(X,\QQ)$ be the $F$-linear contraction with $\Theta$. We may further assume that $\{y_1,\dots, y_{d/2}\}$ is a basis for an $F$-subspace isotropic with respect to $\Theta$.
Then $y_i(\theta(y_j))=0,$ for $1\leq i,j\leq d/2$.
The element $g_0$
of $\Spin(V_K)$ with $m_{g_0}(\bullet)=\exp(\sqrt{-q}\Theta)\cup(\bullet)$ acts on $V_{\hat{\eta}}\otimes_FK$ sending $(w,y)$ to $(w-\theta(y)\otimes \sqrt{-q},y)$ (see \cite[III.1.7]{chevalley}). 
So
\begin{eqnarray}
\nonumber
(0,y)&=&(1/2)[m_{g_0}(0,y)+\overline{m_{g_0}}(0,y)],
\\
\label{eq-m-g-0-H-1-hat-X}
m_{g_0}(0,y)&=&(-\theta(y)\otimes \sqrt{-q},y),
\\
\label{eq-bar-m-g-0-H-1-hat-X}
\overline{m_{g_0}}(0,y)&=&(\theta(y)\otimes \sqrt{-q},y),
\\
\nonumber
\eta_t(0,y)&=& (1/2)[tm_{g_0}(0,y)-t\overline{m_{g_0}}(0,y)]
\\
\nonumber
&=&(1/2)[(-\theta(y)\otimes t\sqrt{-q},y\otimes t)+(-\theta(y)\otimes t\sqrt{-q},-y\otimes t)
\\
\nonumber
&=& (-\theta(y)\otimes t\sqrt{-q},0)
],
\\
\nonumber
(\eta_t(0,y_i),(0,y_j))_{V_{\hat{\eta}}}&=&-\Theta(y_i,y_j)t\sqrt{-q}.
\end{eqnarray}

Note that $\B:=\{(0,y_1), \dots, (0,y_d)\}$ is a $K$-basis of $V_K$. 
It suffices to show that the subspace spanned by $\{(0,y_1),\dots, (0,y_{d/2})\}$ is isotropic with respect to $H_t$, by Lemma \ref{lemma-invariants-of-hermitian-forms}(\ref{lemma-item-two-characterizations-of-split-type}).
Now, 
\begin{eqnarray*}
H_t((0,y_i),(0,y_j))&=&(-t^2)((0,y_i),(0,y_j))_{V_{\hat{\eta}}}+t(\eta_t(0,y_i),(0,y_j))_{V_{\hat{\eta}}}
\\
&=& (-t^2)((0,y_i),(0,y_j))_{V_{\hat{\eta}}}-\Theta(y_i,y_j)t^2\sqrt{-q}.
\end{eqnarray*}
The first summand vanishes, since $H^1(\hat{X},\QQ)$ is an isotropic subspace of $V_{\hat{\eta}}$ and the second summand vanishes for $i,j\leq d/2$, since 
$\{y_1,\dots, y_{d/2}\}$ is a basis for an $F$-subspace isotropic with respect to $\Theta$.
\end{proof}
%*******************************************************************************************************************
%
%*******************************************************************************************************************
\subsection{A criterion for a class $\alpha_0\otimes \beta$ in $B\otimes B$ to map to a non-zero class in $HW\!(\!X\!\!\times\!\! \hat{X}\!,\!\eta\!)$}
\label{sec-example-of-a-class-c}

In this subsection we consider the special case where the element $q$ of $F$ actually belongs to $\QQ$ (and is positive as $K$ is totally complex), though $e/2:=\dim_\QQ(F)>1$. 
Choose $\sqrt{-q}$ to be the square root in the upper half plane. 
The set (\ref{eq-set-of-pure-spinnors-spanning-B}) of $2^{e/2}$ pure spinors becomes
\begin{equation}
\label{eq-set-of-pure-spinnors-spanning-B-case-q-in-Q}
\exp\left(\sqrt{-q}\sum_{\hat{\sigma}\in \hat{\Sigma}}\pm\Theta_{\hat{\sigma}}\right).
\end{equation}
Two of those, namely $\exp\left(\pm\sqrt{-q}\Theta\right)$, are defined over the quadratic imaginary number field $K_0:=\QQ(\sqrt{-q})$ contained in $K$ and 
span a secant plane $P_{K_0}$ in $S^+_\CC$ defined over $\QQ$. In 
Lemma \ref{lemma-a-criterion-for-alpha-0-tensor-beta-to-map-via-orlov-equivalence-to-a-class-inducing-non-zero-HW-class}
in section \ref{sec-Gal-K-is-cartesion-product-of-Gal-F-and-Z-2} 
we provide a criterion for a class $\alpha_0\in P_{K_0}$ and a class $\beta\in B$ so that $\alpha_0\otimes\beta$ satisfies the hypothesis of Proposition \ref{prop-kappa-class-of-image-of-secant-class-yields-a-HW-class} and hence 
maps to a non-zero class in $HW(X\!\times \hat{X},\eta)$.

We begin in Section \ref{sec-example-CM-field-with-order-4-Galois-group} with the special case where $[F:\QQ]=2$.
If $X$ is the Jacobian of a genus $4$ curve and $\Theta$ is the natural principal polarization, then we have seen examples of secant sheave $E_0$ with Chern character $\alpha_0$ in $P_{K_0}$ in \cite[Examples 8.2.3 and 8.2.4]{markman-sixfolds}.
%\ref{example-secant-sheaf-on-jacobian-of-genus-4-curve} and \ref{example-torsion-secant-sheaf-on-jacobian-of-genus-4-curve}. 
Examples of Jacobians with real multiplication can be found in \cite{ellenberg,shimada}. 
In Lemma \ref{lemma-box-tensor-product-of-sheaves-maps-to-one-with-equivariant-kappa-class-case-of-degree-4-extension} we provide an example of a coherent sheaf $E'$ on $X$, such that $ch(E_0)\otimes ch(E')$ satisfies the hypothesis of Proposition \ref{prop-kappa-class-of-image-of-secant-class-yields-a-HW-class} and hence 
maps to a non-zero class in $HW(X\!\times \hat{X},\eta)$.

%*******************************************************************************************************************
%
%*******************************************************************************************************************
\subsubsection{The case $[K:\QQ]=4$  with $\mbox{Gal}(K/\QQ)\cong\ZZ/2\ZZ\times\ZZ/2\ZZ$}
\label{sec-example-CM-field-with-order-4-Galois-group}
Following is an example of a class $c$ in $B\otimes B$ satisfying the hypothesis of Proposition \ref{prop-kappa-class-of-image-of-secant-class-yields-a-HW-class} when 
%$[F:\QQ]=2$.
%Keep the notation and assumptions of Example \ref{example-K-is-extention-of-F-by-sqrt-of-a-negative-rational-number} and 
%Assume furthermore that 
$F=\QQ(\sqrt{t})$, where $t$ is a positive integer, which is not a perfect square.
In this case $K$ is a Galois extension of $\QQ$, being the splitting field of $(x^2+q)(x^2-t)$, and 
%\footnote{For example, the jacobians $X$ of the family of genus $4$ curves $y^3=x^5-5x^3+5x+s$ depending on the parameter $s$ 
%admit an embedding of $K=\QQ(\sqrt{-3},\sqrt{5})$  in $\End_\QQ(X)$ (see \cite{tautz})} 
$Gal(K/\QQ)$ is isomorphic to $(\ZZ/2\ZZ)^2$, since it has order $4$ and exponent $2$. Let $\gamma\in Gal(K/\QQ)$ send $\sqrt{t}$ to $-\sqrt{t}$ and
$\sqrt{-q}$ to itself. Then $\gamma$ interchanges $\Theta_{\hat{\sigma}_1}$ and $\Theta_{\hat{\sigma}_2}$, since $\Theta=\Theta_{\hat{\sigma}_1}+\Theta_{\hat{\sigma}_2}$ is $\gamma$-invariant. Hence,
$\sqrt{t}(\Theta_{\hat{\sigma}_1}-\Theta_{\hat{\sigma}_2})$ is rational, 
$\sqrt{-q}(\Theta_{\hat{\sigma}_1}-\Theta_{\hat{\sigma}_2})$ is $\iota\circ\gamma$ invariant, and the pure spinors 
$\exp(\pm \sqrt{-q}(\Theta_{\hat{\sigma}_1}-\Theta_{\hat{\sigma}_2}))$ are defined over the imaginary quadratic number subfield $K_1:=\QQ(\sqrt{-tq})$ of $K$.
Let $P_{K_1}$ be the plane spanned by these two pure spinors.
The $\QQ$-linear  $Gal(K/\QQ)$-action on the second factor $K$ of $B\otimes_\QQ K$ permutes\footnote{This permutation action on 
the four $K$-basis elements of $B\otimes_\QQ K$
defines 
a new $Gal(K/\QQ)$-action on $B\otimes_\QQ K$ via $K$-linear transformations. The latter is not natural, as it changes if we rescale the basis elements by a scalar in $K$.
} 
the four pure spinors (\ref{eq-set-of-pure-spinnors-spanning-B-case-q-in-Q}). 
This is a lift of the action of $Gal(K/\QQ)$ on the set $\T_K$ of CM-types in the sense that the map $T\mapsto \ell_T$ is $Gal(K/\QQ)$-equivariant.
The set $\T_K$ consists of two $Gal(K/\QQ)$-orbits.
%The secant space $B$ is $4$-dimensional, 
The two latter pure spinors constitute the orbit of elements fixed by $\iota\circ\gamma$. 
%subspace $(B\otimes_\QQ K)^{\iota\circ\gamma}$, 
The secant plane $P_{K_0}$ 
%in Example \ref{example-K-is-extention-of-F-by-sqrt-of-a-negative-rational-number} 
is 
spanned by the $Gal(K/\QQ)$-orbit of pure spinors fixed by $\gamma$. 

Let $\sigma_i\in\Sigma_i$ be the extension of $\hat{\sigma_i}$ mapping $\sqrt{-q}$ to the upper half plane in $\CC$. The four CM-types are
\[
\left(
\begin{array}{cc}
\hat{\sigma}_1 & \hat{\sigma}_2
\\
\sigma_1 & \sigma_2
\end{array}
\right),
\ \ 
\left(
\begin{array}{cc}
\hat{\sigma}_1 & \hat{\sigma}_2
\\
\bar{\sigma}_1 & \bar{\sigma}_2
\end{array}
\right),
\ \ 
\left(
\begin{array}{cc}
\hat{\sigma}_1 & \hat{\sigma}_2
\\
\bar{\sigma}_1 & \sigma_2
\end{array}
\right),
\ \ 
\left(
\begin{array}{cc}
\hat{\sigma}_1 & \hat{\sigma}_2
\\
\sigma_1 & \bar{\sigma}_2
\end{array}
\right).
\]

Let $\alpha_0$ be a non-zero rational class in $P_{K_0}$. Let $\beta$ be a class in $B$, which does not belong to neither $P_{K_0}$ nor $P_{K_1}$.
Set $c:=\alpha_0\otimes \beta$. Write $\beta=\beta_0+\beta_1$, with $\beta_i\in P_{K_i}$. This is possible, since $B$ is the direct sum of $P_{K_0}$ and $P_{K_1}$.
%Let $\beta_0$ be a non-zero class in $P_{K_0}$. 
Assume that $\check{\phi}(\alpha_0\otimes\beta_0)_0\neq 0.$

\begin{lem}
\label{lemma-criterion-for-c-to-satisfy-hypotheses-of-Prop-case-degree-4-extension}
The class $c$ satisfies the hypothesis of Proposition \ref{prop-kappa-class-of-image-of-secant-class-yields-a-HW-class}, namely $\check{\phi}(c)_0\neq 0$ and 
$c_1$ does not belong to $KB_1$.
\end{lem}

\begin{proof}
Let $\T_K^\gamma$ be the set of CM-types fixed by $\gamma$ and define $\T_K^{\iota\circ\gamma}$ similarly. 
The set $\T_K^\gamma$ consists of the first two CM-types displayed above and $\T_K^{\iota\circ\gamma}$ consists of the last two. We see that 
$c_1=\alpha_0\otimes\beta_1$ and $\alpha_0\otimes\beta_0$ belongs to $BB_0\oplus BB_2$.
For each $\sigma\in\Sigma$ there exists a unique CM-type $T\in \T_K^\gamma$ of which $\sigma$ is a value. The same holds for $\T_K^{\iota\circ\gamma}$. 
Hence, the left hand side of each of the linear homogeneous equations of $KB_1$ in Equation (\ref{eq-explicit-equations-for-KB-1}), for each $\sigma\in\Sigma$, is a sum consisting of precisely one non-zero term $c_{(T,T')}$. 
Finally, $\check{\phi}(c)_0=\check{\phi}(c_0)_0=\check{\phi}(\alpha_0\otimes\beta_0)_0$, by Proposition \ref{prop-kappa-class-of-image-of-secant-class-yields-a-HW-class}, which is assumed not to vanish.
\end{proof}
%\EndProof
%\end{example}

\begin{question}
Let $X$ be the Jacobian of a genus $4$ curve. 
Let $\alpha_0$ be the Chern character $ch(F_1)$ of one of the $K_0$-secant sheaves in \cite[Examples 8.2.3 or 8.2.4]{markman-sixfolds}.
%\ref{example-secant-sheaf-on-jacobian-of-genus-4-curve} or \ref{example-torsion-secant-sheaf-on-jacobian-of-genus-4-curve}.
Find a coherent sheaf $F_2$ on $X$ with Chern character $\beta$ in $B$, such that $\alpha_0$ and $\beta$ satisfy the condition of Lemma 
\ref{lemma-criterion-for-c-to-satisfy-hypotheses-of-Prop-case-degree-4-extension}, and such that the semiregularity map restricts to the image of $at_{F_2}:HT^2(X)\rightarrow \Ext^2(F_2,F_2)$ as an injective map.
\end{question}

We proceed to provide an example of a coherent sheaf $F_2$ verifying the hypotheses in the above question, except possibly for the condition on $at_{F_2}$.
%\begin{rem}
Note that the $4$ dimensional space $B$ of the $2^7$-dimensional $S^+_\QQ$ in Section \ref{sec-example-CM-field-with-order-4-Galois-group} depends on $W$, which depends on the choice of the element of $\Spin(V_{\hat{\eta}}\otimes_FK)$ of cup product with $\exp(\sqrt{-q}\Theta)$. In particular, it depends on the choice of $\Theta$ is the $2$-dimensional  $H^{1,1}(X,\QQ)=\wedge^2_FH^1(X,\QQ)$. Denote it thus by $B_{\sqrt{-q}\Theta}$.
The subspace $B_{\sqrt{-q}\Theta}$ is $\tau$-invariant, where $\tau$ is the duality involution (\ref{eq-tau}),
and it is contained in the subalgebra $\langle H^{1,1}(X,\QQ)\rangle$ 
of $S^+_\QQ=H^{ev}(X,\QQ)$ generated by $H^{1,1}(X,\QQ)$. Write $\Theta=\Theta_{\hat{\sigma}_1}+\Theta_{\hat{\sigma}_2}$ as above. Then $\Theta_{\hat{\sigma}_i}^3=0$, since $H^1_{\hat{\sigma}_i}(X)$ is $4$-dimensional. Thus, the subalgebra $\langle H^{1,1}(X,\QQ)\rangle$ is $9$-dimensional with the basis 
\[
\{1,\Theta,\sqrt{t}\tilde{\Theta}, \Theta^2,\Theta\tilde{\Theta},\tilde{\Theta}^2, \Theta^3,\sqrt{t}\tilde{\Theta}^3,[pt]
\},
\]
where $\tilde{\Theta}:=\Theta_{\hat{\sigma}_1}-\Theta_{\hat{\sigma}_2}$ and the multiplication by the positive $\sqrt{t}\in\RR\subset\CC$ above is via scalar multiplication in the $\CC$-vector space structure of $H^2(X,\CC)$. 
Note that $\sqrt{t}\tilde{\Theta}$ is the scalar product\footnote{In contrast, the rational correspondence $[\hat{\eta}(\sqrt{t})]$ sends $\Theta$ to $t\Theta$, as it acts on the second cohomology as the second exterior power of its action on the first cohomology.
} 
of the element $\Theta$ of the $F$-vector space $\wedge^2_F H^1(X,\QQ)$ by the scalar $\sqrt{t}$ in $F$,
if we let the embedding $\hat{\sigma}_1$ of $F$ send $\sqrt{t}$ to the positive square root and $\hat{\sigma}_2$ send $\sqrt{t}$ to the negative square root.
The group of units $U:=\{a+bt \ : \ a^2-tb^2=\pm 1\}$ in the subring of algebraic integers in $F$ is the product of an infinite cyclic group and $\{\pm 1\}$. If we choose $X$, so that $\hat{\eta}(U)\cap\End(X)$ 
contains an automorphism $\hat{\eta}(f)$ of $X$ of infinite order, then the rays $\RR_{\geq 0}\Theta_{\hat{\sigma}_i}$, $i=1,2$, are the two limits  $\lim_{n\rightarrow \pm \infty}\RR_{\geq 0}\hat{\eta}(f^n)^*(\Theta)$. More generally, we obtains these rays as limits of rays in the ample cone, by pulling back $\Theta$ via self isogenies 
in $\hat{\eta}(F)\cap\End(X)$.
%If, furthermore, $U=\Aut(X)$, then 
These two lines form the boundary of the ample cone, since $\Theta_{\hat{\sigma}_i}^4=0$. We conclude that neither $\sqrt{t}\tilde{\Theta}$ nor $-\sqrt{t}\tilde{\Theta}$ are ample.

Set
$\alpha:=1-\frac{q}{2}\Theta^2+\frac{q^2}{4!}\Theta^4$ and $\beta:=\Theta-\frac{q}{3!}\Theta^3$. Then $\exp(\sqrt{-q}\Theta)=\alpha+\sqrt{-q}\beta$.
Define $\tilde{\alpha}$ and $\tilde{\beta}$ similarly in terms of $\tilde{\Theta}$.  Then $B_{\sqrt{-q}\Theta}$ has the basis
\[
\{\alpha,\beta,\tilde{\alpha},\sqrt{t}\tilde{\beta}
\}.
\]
Note that $\alpha-\tilde{\alpha}=-2q\Theta_{\hat{\sigma}_1}\Theta_{\hat{\sigma}_2}$ is an element of $B_{\sqrt{-q}\Theta}$ spanning a line which is independent of the choice of $\Theta$.  
%If $\Theta'\neq \pm\Theta$, then $B_{\sqrt{-q}\Theta}\cap B_{\sqrt{-q}\Theta'}=\span_\QQ\{\Theta_{\hat{\sigma}_1}\Theta_{\hat{\sigma}_2}\}$.

Given an element $f\in F$ and a class $\delta$ in the $F$-vector space $\wedge^*_FH^1(X,\QQ)$, denote by $f\cdot \delta$ the scalar multiplication of $\delta$ by $f$, while given $c\in\CC$ we denote by $c\delta$ the scalar product in $H^2(X,\CC)$. We choose $\hat{\sigma}_1$ to map $\sqrt{t}$ to the positive square root and denote $\hat{\sigma}_1(\sqrt{t})$ by $\sqrt{t}$ for short.
So $(a+b\sqrt{t})\Theta=(a+b\sqrt{t})\Theta_{\hat{\sigma}_1}+(a+b\sqrt{t})\Theta_{\hat{\sigma}_2}$,
while $(a+b\sqrt{t})\cdot\Theta=(a+b\sqrt{t})\Theta_{\hat{\sigma}_1}+(a-b\sqrt{t})\Theta_{\hat{\sigma}_2}$.
The equality
$
B_{\sqrt{-q} f\cdot\Theta}=B_{\sqrt{-qf^2}\Theta}
$
holds for all $f\in F$ satisfying $f^2\in \QQ$, since $\sqrt{-q} f\cdot\Theta=\sqrt{-qf^2}\Theta$ in $\wedge^2_F H^1(X,\QQ)\otimes_F K$. In particular,
$B_{\sqrt{-q} c\sqrt{t}\cdot\Theta}=B_{\sqrt{-qc^2t}\Theta}$, for all $c\in \QQ^\times$.

Given $f\in F$, let $M_f\in \End[\wedge^2_FH^1(X,\QQ)\oplus H^6(X,\QQ)]$ send $(x,y)$ to $(f^2\cdot x,Nm(f^2)y)$. 

\begin{lem}
\label{lemma-M-f}
The endomorphism $M_f$ maps
$B_{\sqrt{-q}\Theta}\cap[\wedge^2_FH^1(X,\QQ)\oplus H^6(X,\QQ)]$ to $B_{\sqrt{-q}f\cdot\Theta}\cap[\wedge^2_FH^1(X,\QQ)\oplus H^6(X,\QQ)].$
\end{lem}

\begin{proof}
We have the equality
\[
-\frac{q}{3!}(a\Theta^3+b\sqrt{t}\tilde{\Theta}^3)=-\frac{q}{2}\Theta_{\hat{\sigma}_1}\Theta_{\hat{\sigma}_2}[(a-b\sqrt{t})\Theta_{\hat{\sigma}_1})+(a+b\sqrt{t})\Theta_{\hat{\sigma}_2})]=-\frac{q}{2}\Theta_{\hat{\sigma}_1}\Theta_{\hat{\sigma}_2}(a-b\sqrt{t})\cdot \Theta.
\]
The intersection of $B_{\sqrt{-q}\cdot\Theta}$ with $H^2(X,\QQ)\oplus H^6(X,\QQ)$ is 
\[
\{(a+b\sqrt{t})\cdot \Theta +
(-q/2)\Theta_{\hat{\sigma}_1}\Theta_{\hat{\sigma}_2}(a-b\sqrt{t})\cdot  \Theta
\ : \ a+b\sqrt{t}\in F
\}.
\]
The intersection of $B_{\sqrt{-q}f\cdot\Theta}$ with $H^2(X,\QQ)\oplus H^6(X,\QQ)$ is 
\[
\{(a+b\sqrt{t})\cdot f\cdot \Theta +
(-q/2)Nm(f)\Theta_{\hat{\sigma}_1}\Theta_{\hat{\sigma}_2}(a-b\sqrt{t})\cdot f\cdot \Theta
\ : \ a+b\sqrt{t}\in F
\}.
\]
The intersection is thus the graph of the linear transformation
\begin{eqnarray*}
\B_{\sqrt{-q}f\cdot\Theta}:\wedge^2_FH^2(X,\QQ)&\rightarrow&  H^6(X,\QQ)
\\
(a+b\sqrt{t})\cdot f\cdot\Theta&\mapsto & (-q/2)Nm(f)\Theta_{\hat{\sigma}_1}\Theta_{\hat{\sigma}_2}(a-b\sqrt{t})\cdot  f\cdot\Theta
\end{eqnarray*}
Thus,
\begin{eqnarray*}
\B_{\sqrt{-q}f\cdot\Theta}((a+b\sqrt{t})\cdot f^2\cdot\Theta)&=&
(-q/2)Nm(f)\Theta_{\hat{\sigma}_1}\Theta_{\hat{\sigma}_2}(a-b\sqrt{t})\cdot\gamma(f)\cdot f\cdot\Theta
\\
&=&(-q/2)Nm(f^2)\Theta_{\hat{\sigma}_1}\Theta_{\hat{\sigma}_2}(a-b\sqrt{t})\cdot\Theta
\\
&=& Nm(f^2)\B_{\sqrt{-q}\cdot\Theta}((a+b\sqrt{t})\cdot\Theta)
\end{eqnarray*}
We see that $\B_{\sqrt{-q}f\cdot\Theta}\circ (f^2\cdot(\bullet))=Nm(f^2)\B_{\sqrt{-q}\Theta}$.
%************
% Hide
%************
\hide{
 and so the two are (???) equal,\footnote
{
In fact, the projective line $\PP\left(B_{\sqrt{-q}\Theta}\otimes_\QQ\CC\cap[H^2(X,\CC)\oplus H^6(X,\CC)]\right)$ is contained in the projection of the even spinorial variety to $\PP(H^2(X,\CC)\oplus H^6(X,\CC))$. Indeed, given $a$, $b\in\CC$
the projection of the pure spinor $\exp(\sqrt{-q}[a\Theta+b\sqrt{t}\tilde{\Theta}])$ to $H^2(X,\CC)\oplus H^6(X,\CC)$ is $\sqrt{-q}$ times
\[
(a\Theta+b\sqrt{t}\tilde{\Theta})-\frac{q}{2}(a^2-tb^2)\Theta_{\hat{\sigma}_1}\Theta_{\hat{\sigma}_2}(a\Theta+b\sqrt{t}\tilde{\Theta})=
(a+b\sqrt{t})\cdot\Theta-\frac{q}{2}(a^2-tb^2)\Theta_{\hat{\sigma}_1}\Theta_{\hat{\sigma}_2}(a+b\sqrt{t})\cdot\Theta
\] 
and it is (??? NOT ???) contained in the said projective line if $a^2-b^2t=1$. So the even spinorial variety contains a double cover of the said projective line.
} 
if $Nm(f)=1$. 
%************
% End Hide
%************
}
\end{proof}

\begin{lem}
\label{lemma-Bs-intersect-trivially}
$B_{\sqrt{-q}\Theta}\cap B_{\sqrt{-q}f\cdot\Theta}\cap[\wedge^2_FH^1(X,\QQ)\oplus H^6(X,\QQ)]=(0),$ for $f\in F^\times\setminus \{\pm 1\}$,
and $B_{\sqrt{-q}\Theta}= B_{-\sqrt{-q}\Theta}$.
\end{lem}

\begin{proof}
Fix $f\in F^\times$. Let  $x\in\wedge^2_FH^1(X,\QQ)$ be a non-zero element. It can be written in the form
$x=(a+b\sqrt{t})f\cdot\Theta$. Then $\B_{\sqrt{-q}\Theta}(x)=-(q/2)\Theta_{\hat{\sigma}_1}\Theta_{\hat{\sigma}_1}(a-b\sqrt{t})\gamma(f)\cdot \Theta$ and
$\B_{\sqrt{-q}f\cdot\Theta}(x)=-(q/2)\Theta_{\hat{\sigma}_1}\Theta_{\hat{\sigma}_1}Nm(f)(a-b\sqrt{t})f\cdot \Theta$. The equality $\B_{\sqrt{-q}\Theta}(x)=\B_{\sqrt{-q}f\cdot\Theta}(x)$ is thus equivalent to $\gamma(f)=Nm(f)f$, which is equivalent to $Nm(f)=1$ and $\gamma(f)=f$, which is equivalent to $f=\pm 1$.
\end{proof}

Define an $F$-vector space structure on $\langle H^{1,1}(X,\QQ)\rangle^6$ by 
\begin{eqnarray*}
f\cdot \Theta_{\hat{\sigma}_1}^2\Theta_{\hat{\sigma}_2}&=&\hat{\sigma}_1(f)\Theta_{\hat{\sigma}_1}^2\Theta_{\hat{\sigma}_2},
\\
f\cdot \Theta_{\hat{\sigma}_1}\Theta_{\hat{\sigma}_2}^2&=&\hat{\sigma}_2(f)\Theta_{\hat{\sigma}_1}\Theta_{\hat{\sigma}_2}^2,
\end{eqnarray*} 
so that $(f\cdot\Theta)^3=Nm(f)f\cdot\Theta^3.$ 

\begin{lem}
\label{lemma-a-linear-combination-belongs-to-odd-B}
A non-zero element $f_1\cdot\Theta +f_2\cdot \Theta^3$ of $\langle H^{1,1}(X,\QQ)\rangle^2\oplus\langle H^{1,1}(X,\QQ)\rangle^6$ belongs to $\cup_{f\in F\setminus\{0\}}B_{\sqrt{-q}f\cdot\Theta}$ if and only if
$f_1\neq 0$ and $\frac{f_2}{(-q/6)\gamma(f_1)}=f^2,$ for some non-zero $f\in F$, in which case it belongs to $B_{\sqrt{-q}f\cdot\Theta}$.
\end{lem}

\begin{proof}
The statement follows from the equality 
\begin{eqnarray*}
B_{\sqrt{-q}f\cdot\Theta}&=&
\{(a+b\sqrt{t})f^{-1} f\cdot\Theta -(q/6)Nm(f)\gamma((a+b\sqrt{t})f^{-1})f\cdot\Theta^3
\}
\\
&=&\{(a+b\sqrt{t})\Theta-(q/6)(a-b\sqrt{t})f^2\cdot\Theta^3
\}.
\end{eqnarray*}
\end{proof}

\begin{cor}
\label{cor-pulling-back-one-summand-results-in-an-element-of-B}
Let $f\in F$ be a non-zero element satisfying $\hat{\eta}(f)=g^*\in \End_{Hdg}(H^1(X,\QQ))$, for some $g\in \End(X)$. 
\begin{enumerate}
\item
\label{cor-item1-linear-combination-of-power-of-thetain-B}
The class $\Theta-\frac{q}{6}g^*(\Theta^3)$ belongs to $B_{\sqrt{-q}Nm(f)f\cdot\Theta}.$ 
\item
\label{cor-item2-linear-combination-of-power-of-thetain-B}
The class $g^*\Theta-\frac{q}{6}\Theta^3$ belongs to $B_{\sqrt{-q}\gamma(f^{-1})\cdot\Theta}.$ 
\item
\label{cor-item3-linear-combination-of-power-of-thetain-B}
If $g$ is an automorphism and $Nm(f)=1$, then the class $g^*\Theta -\frac{q}{6}(g^{-1})^*(\Theta^3)$ belongs to $B_{\sqrt{-q}\Theta}$.
\end{enumerate}
\end{cor}

\begin{proof}
(\ref{cor-item1-linear-combination-of-power-of-thetain-B})
We have the equality $g^*(\Theta^3)=(f^2\cdot\Theta)^3=Nm(f)^2f^2\cdot\Theta^3$. The statement thus follows from Lemma \ref{lemma-a-linear-combination-belongs-to-odd-B}.

(\ref{cor-item2-linear-combination-of-power-of-thetain-B}) We have the equality $g^*(\Theta)=f^2\cdot\Theta$. The statement thus follows from Lemma \ref{lemma-a-linear-combination-belongs-to-odd-B}.

(\ref{cor-item3-linear-combination-of-power-of-thetain-B}) $(g^2)^*\Theta-\frac{q}{6}\Theta^3$ belongs to $B_{\sqrt{-q}\gamma(f^{-2})\cdot\Theta}$,
by part (\ref{cor-item2-linear-combination-of-power-of-thetain-B}), $\gamma(f^{-2})=f^2$,  and $B_{\sqrt{-q}f^2\cdot\Theta}=B_{\sqrt{-q}g^*(\Theta)}$. The statement follows by applying $(g^{-1})^*$ to both sides of the statement $(g^2)^*\Theta-\frac{q}{6}\Theta^3\in B_{\sqrt{-q}g^*(\Theta)}$.
\end{proof}

Assume that $X$ is the Jacobian of a non-hyperelliptic genus $4$ curve with two $g^1_3$'s with real multiplication by $F$. 
Note that the Bring curve $C$ is non-hyperelliptic and has two $g^1_3$'s (see Equation (1) and the paragraph after Proposition 2.11 in \cite{braden-hoggs}).
Every real quadratic number field embeds in $\End_{Hdg}(J_C)$, since $J_C$ is isogenous to the cartesian $4$-th power of an elliptic curve and so $\End_{Hdg}(H^1(J_C,\QQ))$ contains the algebra of $4\times 4$ matrices with rational entries.
Let $\Theta$ be the principal polarization of $X$. 
Assume that $f\in F$ satisfies $Nm(f)=1$, $f^2\neq 1$, and $\hat{\eta}(f)=g^*$, where $g$ is an automorphism of $X$. 
We have seen in  \cite[Example 8.2.4]{markman-sixfolds}
%Example \ref{example-torsion-secant-sheaf-on-jacobian-of-genus-4-curve} 
that $\Theta-\frac{q}{6}\Theta^3$ is an element of $P_{K_0}$, which is $ch(E_0)$ for a coherent sheaf $E_0$. 

\begin{example}
\label{example-coherent-sheaf-on-Jacobian-of-genus-4-with-ch-a-multiple-of-beta-prime}
Following is the construction of a simple coherent sheaf $E'$ on $X$ with Chern class an integer multiple of  $\beta':=g^*\Theta-\frac{q}{6}(g^{-1})^*(\Theta^3)$.
Choose a positive integer $d$ sufficiently large, so that there exists a positive integer $N$ such that 
\[
(N/6)\left(dg^*(\Theta^3)-q(g^{-1})^*(\Theta^3)
\right)
\]
is the class of a curve $C'$ in $X$. The class $\beta:=\Theta-\frac{d}{6}\Theta^3$ is the Chern class of a simple coherent sheaf $F'$, 
by \cite[Example 8.2.4]{markman-sixfolds}.
%\ref{example-torsion-secant-sheaf-on-jacobian-of-genus-4-curve}. 
The difference $g^*ch(F')-\beta'$ is $[-dg^*(\Theta^3)+q(g^{-1})^*(\Theta^3)]/6$. Choose $N$ generic translates of $g^*(F')$ 
so that $C'$ intersects the support of each only at points where the sheaf is a line bundle over its support. Glue the direct sum of these $N$ translates to a line bundle $L$ over $C'$ by choosing isomorphisms of the fibers at each of the points of intersection. 
For a suitable degree of $L$ the resulting coherent sheaf $E'$ will satisfy $\chi(E')=0$ and so $ch(E')=N\beta'$.
\end{example}

Let $E'$ be an object in $D^b(X)$ with $ch(E')$ a non-zero integer multiple of $g^*\Theta-\frac{q}{6}(g^{-1})^*(\Theta^3)$. Note that $ch(E')^\vee=-ch(E')$. We may thus use $\Phi(E_0\boxtimes E')$ rather than $\Phi(E_0\boxtimes (E')^\vee)$ below to simplify the notation.

\begin{lem}
\label{lemma-box-tensor-product-of-sheaves-maps-to-one-with-equivariant-kappa-class-case-of-degree-4-extension}
The object $\Phi(E_0\boxtimes E')$ has non-zero rank and $\kappa(\Phi(E_0\boxtimes E'))$ remains of Hodge type under any deformation of $(X\times\hat{X},\eta)$ as an abelian variety with complex multiplication by $K$. Furthermore, the class $\kappa_2(\Phi(E_0\boxtimes E'))$ belongs to the direct sum $\langle H^{1,1}(X\times\hat{X},\QQ)\rangle^4\oplus HW(X\times\hat{X},\eta)$, but not to $\langle H^{1,1}(X\times\hat{X},\QQ)\rangle^4.$
\end{lem}

\begin{proof}
It suffices to verify that the classes $\alpha_0:=\Theta-\frac{q}{6}\Theta^3$ and $\beta':=g^*\Theta-\frac{q}{6}(g^{-1})^*(\Theta^3)$
satisfy the hypotheses of Lemma \ref{lemma-criterion-for-c-to-satisfy-hypotheses-of-Prop-case-degree-4-extension}. 
%We have $g^*\Theta=f^2\cdot \Theta$,
%\[
%g^*\Theta_{\hat{\sigma}_1}\Theta_{\hat{\sigma}_2}=Nm(f^2)\Theta_{\hat{\sigma}_1}\Theta_{\hat{\sigma}_2},
%\]
The class $ch(E')$ belongs to $B_{\sqrt{-q}\Theta}$, 
by Corollary \ref{cor-pulling-back-one-summand-results-in-an-element-of-B},
 but it belongs to $P_{K_0}$ only if $f^2= 1$, by Lemma \ref{lemma-Bs-intersect-trivially}. We assumed that $f^2\neq 1$. 
Hence, $\beta'$ does not belong to $P_{K_0}$.
The class $\beta'$ does not belong to $P_{K_1}$, since $P_{K_1}\cap [H^2(X,\QQ)\oplus H^6(X,\QQ)]=\span_\QQ\{\sqrt{t}(\tilde{\Theta}-(q/3!)\tilde{\Theta}^3)\}$
and so the summand $f^2\cdot\Theta$ in $H^2(X,\QQ)$ is linearly independent (over $\QQ$) from $\sqrt{t}\tilde{\Theta}=\sqrt{t}\cdot\Theta$.
It remains to verify that the rank of $\Phi(E_0\boxtimes E')$ is non-zero. The rank is equal to 
$\chi(E_0\stackrel{L}{\otimes}E')$. The latter is
\begin{eqnarray*}
\int_X(\Theta-\frac{q}{2}\Theta_{\hat{\sigma}_1}\Theta_{\hat{\sigma}_2}\Theta)
(f^2\cdot\Theta-\frac{q}{2}\Theta_{\hat{\sigma}_1}\Theta_{\hat{\sigma}_2}f^{-2}\cdot\Theta)&=&
\\
-\frac{q}{2}(\hat{\sigma}_1(f)^{-2}+\hat{\sigma}_2(f)^{-2}+\hat{\sigma}_1(f)^2+\hat{\sigma}_2(f)^2)\int_X (\Theta_{\hat{\sigma}_1}\Theta_{\hat{\sigma}_2})^2&=&
\\
-\frac{q}{12}(\hat{\sigma}_1(f)^{-2}+\hat{\sigma}_2(f)^{-2}+\hat{\sigma}_1(f)^2+\hat{\sigma}_2(f)^2)\int_X \Theta^4
\neq 0.
\end{eqnarray*}
\end{proof}

%************
% Hide
%************
\hide{
Assume that $X$ is a Jacobian, $\Theta$ is its principal polarization, $f\in F$ satisfies $Nm(f)=1$, $f^2\neq 1$, and $\hat{\eta}(f)$ is an automorphism $g$ of $X$. 
We have seen in Example \ref{example-torsion-secant-sheaf-on-jacobian-of-genus-4-curve} that $\Theta-\frac{q}{6}\Theta^3$ is an element of $P_{K_0}$, which is $ch(E)$ for a coherent sheaf $E$. 
%if $X$ is a Jacobian and $\Theta$ is its principal polarization. 
%If $Nm(f)=1$ and $\hat{\eta}(f)$ is an automorphism $g$ of $X$, then

\begin{lem}
The object $\Phi(E\boxtimes g^*E)$ has non-zero rank and $\kappa(\Phi(E\boxtimes g^*E))$ remains of Hodge type under any deformation of $(X\times\hat{X},\eta)$ as an abelian variety with complex multiplication by $K$. Furthermore, $\kappa_2(\Phi(E\boxtimes g^*E))$ belongs to $\R^4\oplus HW(X\times\hat{X},\eta)$, but not to $\R^4.$
\end{lem}

\begin{proof}
It suffices to verify that the classes $\alpha_0:=\Theta-\frac{q}{6}\Theta^3$ and $\beta':=g^*(\alpha_0)$
satisfy the hypotheses of Lemma \ref{lemma-criterion-for-c-to-satisfy-hypotheses-of-Prop-case-degree-4-extension}. 
We have $g^*\Theta=f^2\cdot \Theta$,
\[
g^*\Theta_{\hat{\sigma}_1}\Theta_{\hat{\sigma}_2}=Nm(f^2)\Theta_{\hat{\sigma}_1}\Theta_{\hat{\sigma}_2}=\Theta_{\hat{\sigma}_1}\Theta_{\hat{\sigma}_2},
\]
 and so $ch(g^*E)=f^2\cdot\Theta +\B_{\sqrt{-q}\Theta}(\gamma(f^2)\cdot\Theta)$ belongs (??? NO ???) to $B_{\sqrt{-q}\Theta}$, 
by Lemma \ref{lemma-M-f},
 but it belongs to $P_{K_0}$ only if $f^2=1$. We assumed that $f^2\neq 1$. 
Hence, $\beta'$ does not belong to $P_{K_0}$.
The class $\beta'$ does not belong to $P_{K_1}$, since $P_{K_1}\cap [H^2(X,\QQ)\oplus H^6(X,\QQ)]=\span_\QQ\{\sqrt{t}(\tilde{\Theta}-(q/3!)\tilde{\Theta}^3)\}$
and so the summand $f^2\cdot\Theta$ in $H^2(X,\QQ)$ is linearly independent (over $\QQ$) from $\sqrt{t}\tilde{\Theta}=\sqrt{t}\cdot\Theta$.
It remains to verify that the rank of $\Phi(E\boxtimes g^*E)$ is non-zero. The rank is equal to 
$\chi(E\stackrel{L}{\otimes}g^*E)$. The latter is
\[
\int_X(\Theta-\frac{q}{2}\Theta_{\hat{\sigma}_1}\Theta_{\hat{\sigma}_2}\Theta)
(f^2\cdot\Theta-\frac{q}{2}\Theta_{\hat{\sigma}_1}\Theta_{\hat{\sigma}_2}f^2\cdot\Theta)=
-q(\hat{\sigma}_1(f)^2+\hat{\sigma}_2(f)^2)\int_X (\Theta_{\hat{\sigma}_1}\Theta_{\hat{\sigma}_2})^2\neq 0.
\]
\end{proof}
%************
% End Hide
%************
}

%*******************************************************************************************************************
%
%*******************************************************************************************************************
\subsubsection{The case $\mbox{Gal}(\tilde{K}/\QQ)\cong \mbox{Gal}(\tilde{F}/\QQ)\times\ZZ/2\ZZ$}
\label{sec-Gal-K-is-cartesion-product-of-Gal-F-and-Z-2}
%\begin{example}
\label{sec-main-example}
Lemma \ref{lemma-criterion-for-c-to-satisfy-hypotheses-of-Prop-case-degree-4-extension} generalizes to CM-fields $K$ with maximal totally real subfield $F$, 
such that $K=F(\sqrt{-q})$, $q\in \QQ$ and $q>0$. In this case the Galois group $\mbox{Gal}(\tilde{K}/\QQ)$ of the Galois closure $\tilde{K}$ of $K$ is the cartesian product\footnote{
For more general CM-field, the restriction homomorphism $\mbox{Gal}(\tilde{K}/\QQ)\rightarrow\mbox{Gal}(\tilde{F}/\QQ)$ is surjective, and its kernel is isomorphic to $(\Integers/2\Integers)^\nu$, where $1\leq \nu\leq e/2$, by \cite[Sec. 1.1]{dodson}. The assumption that $q$ is in $\QQ$ and $q>0$ implies  that $\nu=1$ and that the short exact sequence 
$0\rightarrow \{1,\iota\}\rightarrow \mbox{Gal}(\tilde{K}/\QQ)\rightarrow \mbox{Gal}(\tilde{F}/\QQ)\rightarrow 0$
splits.
Indeed, if $\tilde{F}$ is the splitting field of a polynomial $f(x)\in \QQ[x]$, then $\tilde{F}(\sqrt{-q})$ is the splitting field of $f(x)(x^2+q)$ and is hence a Galois extension of $\QQ$, and so $\nu=1$. Furthermore, the short exact sequence splits, since (a) the subgroup $G$ of $\mbox{Gal}(\tilde{K}/\QQ)$ stabilizing $\sqrt{-q}$ has index $2$ in $\mbox{Gal}(\tilde{K}/\QQ),$ as the $\mbox{Gal}(\tilde{K}/\QQ)$-orbit of $\sqrt{-q}$ is $\pm\sqrt{-q}$, and (b) $\iota\not\in G$.
} 
of $\mbox{Gal}(\tilde{F}/\QQ)$ and the involution $\iota$ in $\mbox{Gal}(\tilde{K}/\tilde{F})$, such that $\sigma(\iota(x))=\overline{\sigma(x)}$, for every embedding $\sigma:\tilde{K}\rightarrow\CC$
\[
\mbox{Gal}(\tilde{K}/\QQ)\cong \mbox{Gal}(\tilde{F}/\QQ)\times \{1,\iota\}.
\]
Let $T_1$ be the CM-type, such that $T_1(\hat{\sigma})$ maps $\sqrt{-q}$ to the upper-half-plane in $\CC$, for all $\hat{\sigma}\in\hat{\Sigma}$.
Let $G\subset \mbox{Gal}(\tilde{K}/\QQ)$ be the subgroup stabilizing $\sqrt{-q}$. The restriction homomorphism from $G$ to $\mbox{Gal}(\tilde{F}/\QQ)$ is an isomorphism, by our assumption on $K$, and both $T_1$ and $\bar{T}_1$ are $G$-invariant. Hence, each of $\ell_{T_1}$ 
and $\ell_{\bar{T}_1}$ is defined over the subfield $\QQ(\sqrt{-q})$ of $\tilde{K}$ fixed by $G$ and
the plane $P$ spanned in $S^+_\CC$ by $\ell_{T_1}$ and $\ell_{\bar{T}_1}$ is rational. 

Let $B_i\subset S^+_\CC$ be the subspace of $B$ spanned by
$\{\ell_T \ : \ |T\cap T_1|=i\}$. So $P=B_0\oplus B_{e/2}$.
We get the direct sum decomposition
\[
B_\CC=\oplus_{i=0}^{e/2}B_i.
\]
Note that $\bar{B}_i=B_{e/2-i}$ and $B_i\oplus B_{e/2-i}$ is $\mbox{Gal}(\tilde{K}/\QQ)$-invariant, hence it is defined over $\QQ$. We have the inclusions
\begin{eqnarray*}
B_0\otimes B_i & \subset & BB_{e/2-i},
\\
B_{e/2}\otimes B_i& \subset & BB_i,
\\
(B_0\oplus B_{e/2})\otimes (B_i\oplus B_{e/2-i}) & \subset & BB_i + BB_{e/2-i}.
\end{eqnarray*}

Let $\alpha_0$ be a non-zero rational class in $P=B_0\oplus B_{e/2}$. Note that $P$ is 
spanned by $\exp(\sqrt{-q}\Theta)$ and its complex conjugate. 
%given in terms of $\Theta$ and $q$ in equation (\ref{eq-P-introduction}) (replacing $d$ by $q$).
Let $\beta$ be a rational class in $B$. Then $\beta=\sum_{i=0}^{\lceil e/4\rceil}\beta_i$,
where $\beta_i\in B_i+ B_{e/2-i}$. 

\begin{lem}
\label{lemma-a-criterion-for-alpha-0-tensor-beta-to-map-via-orlov-equivalence-to-a-class-inducing-non-zero-HW-class}
Assume that $\beta_0\neq 0$ and $\beta_1\neq 0$ and $\check{\phi}(\alpha_0\otimes\beta_0)_0\neq 0$. The class $c:=\alpha_0\otimes \beta$
satisfies the assumptions of Proposition \ref{prop-kappa-class-of-image-of-secant-class-yields-a-HW-class}. 
\end{lem}

\begin{proof}
The class $\alpha_0\otimes\beta_i$ belongs to $BB_i+BB_{e/2-i}$, 
$c_0$ is determined by $\alpha_0\otimes\beta_0$ and $c_1$ by $\alpha_0\otimes\beta_1\in BB_1+BB_{e/2-1}$. The rationality and non-vanishing of each of $\alpha_0$ and $\beta_1$ implies that $c_1\neq 0$. Indeed, if $e=4$, then $c_1=\alpha_0\otimes\beta_1$ and if $e>4$, then 
$\alpha_0=\alpha_0'+\alpha_0''$, where $\alpha_0'\in B_0$, $\alpha_0''\in B_{e/2}$ and $\bar{\alpha}_0'=\alpha_0''\neq 0$.
$\beta_1=\beta_1'+\beta_1''$, where $\beta_1'\in B_1$, $\beta_1''\in B_{e/2-1}$, and $\bar{\beta}_1'=\beta_1''\neq 0$. 
So $c_1=\alpha_0'\otimes\beta_1''+\alpha_0''\otimes\beta_1'$,
where the first summand is in $B_0\otimes B_{e/2-1}$ and the second in $B_{e/2}\otimes B_1$. Hence, $c_1\neq 0$.
Again the left hand side of each of the linear homogeneous equations of $KB_1$ in Equation (\ref{eq-explicit-equations-for-KB-1}), for each $\sigma\in\Sigma$, is a sum consisting of precisely one non-zero term $c_{(T,T')}$. Hence,
the intersection $BB_1\cap  [B_0\oplus B_{e/2}]\otimes [B_1+B_{e/2-1}]$ is $e$-dimensional and it is complementary of $KB_1$. 
We conclude that $c_1$ does not belong to $KB_1$.
Finally, $\check{\phi}(c)_0=\check{\phi}(c_0)_0=\check{\phi}(\alpha_0\otimes\beta_0)_0$, by Proposition \ref{prop-kappa-class-of-image-of-secant-class-yields-a-HW-class}, which is assumed not to vanish.
We conclude that $c$ satisfies the assumptions of Proposition \ref{prop-kappa-class-of-image-of-secant-class-yields-a-HW-class}.
\end{proof}
%\EndProof
%***********
% Hide
%***********
\hide{
Set $O_1:=\{T_1,\bar{T}_1\}$. Let $G\subset \mbox{Gal}(\tilde{K}/\QQ)$ be the subgroup stabilizing $\sqrt{-q}$. The restriction homomorphism from $G$ to $\mbox{Gal}(\tilde{F}/\QQ)$ is an isomorphism, by our assumption on $K$, and both $T_1$ and $\bar{T}_1$ are $G$-invariant. Hence, each of $\ell_{T_1}$ 
and $\ell_{\bar{T}_1}$ is defined over the subfield $\QQ(\sqrt{-q})$ of $\tilde{K}$ fixed by $G$ and
the plane spanned in $S^+_\CC$ by $\ell_{T_1}$ and $\ell_{\bar{T}_1}$ is rational. 
Let $O'_2$ be the set of CM-types 
$T_2$, such that $T_2(\hat{\sigma})(\sqrt{-q})$ is in the upper-half-plane for precisely one $\hat{\sigma}\in\hat{\Sigma}$. Let $O_2$ be the union of $O'_2$ and $\{\bar{T} \ : \ T\in O'_2\}$. Note that $O_2$ consists of $e$ elements, if $e>4$, and $e/2$ elements if $e=4$, as in the latter case $O_2=O'_2$. 
The set $O'_2$ is $G$-invariant and the set $O_2$ is $\mbox{Gal}(\tilde{K}/\QQ)$-invariant. Hence, 
the linear subspace of $S^+_\CC$ spanned by $\ell_T$, $T\in O_2$, is defined over $\QQ$.
%The case $e=4$ is Example \ref{sec-example-CM-field-with-order-4-Galois-group}. 
%Assume that $e>4$.  
If $e>4$, then 
the  cartesian product $O_1\times O_2$ consists of two subsets, each consisting of $e$ elements. One subset consists of pairs $(T,T')$ with $|T\cap T'|=1$, the other of pairs with $|T\cap T'|=e/2-1.$ 
If $e=4$, then the  cartesian product $O_1\times O_2$ consists of $e$ pairs $(T,T')$ with $|T\cap T'|=1$.
Let $\alpha_i$ be a non-zero rational class in the rational subspace spanned by $\ell_T$, $T\in O_i$.
Let $\beta_1$ be a non-zero rational class in the rational subspace spanned by $\ell_T$, $T\in O_1$. Assume that $\phi(\alpha_1\otimes\beta_1)_0\neq 0$. Set $c=\alpha_1\otimes(\alpha_2+\beta_1)$. 
Let $\langle O_1\times O_2\rangle$ be the subspace  of $B\otimes B$ spanned by $\ell_T\otimes\ell_{T'}$, $T\in O_1$ and $T'\in O_2$. 
Then $\langle O_1\times O_2\rangle$ is contained in $BB_1\oplus BB_{e/2-1}$.
The class $\alpha_1\otimes\alpha_2$ belongs to $\langle O_1\times O_2\rangle$. The class $\alpha_1\otimes\beta_1$ belongs to $BB_0\oplus BB_{e/2}$ and so determines $c_0$, but $c_1$ is independent of $\beta_1$. The rationality of each of $\alpha_1$ and $\alpha_2$ implies that $c_1\neq 0$.
Again the left hand side of each of the linear homogeneous equations of $KB_1$ in Equation (\ref{eq-explicit-equations-for-KB-1}), for each $\sigma\in\Sigma$, is a sum consisting of precisely one non-zero term $c_{(T,T')}$. Hence,
the intersection $BB_1\cap \langle O_1\times O_2\rangle$ is $e$-dimensional and it is complementary of $KB_1$. 
Hence, $c_1$ does not belong to $KB_1$.
Finally, $\phi(c)_0=\phi(c_0)_0=\phi(\alpha_1\otimes\beta_1)_0$, by Proposition \ref{prop-kappa-class-of-image-of-secant-class-yields-a-HW-class}, which is assumed not to vanish.
We conclude that $c$ satisfies the assumptions of Proposition \ref{prop-kappa-class-of-image-of-secant-class-yields-a-HW-class}.
%***********
% End Hide
%***********
}
%\end{example}

%************
% Hide
%************
\hide{
\begin{example}
Example \ref{sec-example-CM-field-with-order-4-Galois-group} can be generalized to CM-fields $K$ with maximal totally real subfield $F$, such that $F$ is Galois over $\QQ$ with $G:=\mbox{Gal}(F/\QQ)$ admitting a subgroup of index $2$, so that $\tilde{G}:=\mbox{Gal}(K/\QQ)\cong G\times \{1,\iota\}$. 
Let $h:G\rightarrow \Integers/2\Integers$ be a surjective homomorphism. Define
$\tilde{h}:K\rightarrow K$, by $\tilde{h}(g):=\iota^{h(g)}g$. Again $K=F(\sqrt{-q})$, where $q\in\QQ$ is positive. The group $G$ acts transitively on $\hat{\Sigma}$, as $F$ is a Galois extension of $\QQ$. Let $T_1$ be the CM-type, such that $T_1(\hat{\sigma})$ maps $\sqrt{-q}$ to the upper-half-plane in $\CC$.
Then $O_1:=\{T_1,\bar{T}_1\}$ is a $\tilde{G}$-orbit. Let $T_2$ be a CM-type, such that $T_2(\hat{\sigma})$ is in the upper-half-plane for precisely one $\hat{\sigma}\in\hat{\Sigma}$. The case $e=4$ is Example \ref{sec-example-CM-field-with-order-4-Galois-group}. Assume that $e>4$.  Let $O_2$ be the $\tilde{G}$-orbit of $T_2$. Note that $O_2$ consists of $e$ elements.
The diagonal action of $\tilde{G}$ on the cartesian product $O_1\times O_2$ consists of two orbits, each consisting of $e$ elements. One orbit consists of pairs $(T,T')$ with $|T\cap T'|=1$, the other of pairs with $|T\cap T'|=e/2-1.$ Let $\alpha_i$ be a non-zero rational class in the rational subspace spanned by $\ell_T$, $T\in\T_i$.
Let $\beta_1$ be a non-zero rational class in the rational subspace spanned by $\ell_T$, $T\in\T_1$. Assume that $\phi(\alpha_1\otimes\beta_1)_0\neq 0$. Set $c=\alpha_1\otimes(\alpha_2+\beta_1)$. The argument of Example \ref{sec-example-CM-field-with-order-4-Galois-group} shows that $c$ satisfies the assumptions of Proposition \ref{prop-kappa-class-of-image-of-secant-class-yields-a-HW-class}, since 
again the left hand side of each of the linear homogeneous equations of $KB_1$ in Equation (\ref{eq-explicit-equations-for-KB-1}), for each $\sigma\in\Sigma$, is a sum consisting of precisely one non-zero term $c_{(T,T')}$. 
\end{example}
%*************
% End Hide
%*************
}

%*****************************************************************
%
%*****************************************************************

{\bf Acknowledgements:}
I thank the three referees for their insightful comments and suggestions.
This work was partially supported by a grant from the Simons Foundation Travel Support for Mathematicians program (\#962242).

%****************************************************************
% Bibliography:
%****************************************************************

\end{document}